\documentclass[12pt,reqno]{amsart}
\usepackage{amssymb}
\usepackage{amsxtra}
\usepackage{mathrsfs}
\usepackage[all]{xy}
\newtheorem{theorem}{Theorem}[section]
\newtheorem{lemma}[theorem]{Lemma}
\newtheorem{prop}[theorem]{Proposition}
\newtheorem{corollary}[theorem]{Corollary}

\theoremstyle{definition}
\newtheorem{definition}{Definition}[section]
\theoremstyle{remark}
\newtheorem{example}{Example}[section]
\newtheorem{remark}{Remark}[section]
\DeclareMathOperator{\Hom}{Hom}
\DeclareMathOperator{\End}{End}
\DeclareMathOperator{\Tor}{Tor}
\DeclareMathOperator{\Ext}{Ext}
\newcommand*{\ptens}[1]{\mathop{\widehat\otimes}_{#1}}
\newcommand*{\bptens}[1]{\mathop{\widehat\bigotimes}_{#1}}
\newcommand*{\Ptens}{\mathop{\widehat\otimes}}
\newcommand*{\tens}[1]{\mathop{\otimes}_{#1}}
\newcommand*{\Tens}{\mathop{\otimes}}
\newcommand*{\pperp}{\mathop{\underline\perp}}

\newcommand*{\pbigwedge}[1]{\mathop{\widehat\bigwedge}_{#1}}
\newcommand*{\lmod}{\mbox{-}\!\mathop{\mathbf{mod}}}
\newcommand*{\rmod}{\mathop{\mathbf{mod}}\!\mbox{-}}
\newcommand*{\bimod}{\mbox{-}\!\mathop{\mathbf{mod}}\!\mbox{-}}
\newcommand*{\id}{1}
\newcommand*{\op}{\mathrm{op}}
\renewcommand*{\dh}{\mathop{\mathrm{dh}}}
\newcommand*{\db}{\mathop{\mathrm{db}}}
\newcommand*{\dg}{\mathop{\mathrm{dg}}}

\newcommand*{\injdh}{\mathop{\mathrm{inj.dh}}}
\newcommand*{\wh}{\widehat}
\newcommand*{\wt}{\widetilde}
\newcommand*{\tp}{\mathrm{top}}
\newcommand*{\la}{\langle}
\newcommand*{\ra}{\rangle}

\newcommand{\lriso}{\stackrel{\textstyle\sim}{\smash\longrightarrow%
\vphantom{\scriptscriptstyle{_1}}}}
\newcommand{\lliso}{\stackrel{\textstyle\sim}{\smash\longleftarrow%
\vphantom{\scriptscriptstyle{_1}}}}

\newcommand*{\h}{\mathbf h}
\newcommand*{\LCS}{\mathbf{LCS}}
\newcommand*{\CC}{\mathbb C}
\newcommand*{\cC}{\mathcal C}
\newcommand*{\A}{\mathscr A}
\newcommand*{\B}{\mathscr B}
\renewcommand*{\O}{\mathscr O}
\newcommand*{\E}{\mathscr E}
\newcommand*{\N}{\mathbb N}
\newcommand*{\Z}{\mathbb Z}
\newcommand*{\R}{\mathbb R}
\renewcommand*{\L}{\mathscr L}
\renewcommand*{\P}{\mathscr P}
\renewcommand*{\H}{\mathscr H}
\newcommand*{\HH}{\mathrm H}
\newcommand*{\M}{\mathscr M}
\newcommand*{\F}{\mathscr F}
\newcommand*{\G}{\mathscr G}
\newcommand*{\I}{\mathscr I}
\newcommand*{\J}{\mathscr J}
\newcommand*{\fg}{\mathfrak g}
\newcommand*{\fh}{\mathfrak h}
\newcommand*{\fF}{\mathfrak F}
\newcommand*{\m}{\mathfrak m}
\newcommand*{\f}{\mathbf f}
\DeclareMathOperator{\Der}{Der}
\DeclareMathOperator{\Ker}{Ker}
\DeclareMathOperator{\Coh}{Coh}
\DeclareMathOperator{\St}{St}
\DeclareMathOperator{\Vect}{Vect}
\newcommand*{\shIm}{\mathop{\mathscr I\!m}}
\renewcommand*{\Im}{\mathop{\mathrm{Im}}}
\newcommand*{\spn}{\mathrm{span}}
\newcommand*{\Lie}{\mathrm{Lie}}
\newcommand*{\eps}{\varepsilon}
\newcommand*{\lla}{\longleftarrow}
\newcommand*{\lar}{\leftarrow}
\newcommand*{\rar}{\rightarrow}
\newcommand*{\xla}{\xleftarrow}
\newcommand*{\xra}{\xrightarrow}
\begin{document}
\title[Completions
of universal enveloping algebras]{Stably flat
completions of universal enveloping algebras}
\author{A. Yu. Pirkovskii}
\thanks{Partially supported by the RFBR grants 03-01-06392, 02-01-00928
and 01-01-00490}
\date{}
\begin{abstract}
We study localizations (in the sense of J.~L.~Taylor \cite{T2})
of the universal enveloping algebra, $U(\fg)$, of a complex
Lie algebra $\fg$.
Specifically, let $\theta\colon U(\fg)\to H$ be a homomorphism to some
well-behaved \cite{BFGP} topological Hopf algebra $H$.
We formulate some conditions on the dual algebra, $H'$, that are
sufficient for $H$ to be stably flat \cite{NR} over $U(\fg)$
(i.e., for $\theta$ to be a localization).
As an application, we prove that the Arens-Michael envelope,
$\wh{U}(\fg)$, of $U(\fg)$ is stably flat over $U(\fg)$ provided $\fg$ admits
a positive grading. We also show that Goodman's weighted completions
\cite{Good_hreps} of $U(\fg)$ are stably flat over $U(\fg)$
for each nilpotent Lie algebra $\fg$,
and that Rashevskii's hyperenveloping algebra \cite{Rash} is stably flat
over $U(\fg)$ for arbitrary $\fg$.
Finally, Litvinov's algebra $\A(G)$ of analytic functionals
\cite{Litv_FAA,Litv_dokl,Litv} on the corresponding
connected, simply connected complex Lie group $G$ is shown to be stably
flat over $U(\fg)$ precisely when $\fg$ is solvable.
\end{abstract}
\maketitle

One of the most important problems of modern analysis is to construct
a functional calculus of several noncommuting operators. This problem
goes back to von Neumann \cite{vN_book} and has its origin in mathematical
foundations of quantum mechanics. Functions of noncommuting variables
also appear naturally in the theory of partial differential and pseudodifferential
operators and in some problems of algebra, geometry, and mathematical
physics; see, e.g., \cite{NSS} and references therein.

A possible way to define the value of a function $f$ at an $n$-tuple
$(a_1,\ldots ,a_n)$ of linear operators is provided by the so-called
ordered representation method \cite{NSS} that was introduced
by Feynman and developed by Maslov \cite{Maslov} (see also related
papers \cite{Litv_Lapl,Litv} by Litvinov).
An essential difference of this method with the single-variable case is that
the assignment $f\mapsto f(a_1,\ldots ,a_n)$ is no longer an algebra
homomorphism. Another approach to the functional calculus problem is
based on the philosophy of noncommutative geometry: we may
change the concept of function itself and replace the commutative algebra of
functions by some noncommutative algebra.

In a coordinate-free language, a tuple of noncommuting linear operators
on a Banach space $E$ is a representation of some
finitely generated associative algebra, $A$,
that can be viewed as an ``algebra of polynomial functions on a noncommutative
space''. Therefore a noncommutative analogue of the classical
(i.e., single-variable) functional calculus problem can be formulated as follows:
Is is possible to extend the given representation $A\to\mathscr B(E)$
to some larger algebra $B$ containing $A$? Depending on their properties,
such algebras $B$ can be considered as noncommutative versions of algebras
of holomorphic functions, smooth functions, continuous functions,
Borel functions, etc. Note that this problem also makes sense when $A$
is not a subalgebra of $B$; it is sufficient that a homomorphism
$A\to B$ be fixed.

The notion of spectrum plays a key role in the functional calculus problem.
To give a simple example, recall that a bounded operator $T$ on a Banach space
$E$ has a holomorphic functional calculus on an open set $U\subset\CC$
(i.e., the homomorphism $\CC[t]\to\B(E),\; t\mapsto T$ extends to
a continuous homomorphism from the algebra $\O(U)$ of holomorphic functions
to $\B(E)$) if and only if $U$ contains
the spectrum of $T$. It is therefore natural to look for a reasonable
analogue of the notion of spectrum for several (possibly noncommuting) operators.
A general approach to this problem was suggested by J.~L.~Taylor \cite{T2}.
If $\pi\colon A\to\B(E)$ is a representation of an algebra $A$ on
a Banach space $E$, then the spectrum $\sigma(\pi,A)$ is a part of
a suitably chosen set (a ``structure space'') $\Omega_A$ of
(isomorphism classes of) locally convex $A$-modules, and $F\in\Omega_A$
does not belong to $\sigma(\pi,A)$ if and only if $\Tor_n^A(F,E)=0$
for all $n\ge 0$. (Here $\Tor_n^A$ denotes the $n$th derived functor of the
projective tensor product; see \cite{X1,T1} and Section~1 below).
For example, in the case $A=\CC[t]$ one can take
$\Omega_A$ to be the set of all $1$-dimensional $A$-modules. This set
is naturally parametrized by points of the complex plane, and the Taylor
spectrum, $\sigma(\pi,A)$, coincides with the usual spectrum, $\sigma(T)$,
of the operator $T=\pi(t)$. In the same manner, if $A=\CC[t_1,\ldots ,t_n]$,
then $\sigma(\pi,A)$ is a subset of $\CC^n$. In this case, representations
of $A$ are in bijective correspondence with $n$-tuples of commuting
operators, and $\sigma(\pi,A)$ is what is now called
the {\em Taylor joint spectrum} of the $n$-tuple $(T_1,\ldots ,T_n),\;
T_i=\pi(t_i)$. In his famous papers \cite{T00}--\cite{T2}, Taylor
established a number of remarkable properties of the joint spectrum
and constructed a multivariable version of an analytic functional
calculus. For a modern treatment of this theory, see \cite{Eschm_Put}.

The definition of $\sigma(\pi,A)$ suggested by Taylor
depends not only on the image of the representation $\pi\colon A\to\mathscr B(E)$
(i.e., not only on the given $n$-tuple of operators), but also on the
algebra $A$. Therefore, if $\pi$ can be extended to a representation $\rho$
of a larger algebra $B\supset A$, then one cannot expect that
$\sigma(\pi,A)=\sigma(\rho,B)$ in general. On the other hand, the equality still holds
in many important cases (e.g., in the above-mentioned case
$A=\CC[t],\; B=\O(U)$). Therefore it seems natural to consider only those
algebras $B\supset A$ which have the property that if some representation
$\pi$ of $A$ extends to a representation $\rho$ of $B$, then
$\sigma(\pi,A)=\sigma(\rho,B)$. More generally, if $A$ is not a subalgebra
of $B$, but a homomorphism $\theta\colon A\to B$ is given, then it is natural
to require that $\sigma(\pi,A)=\theta^*(\sigma(\rho,B))$ where
$\theta^*\colon\Omega_B\to\Omega_A$ denotes the pullback along $\theta$.

Taylor \cite{T2} introduced an appropriate class of algebra homomorphisms satisfying
the above requirement and called them {\em localizations}.
Roughly speaking, a topological algebra
homomorphism $A\to B$ is a localization if it identifies
the category of topological $B$-modules with a full subcategory
of the category of topological $A$-modules, and if the homological
relations between $B$-modules do not change when the modules
are considered as $A$-modules.
Since Taylor's objective was to
construct a holomorphic functional calculus of several {\em commuting}
operators, he considered mainly the case where $A=\CC[t_1,\ldots ,t_n]$,
the polynomial algebra endowed with the finest locally convex topology.
Taylor has proved that the canonical homomorphism of
$\mathbb C[t_1,\ldots ,t_n]$ to $\O(U)$, the Fr\'echet algebra of
holomorphic functions on an open set $U\subset\mathbb C^n$,
is a localization provided $U$ is a domain of holomorphy.
He has also shown that the canonical homomorphisms of
$\mathbb C[t_1,\ldots ,t_n]$ to the algebra $C^\infty(V)$ of smooth functions
(where $V\subset\R^n$ is an open set) and to
the algebra $\mathscr E'(\R^n)$ of compactly supported distributions
are localizations.

Thus the polynomial algebra $\mathbb C[t_1,\ldots ,t_n]$ has a rich supply
of localizations. Motivated by this example, Taylor suggested a general
scheme for constructing a noncommutative functional calculus,
a scheme where the notion of localization plays a fundamental role.
The first step of this scheme is as follows. Suppose $A$ is a fixed
finitely generated algebra (the ``base algebra'') endowed with the finest
locally convex topology. The problem is to construct a sufficiently large
family of localizations  of $A$ with values in some topological algebras
having a richer structure. Having constructed such a family, one can
hope to develop a reasonable spectral theory for representations of $A$.

As was said above, Taylor defined localizations in the topological algebra
setting. In pure algebra, a notion analogous to that
of localization was introduced by W.~Geigle and H.~Lenzing \cite{GL}
under the name ``homological epimorphism''. This notion turned out to
be useful in the representation theory of finite-dimensional algebras
(see \cite{Crawley}). Recently, A.~Neeman and A.~Ranicki \cite{NR} applied
homological epimorphisms to some problems of algebraic $K$-theory.
They use a different terminology;
namely, in the case where $\theta\colon A\to B$ is a homological epimorphism,
Neeman and Ranicki say that {\em $B$ is stably flat over $A$},
while the word ``localization'' is used by them in a different
(rather ring-theoretical than homological) sense.
We adopt both the languages here and use the phrases
``$\theta\colon A\to B$ is a localization'' (in Taylor's sense)
and ``$B$ is stably flat over $A$'' as synonyms.
The reason is that the word ``localization''
is used in modern mathematics in many different senses, and the
terminology of \cite{T2} is not the most common one.
On the other hand, it is convenient to use
Taylor's terminology when it is necessary to emphasize the role of the
homomorphism~$\theta$.

Taylor \cite{T2} has pointed out that
a possible candidate for an algebra $B$ which often seems to be
stably flat over $A$ is its Arens-Michael envelope
(the completed l.m.c. envelope, in the terminology of \cite{T2}), which is
defined as the completion of $A$
w.r.t. the family of all submultiplicative seminorms on $A$.
From the viewpoint of operator theory, an important property of $\wh{A}$
(which uniquely characterizes $\wh{A}$ within the class of Arens-Michael algebras)
is that $A$ and $\wh{A}$ have the same set of continuous
Banach space representations. If $A=\CC[t_1,\ldots ,t_n]$, then
$\wh{A}$ is isomorphic to the algebra $\O(\CC^n)$ of entire functions
(and hence is stably flat over $A$).
Thus the Arens-Michael envelope of a noncommutative finitely generated
algebra can be viewed as an ``algebra of noncommutative entire functions''.

Apart from the polynomial algebra, Taylor \cite{T2,T3} has also studied
localizations of the free algebra $F_n$ on $n$ generators.
In particular, he proved that the canonical homomorphism of $F_n$ to its
Arens-Michael envelope $\wh{F}_n$ is a localization
(i.e., $\wh{F}_n$ is stably flat over $F_n$).
Some results on localizations of $F_n$ were also obtained
by Luminet \cite{Luminet_PI}.

Another important class of noncommutative algebras considered by Taylor
is that of universal enveloping algebras. Let $\fg$ be a complex Lie algebra
and $U(\fg)$ its universal enveloping algebra.
Taylor \cite{T2} proved that if $\fg$ is semisimple, then
the Arens-Michael envelope $\wh{U}(\fg)$ of $U(\fg)$
fails to be stably flat over $U(\fg)$, in contrast to the abelian case.
On the other hand, Dosiev \cite{Dos_hprop} has recently
proved that $\wh{U}(\fg)$ is stably flat over $U(\fg)$ provided
$\fg$ is metabelian (i.e., $[\fg,[\fg,\fg]]=0$).
A natural conjecture is that $\wh{U}(\fg)$ is stably flat over $U(\fg)$
for each nilpotent Lie algebra $\fg$, but this question is still open.

In this paper we consider some ``standard'' locally convex algebras $H$
that contain $U(\fg)$ as a dense subalgebra, and study the question of
whether or not they are stably flat over $U(\fg)$.
Specifically, we concentrate on the following algebras:
\begin{itemize}
\item $H=\wh{U}(\fg)$, the Arens-Michael envelope of $U(\fg)$;
\item $H=U(\fg)_\M$, Goodman's weighted completion of $U(\fg)$
\cite{Good_filt,Good_hreps};
\item $H=\fF(\fg)$, Rashevskii's hyperenveloping algebra \cite{Rash};
\item $H=\A(G)$, Litvinov's algebra of analytic functionals
on the corresponding
connected, simply connected complex Lie group $G$
\cite{Litv_FAA,Litv_dokl,Litv}.
\end{itemize}
We generalize the above-mentioned result of Dosiev
and show that $\wh{U}(\fg)$ is stably flat over $U(\fg)$ provided $\fg$ admits
a positive grading. The weighted completion $U(\fg)_\M$ is shown to be
stably flat over $U(\fg)$ for each nilpotent Lie algebra $\fg$
and each entire weight sequence $\M$ (for terminology, see \cite{Good_hreps}
and Section~\ref{sect:weight} below). We also prove that
Rashevskii's hyperenveloping algebra $\fF(\fg)$ is stably flat over $U(\fg)$
for every Lie algebra $\fg$. Finally, $\A(G)$ turns out to be stably flat
over $U(\fg)$ if and only if $\fg$ is solvable.

A common feature of the above algebras $H\supset U(\fg)$
is that they are {\em well-behaved
topological Hopf algebras} \cite{BFGP,Litv_Hopf}. This means that they are Hopf algebras
in the tensor category of complete locally convex spaces equipped with the
projective tensor product $\Ptens$, and that their underlying locally convex
spaces are either nuclear Fr\'echet spaces or nuclear (DF)-spaces.
The category of well-behaved topological Hopf algebras has a number
of remarkable properties (see \cite{BFGP,Litv_Hopf}); in particular, it is
anti-equivalent to itself via the strong duality functor.
To answer the question of whether or not a morphism $U(\fg)\to H$
in this category is a localization, we propose a general method
that applies to all of the above-mentioned algebras $H$.
This method is based on the following observation.
Let $\fg$ be a Lie algebra, and let $V_\cdot(\fg)=C_\cdot(\fg,U(\fg))$
denote the Koszul resolution of the trivial $\fg$-module $\CC$.
The classical fact that the augmented complex $0 \lar \CC \lar V_\cdot(\fg)$
is exact is traditionally proved by introducing an appropriate
filtration on this complex and then using an induction
or a spectral sequence argument (see, e.g.,
\cite{CE}, Chap.~XIII,~Theorem 7.1, or \cite{Guichardet}, Chap.~II, Lemme~2.2).
However, if $\fg$ a finite-dimensional Lie algebra over $\CC$, it is
possible to give another proof using the fact that the (topological)
dual of $U(\fg)$ is isomorphic to the Fr\'echet algebra $\CC[[z_1,\ldots ,z_n]]$
of formal power series. The main point is that the complex dual
to $V_\cdot(\fg)$ turns out to be isomorphic to the (formal)
de Rham complex of $\CC[[z_1,\ldots ,z_n]]$.
By the Poincar\'e lemma, the latter complex (augmented by the unit map
$\CC\to\CC[[z_1,\ldots ,z_n]]$) splits as a complex of topological vector spaces.
Taking the topological dual, we conclude that $0 \lar \CC \lar V_\cdot(\fg)$
is exact.

The advantage of this proof is that it carries over to topological
algebras more general than $U(\fg)$. This suggests the following
approach to the above-mentioned localization problem for $U(\fg)$.
Given a well-behaved topological Hopf algebra $H$ and a morphism
$\theta\colon U(\fg)\to H$, we can view $H$ as a right $\fg$-module
via $\theta$. Using a version of the Cartan-Eilenberg ``inverse process''
(see \cite{CE}, Chap.~X), we prove that $\theta$ is a localization
if and only if the standard complex $C_\cdot(\fg,H)$ augmented by the
counit map $H\to\CC$ splits as a complex of topological vector spaces.
Due to the reflexivity of the algebras involved, this happens
precisely when the dual complex $0\to\CC\to C^\cdot(\fg,H')$
splits. Suppose now that $H$ is cocommutative; then the dual algebra, $H'$,
is commutative. Under some additional conditions on $H$, the latter
complex turns out to be isomorphic to the de Rham complex of $A=H'$.
In this situation we say that $A$ is {\em $\fg$-parallelizable}.
Thus the problem of whether or not $\theta$ is a localization
reduces to the question of whether or not the augmented de Rham complex
$0\to\CC\to\Omega(A)$ splits. A sufficient condition for this to be true is
that $A$ be {\em contractible} in the sense of Chen \cite{Chen_hmt}.
Therefore in order to prove that $\theta\colon U(\fg)\to H$
is a localization it is sufficient to show that $H'$
is $\fg$-parallelizable and contractible.

It should be noted that the above method is inspired by the following
result due to Taylor \cite{T2}.
Suppose that $\fg$ is the complexification of the Lie algebra of a real
Lie group $G$. Then $U(\fg)$ is canonically embedded into $\mathscr E'(G)$,
the algebra of compactly supported distributions on $G$. Taylor proved that
this embedding is a localization if and only if the de Rham cohomology
of $G$ vanishes.
The method described above is in fact a generalization of Taylor's proof.

This paper is organized as follows. In Section~1 we recall some basic facts
from topological homology (i.e., the homology theory for locally
convex algebras \cite{X1}).
We also discuss ``continuous versions'' of some concepts from pure algebra
such as DG~algebras, K\"ahler differentials and de~Rham cohomology.
Section~2 is devoted to a version of the Cartan-Eilenberg inverse process
for topological Hopf algebras. As a byproduct, we
describe the Hochschild cohomology groups of the algebras $\ell^1(G)$
(where $G$ is a discrete group) and $\E'(G)$ (where $G$ is a real Lie group)
in terms of the bounded and continuous cohomology groups of $G$.
As another application, we show that a Banach Hopf algebra with invertible antipode
is amenable precisely when it is left amenable in the sense of Lau \cite{Lau_F-alg}.
In Section~3 we discuss the notion of localization for topological algebras
and introduce related concepts of weak localization and strong transversality.
The latter notion is a somewhat stronger version of transversality condition
for Fr\'echet modules that was introduced in \cite{Ast}
and has proved to be extremely useful in complex analytic geometry
and operator theory \cite{KV,Ast,Eschm_Put,Demailly}.
Using results of the previous section, we show that for Hopf $\Ptens$-algebras
with invertible antipode the notions of localization and weak localization coincide.
In Section~4 we recall some portions of Chen's algebraic homotopy theory
\cite{Chen_hmt} in the topological algebra framework, and apply
this theory to localizations of $U(\fg)$
within the category of well-behaved cocommutative
Hopf $\Ptens$-algebras.
Given a morphism $\theta\colon U(\fg)\to H$ in this category
such that $\Im\theta$ is dense in $H$, we
show that $\theta$ is a localization provided $H'$ is $\fg$-parallelizable
and contractible. In Section~5 we concentrate on nilpotent Lie algebras
and show that the dual of a well-behaved Hopf $\Ptens$-algebra $H$
containing $U(\fg)$ is $\fg$-parallelizable provided $H$ is contained
in the formal power series completion $[U(\fg)]$ of $U(\fg)$.
In Section~6 we discuss some general properties of Arens-Michael envelopes
and describe the Arens-Michael envelopes of graded algebras as certain
vector-valued K\"othe spaces. As a corollary, we show that the dual of
the Arens-Michael envelope of $U(\fg)$ is $\fg$-parallelizable provided
$\fg$ admits a positive grading. Next we introduce a notion of contractible
Lie algebra. By definition, a Lie algebra $\fg$ is contractible if there is
a smooth path in the set of all endomorphisms of $\fg$ connecting the zero
endomorphism and the identity endomorphism of $\fg$. We show that if $\fg$
is contractible, then $\wh{U}'(\fg)$, the dual of the Arens-Michael envelope
of $U(\fg)$, is contractible in the sense of Chen. This result is then used to prove that
the Arens-Michael envelope of $U(\fg)$ is stably flat over $U(\fg)$
for each positively graded $\fg$. As a byproduct, we show that the
injective homological dimension of each nonzero $\wh{U}(\fg)$-$\Ptens$-module
is equal to the dimension of $\fg$. In Sections~7 and 8 we prove the
above-mentioned results on the stable flatness of weighted completions of $U(\fg)$,
hyperenveloping algebras, and algebras of analytic
functionals. Finally, in Section~9 we explain how the completions
of $U(\fg)$ considered above are related to one another, and formulate
some open problems.

\medskip\noindent
\textbf{Remark.} A. Dosiev has kindly informed the author that he
proved the stable flatness of the Arens-Michael envelope $\wh{U}(\fg)$ over
$U(\fg)$ under the condition that $\fg$ is a nilpotent Lie algebra
with {\em normal growth}. Roughly speaking, $\fg$ has normal growth
if for each embedding of $\fg$ into a Banach algebra norms
of powers of elements from $[\fg,\fg]$ decrease sufficiently rapidly.
The class of Lie algebras with normal growth contains all metabelian Lie
algebras, but it is not clear how this class is related to that
of positively graded Lie algebras.

\medskip\noindent\textbf{Acknowledgments. }
The author is grateful to
A.~Ya.~Helemskii and A.~Dosiev for valuable discussions,
and to G.~L.~Litvinov for helpful comments.

\section{Preliminaries}
\label{sect:prelim}
We shall work over the complex numbers $\CC$. All associative algebras
are assumed to be unital.

\subsection{Topological algebras and modules}
\label{subsect:topalg}
In this subsection we recall some basic notions from topological homology
(homology theory for topological algebras). For more details, see
\cite{X1}, \cite{X_HOA}, and \cite{T1}.

We refer to \cite{Sch} and \cite{Treves} for general facts on topological
vector spaces. Given topological vector spaces $E$ and $F$, we denote
by $\L(E,F)$ the space of all linear continuous maps from $E$ to $F$.
We endow $\L(E,F)$ with the topology of uniform convergence on
bounded subsets of $E$. Unless otherwise specified, $E'=\L(E,\CC)$
denotes the strong dual of $E$.
The completion of $E$ is denoted by $E^\sim$. If $E$ and $F$ are locally
convex spaces (l.c.s.'s), then $E\Ptens F$ stands for their completed projective
tensor product.

By a topological algebra we mean a topological vector space $A$
together with the structure of associative algebra such that
the multiplication map $A\times A\to A$ is separately continuous.
A complete, Hausdorff, locally convex topological algebra with
jointly continuous multiplication is called a
\textit{$\Ptens$-algebra} (see \cite{T1,X1}). If $A$ is a $\Ptens$-algebra, then
the multiplication $A\times A\to A$ extends
to a linear continuous map from the projective tensor
product $A\Ptens A$ to $A$. In other words, a
$\Ptens$-algebra is just an
algebra in the tensor category $(\mathbf{LCS},\Ptens)$ of complete l.c.s.'s
(cf. Section~\ref{sect:inverse} below).

Recall that a seminorm $\|\cdot\|$ on an algebra $A$ is called
{\em submultiplicative} if $\| ab\|\le\| a\|\| b\|$ for all $a,b\in A$.
A $\Ptens$-algebra $A$ is called
an \textit{Arens-Michael algebra} (or a \textit{locally $m$-convex} algebra)
if its topology can be defined by a family of submultiplicative seminorms
(see \cite{Michael,X2}).
In particular, any Banach algebra is an Arens-Michael algebra.
By a {\em Fr\'echet algebra} we mean a metrizable (not necessarily
locally $m$-convex) $\Ptens$-algebra.

Each associative $\CC$-algebra $A$ becomes a topological algebra
w.r.t. the finest locally convex topology. We denote the resulting
topological algebra by $A_s$. If $A$ has countable dimension as a vector
space, then $A_s$ is a $\Ptens$-algebra \cite{BFGP}. In particular,
this condition is satisfied whenever $A$ is finitely generated.

A left \textit{$\Ptens$-module} over
a $\Ptens$-algebra $A$ (a left $A$-$\Ptens$-module for short)
is a complete Hausdorff locally convex space $X$ together
with the structure of left unital $A$-module such that the map
$A\times X\to X,\; (a,x)\mapsto a\cdot x$ is jointly continuous.
As above, this means that $X$ is a left $A$-module in $(\mathbf{LCS},\Ptens)$.
Given two left $A$-$\Ptens$-modules $X$ and $Y$, an \textit{$A$-module morphism}
is a linear continuous map $\varphi\colon X\to Y$ such that
$\varphi(a\cdot x)=a\cdot\varphi(x)$ for all $a\in A,\; x\in X$. The vector space
of all $A$-module morphisms from $X$ to $Y$ is denoted by ${_A}\h(X,Y)$.

Right $A$-$\Ptens$-modules, $A$-$\Ptens$-bimodules, and their morphisms
are defined similarly.
As in pure algebra, $A$-$\Ptens$-bimodules can be regarded as either left
or right $\Ptens$-modules over the algebra
$A^e=A\Ptens A^{\mathrm{op}}$, where $A^{\mathrm{op}}$ stands for the
algebra opposite to $A$.
Given two right $A$-$\Ptens$-modules
(respectively, $A$-$\Ptens$-bimodules) $X$ and $Y$, we use the notation $\h_A(X,Y)$
(respectively, ${_A}\h_A(X,Y)$) to denote the corresponding space of morphisms.
The resulting module categories are denoted by $A\lmod$, $\rmod A$,
and $A\bimod A$, respectively.

If $\theta\colon A\to B$ is a $\Ptens$-algebra homomorphism
(i.e., a unital continuous homomorphism), then each left (resp.
right) $B$-$\Ptens$-module $X$ can be considered as a left (resp. right)
$A$-$\Ptens$-module via $\theta$. Sometimes we will denote the
resulting $A$-$\Ptens$-module by ${_\theta}X$ (resp. $X_\theta$).

If $A$ is a commutative $\Ptens$-algebra, then an $A$-$\Ptens$-bimodule
$X$ is {\em symmetric} if $a\cdot x=x\cdot a$ for all $a\in A,\; x\in X$.
As usual, we identify left modules, right modules,
and symmetric bimodules over a commutative algebra and call them
just ``modules''.

Let $A$ be a $\Ptens$-algebra and $M$ an $A$-$\Ptens$-bimodule.
Recall that a linear continuous map $D\colon A\to M$ is a {\em derivation}
if $D(ab)=Da\cdot b+a\cdot Db$ for all $a,b\in A$.
Denote by $\Der(A,M)$ the set of all continuous derivations from $A$ to $M$.
We also set $\Der A=\Der(A,A)$.
If $A$ is commutative, we may speak about derivations of $A$ with
coefficients in left $A$-$\Ptens$-modules by identifying left modules
with symmetric bimodules (see above).

If $X$ is a right $A$-$\Ptens$-module and $Y$
is a left $A$-$\Ptens$-module, then their $A$-module tensor product
$X\ptens{A}Y$ is defined to be
the completion of the quotient $(X\Ptens Y)/N$, where $N\subset X\Ptens Y$
is the closed linear span of all elements of the form
$x\cdot a\otimes y-x\otimes a\cdot y$
($x\in X$, $y\in Y$, $a\in A$)\footnote[1]{Here we follow
Helemskii's monograph \cite{X1}. To avoid confusion, we note that
this definition of $X\ptens{A} Y$ is different from that given by
Kiehl and Verdier \cite{KV} and Taylor \cite{T1} (and used also in \cite{Ast}
and \cite{Eschm_Put}). More precisely, $X\ptens{A} Y$ is the completion
of the Kiehl-Verdier-Taylor tensor product.}.
As in pure algebra, the $A$-module tensor product can be characterized
by a certain universal property (see \cite{X1} for details).

A morphism $\sigma\colon X\to Y$ of left $A$-$\Ptens$-modules
is said to be an {\em admissible epimorphism} if there exists
a linear continuous map $\tau\colon Y\to X$ such that
$\sigma\tau=\mathbf 1_Y$, i.e., if $\sigma$ is a retraction
when considered in the category of topological vector spaces.
Similarly, a morphism
$\varkappa\colon X\to Y$ is an {\em admissible monomorphism}
if it is a coretraction in the category of topological vector spaces.
Finally, a morphism $\varphi\colon X\to Y$ is
{\em admissible} if it possesses a factorization $\varphi=\varkappa\sigma$
with $\varkappa$ an admissible monomorphism and $\sigma$ an admissible epimorphism.
Geometrically, this means that the kernel and the image of $\varphi$
are complemented subspaces of $X$ and $Y$, respectively, and $\varphi$
is an open map of $X$ onto its image.
A chain complex $X_\bullet=(X_n,d_n)$ of left $A$-$\Ptens$-modules
is called {\em admissible} if it splits as a complex
of topological vector spaces. Equivalently, $X_\bullet$ is admissible if it is
exact, and all the $d_n$'s are admissible morphisms.

\begin{remark}
It can easily be checked that the category $A\lmod$ together with the
class of admissible monomorphisms and epimorphisms satisfies the axioms of
exact category (see \cite{Quillen} and \cite{DCU}),
so that the main constructions of abstract
homological algebra (derived categories, ``total'' derived functors, etc.) make sense
in this setting. However, we shall not use such a general approach here.
\end{remark}

An $A$-module $P\in A\lmod$ is called \textit{projective} if for each
admissible epimorphism $X\to Y$ in $A\lmod$ the induced map
${_A}\h(P,X)\to {_A}\h(P,Y)$ is surjective.
Dually, an $A$-module $Q\in A\lmod$ is called \textit{injective} if for each
admissible monomorphism $X\to Y$ in $A\lmod$ the induced map
${_A}\h(Y,Q)\to {_A}\h(X,Q)$ is surjective.
For each $E\in\mathbf{LCS}$ the projective tensor product $F=A\Ptens E$ has
a natural structure of left $A$-$\Ptens$-module with operation
defined by $a\cdot (b\otimes x)=ab\otimes x$. Such modules are called {\em free}.
In view of the natural isomorphism ${_A}\h(A\Ptens E,Y)\cong\L(E,Y),\;
Y\in A\lmod$, each free module is projective. This implies that
the category $A\lmod$ has {\em enough projectives}, i.e.,
for each $X\in A\lmod$ there exists an admissible epimorphism $P\to X$ with
$P$ projective. To see this, it suffices to set $P=A\Ptens X$ and
to define $A\Ptens X\to X$ by $a\otimes x\mapsto a\cdot x$.

\begin{remark}
\label{rem:noinj}
If $A$ is a Banach algebra, then the category $A\lmod$ has {\em enough
injectives} as well, i.e., each $X\in A\lmod$ can be embedded into an injective
$A$-$\Ptens$-module via an admissible monomorphism \cite{X1}. However,
if $A$ is non-normable, then $A\lmod$ may fail to possess
nonzero injective objects \cite{Pir_inj2,Pir_Nova}; see also
Corollary~\ref{cor:injdim} below.
\end{remark}

Given a left $A$-$\Ptens$-module $X$, a \textit{resolution} of $X$
is a chain complex $P_\bullet=(P_n,d_n)_{n\ge 0}$ of left $A$-$\Ptens$-modules
together with a morphism $\epsilon\colon P_0\to X$
such that the augmented sequence
$$
0\xla{} X\xla{\epsilon} P_0\xla{d_0}\cdots \xla{} P_n\xla{d_n} P_{n+1}\xla{}\cdots
$$
is admissible. A {\em projective resolution} is a resolution
in which all the $P_i$'s are projective $\Ptens$-modules.
Since $A\lmod$ has enough projectives, it follows that
each $A$-$\Ptens$-module has a projective resolution.
Therefore one can define
the derived functors $\Ext$ and $\Tor$ following the general patterns
of relative homological algebra (see \cite{X1}).
Namely, take a projective resolution
$P_\bullet$ of an $A$-module $X\in A\lmod$ and set
$$
\Ext^n_A(X,Y)=\HH^n\bigl( {_A}\h(P_\bullet,Y)\bigr)
$$
for each $Y\in A\lmod$.
(Here ``$\HH^n$'' stands for the $n$th cohomology space).
Similarly,
$$
\Tor_n^A(Y,X)=\HH_n(Y\ptens{A}P_\bullet)
$$
for each $Y\in\rmod A$. Of course,
$\Ext^n_A(X,Y)$ and $\Tor_n^A(Y,X)$ do not depend on the particular
choice of the resolution $P_\bullet$ and possess the usual functorial
properties (see \cite{X1} for details).

A {\em projective bimodule resolution} of $A$ is a projective resolution
of $A$ in $A\bimod A$.

If $M\in A\bimod A$, then the {\em $n$th
Hochschild cohomology} (resp. \emph{homology}) of $A$ with coefficients in $M$
is defined as $\H^n(A,M)=\Ext^n_{A^e}(A,M)$
(resp. $\H_n(A,M)=\Tor_n^{A^e}(M,A)$).

A left $A$-$\Ptens$-module $X$ has {\em projective homological dimension $\le n$}
if $X$ has a projective resolution $P_\bullet$ such that
$P_i=0$ for all $i>n$. The least such integer $n$ is denoted by
$\dh_A X$ and is called the {\em projective homological dimension} of $X$.
Equivalently,
$$
\dh\nolimits_A X=\min\{\, n : \Ext^{n+1}_A(X,Y)=0\;\forall\, Y\in A\lmod\,\}.
$$
If no such $n$ exists, one sets $\dh_A X=\infty$.

The {\em injective homological dimension} of $X\in A\lmod$ is
$$
\injdh\nolimits_A X=\min\{\, n : \Ext^{n+1}_A(Y,X)=0\;\forall\, Y\in A\lmod\,\}.
$$
If $A$ is a Banach algebra, then $\injdh_A X$ can also be defined
as the length of the shortest injective resolution of $X$
(cf. Remark~\ref{rem:noinj} above).

An $A$-$\Ptens$-module $X$ is projective (resp. injective) if and only if
${\dh_A X=0}$ (resp. $\injdh_A X=0$).

The {\em left global dimension} of $A$ is
$$
\dg A=\sup\{\,\dh\nolimits_A X : X\in A\lmod\,\}.
$$

Similarly, one can define homological dimension for right
$A$-$\Ptens$-modules and for $A$-$\Ptens$-bimodules.
The homological dimension of $A$ considered as an $A$-$\Ptens$-bimodule
is called the {\em homological bidimension} of $A$ and is denoted by
$\db A$. For every $\Ptens$-algebra $A$ we have
$\dg A\le\db A$.

\subsection{K\"ahler differentials}
Recall some facts about K\"ahler differentials and de Rham
cohomology for commutative $\Ptens$-algebras. Most of this material is well-known
in the purely algebraic case (see, e.g., \cite{EGA} or \cite{Kunz_new}).
For the $\Ptens$-case, see \cite{Pflaum}, \cite{Mal_DR}.

Let $A$ be a commutative $\Ptens$-algebra.
A pair $(\Omega^1 A,d_A)$ consisting of an $A$-$\Ptens$-module
$\Omega^1 A$ and
a derivation $d_A\colon A\to\Omega^1 A$ is called the {\em module of
K\"ahler differentials} if for each $A$-$\Ptens$-module $M$ and for each derivation
$D\colon A\to M$ there exists a unique $A$-$\Ptens$-module morphism
$\varphi\colon\Omega^1 A\to M$ such that the following diagram is commutative:
$$
\xymatrix@R-20pt{
& \Omega^1 A \ar@{-->}[dd]_{\varphi} \\
A \ar[ur]^{d_A} \ar[dr]_D \\
& M
}
$$
The derivation $d_A$ is called the {\em universal derivation} of $A$.

Obviously, there is a natural isomorphism
${_A}\h(\Omega^1 A,M)\cong\Der(A,M)$ defined by the rule
$\varphi\mapsto\varphi d_A$. In other words, $\Omega^1 A$ represents the functor
$M\mapsto\Der(A,M)$. Hence $\Omega^1 A$ is unique up to a
$\Ptens$-module isomorphism.

The module of K\"ahler differentials can be constructed explicitly as follows.
Let $I$ be the kernel of the product map $A\Ptens A\to A$.
Set $\Omega^1 A=(I/\overline{I^2})^\sim$ and define $d_A\colon A\to\Omega^1 A$
by $d_A a=(a\otimes 1-1\otimes a)+\overline{I^2}$. Then
$(\Omega^1 A,d_A)$ is the module of K\"ahler differentials for $A$
(see, e.g., \cite[Appendix B]{Pflaum} or \cite[\S20]{EGA} for the algebraic case).

If $A=C^\infty(M)$ is the Fr\'echet algebra of smooth functions on a manifold
$M$, then $\Omega^1 A$ is canonically isomorphic with the module
of differential $1$-forms on $M$ (cf. \cite{Pflaum}).
A similar result holds for algebras of holomorphic functions on Stein
manifolds. Since we could not find this fact in the literature, we give a complete
proof below.

Let $(V,\O_V)$ be a complex space. Consider the diagonal map
$\Delta\colon V\to V\times V$, and denote by
$\I\subset\O_{V\times V}$ the ideal sheaf of the subspace
$\Delta(V)\subset V\times V$.
Recall (\cite{Groth_Cart}; cf. also \cite{Hartshorne})
that the {\em sheaf of $1$-differentials}
of $V$ is defined as $\Omega^1_V=\Delta^*(\I/\I^2)$.
There is a canonical morphism of sheaves $d\colon\O_V\to\Omega^1_V$
defined locally as $da=(a\otimes 1-1\otimes a)+\I_{x,x}^2$
for each $a\in\O_{V,x}$, $x\in V$. If $V$ is a complex manifold,
then $\Omega^1_V$ coincides with the usual cotangent sheaf of $V$,
the space of global sections $\Omega^1(V)=\Gamma(V,\Omega^1_V)$
is the space of holomorphic differential $1$-forms on $V$, and the map
$d_V\colon \O(V)\to\Omega^1(V)$ induced by $d$ is precisely the
exterior (de Rham) derivative.

We need some facts on Stein modules \cite{Forster}.
Let $(V,\O_V)$ be a Stein space. For each coherent analytic sheaf
$\F$ on $V$ the space of global sections $\F(V)=\Gamma(V,\F)$ has a canonical
locally convex topology making it into a Fr\'echet $\O(V)$-module.
Modules of this form are called {\em Stein modules}. Denote by
$\Coh (V)$ the category of coherent sheaves of $\O_V$-modules
and by $\St (V)\subset\O(V)\lmod$ the category of Stein $\O(V)$-modules.
The functor of global sections $\Gamma\colon\Coh(V)\to\St(V)$
is exact (Cartan's Theorem B) and fully faithful \cite{Forster}.
Hence $\Gamma$ is an equivalence of $\Coh(V)$ and $\St(V)$.
If $\F,\G\in\Coh(V)$, and at least one of them is locally free, then
there is a canonical isomorphism
$\F(V)\ptens{\O(V)}\G(V)\cong(\F\tens{\O_V}\G)(V)$
(see \cite{Ast} or \cite[4.2.4]{Eschm_Put}).
If $F=\Gamma(V,\F)$ is a Stein module and $G\subset F$ is a closed submodule,
then $G$ is also a Stein module, i.e., $G=\Gamma(V,\G)$ for some
coherent subsheaf $\G\subset\F$. In particular, each closed ideal
$J\subset\O(V)$ has the form $J=\Gamma(V,\J)$ for some coherent sheaf of ideals
$\J\subset\O_V$. In this case, $\overline{J^2}=\Gamma(V,\J^2)$ (see \cite{Forster}).

\begin{lemma}
\label{lemma:Stein_diff}
Let $(V,\O_V)$ be a Stein space. Then the $\O(V)$-$\Ptens$-module
$\Omega^1(\O(V))$ of K\"ahler differentials is canonically isomorphic
with $\Omega^1(V)=\Gamma(V,\Omega_V^1)$. Under this identification,
$d_V\colon\O(V)\to\Omega^1(V)$ becomes a universal derivation.
\end{lemma}
\begin{proof}
Set $V^2=V\times V$, and let $I=\Gamma(V^2,\I)\subset\O(V^2)$
be the ideal of all functions vanishing on $\Delta(V)$.
Identifying $\O(V^2)$ with
$\O(V)\Ptens\O(V)$, we see that $I$ becomes
the kernel of the product map $\O(V)\Ptens\O(V)\to\O(V)$.
Consider the commutative diagram
$$
\xymatrix@R-15pt{
& \Gamma(V^2,\I) \ar[r]^(.45){\tilde q} & \Gamma(V^2,\I/\I^2) \ar@{=}[r] & \Omega^1(V)\\
\O(V) \ar[ur]^D \ar[dr]_D \\
& I \ar@{=}[uu] \ar[r]^q & I/\overline{I^2} \ar[uu]_j \ar@{=}[r] & \Omega^1(\O(V))
}
$$
Here $D(a)=a\otimes 1-1\otimes a$ for all $a\in\O(V)$,
$q$ is the quotient map, and $\tilde q$ is induced by the sheaf quotient map
$\I\to\I/\I^2$. Note also that $\Omega^1(\O(V))=I/\overline{I^2}$ since
$I$ is a Fr\'echet space.
Obviously, $d_V=\tilde q D$, and $d_{\O(V)}=qD$.
Since $V$ is a Stein space, we see that $j$ is an isomorphism
(see the remarks preceding the statement of the lemma). The rest is clear.
\end{proof}

In some important cases the module of K\"ahler differentials
is free and finitely generated. The next lemma gives a simple sufficient condition
for this.

\begin{lemma}
\label{lemma:Omega_free}
Suppose there exist $x_1,\ldots ,x_n\in A$ and $\partial_1,\ldots ,\partial_n\in\Der A$
such that the $x_i$'s generate a dense subalgebra of $A$,
and $\partial_i(x_j)=\delta_{ij}$
for each $i,j$. Then $\partial=(\partial_1,\ldots ,\partial_n)\colon A\to A^n$
is a universal derivation. In particular, $\Omega^1 A$ is isomorphic to $A^n$.
\end{lemma}
\begin{proof}
Let $D\colon A\to M$ be a derivation. Denote by $(u_1,\ldots ,u_n)$ the standard
basis in $A^n$ (i.e., $u_i=(0,\ldots ,1,\ldots ,0)$ with $1$ in the $i$th coordinate,
$0$ elsewhere). We have $\partial(x_i)=u_i$
for each $i=1,\ldots ,n$. Define an $A$-$\Ptens$-module morphism
$\varphi\colon A^n\to M$ by $\varphi(u_i)=D(x_i)$ for $i=1,\ldots ,n$.
Then $(\varphi\partial)(x_i)=D(x_i)$ for each $i$. Since $x_1,\ldots ,x_n$
generate a dense subalgebra of $A$, we conclude that $\varphi\partial=D$.
On the other hand, since $A^n$ is generated (as an $A$-module)
by $\Im\partial$, $\varphi$ is a unique $A$-module
morphism with the above property.
\end{proof}

\subsection{DG $\Ptens$-algebras and de Rham cohomology}
\label{subsect:DGA}
By a {\em graded $\Ptens$-algebra} we mean a sequence $\A=\{ A^n\}_{n\in\Z_+}$
of complete l.c.s.'s together with linear continuous mappings
$$
A^p\Ptens A^q \to A^{p+q},\quad (a,b)\mapsto ab
$$
satisfying the usual associativity conditions. In particular,
$A^0$ is a $\Ptens$-algebra, and each $A^n$ is an $A^0$-$\Ptens$-bimodule.
We will always assume that $\A$ is {\em unital}, i.e., that $A^0$ is unital and
each $A^n$ is a unital $A^0$-bimodule.
A graded $\Ptens$-algebra $\A$ is said
to be {\em graded commutative} if $ab=(-1)^{pq} ba$ for each $a\in A^p,\;
b\in A^q$. If $\A$ is graded commutative, then $A^0$ is commutative
in the usual sense,
and all the $A^0$-$\Ptens$-bimodules $A^n$ are symmetric.

Morphisms of graded $\Ptens$-algebras are defined in an obvious way.

If $\A$ is a graded $\Ptens$-algebra, then $A=\bigoplus_n A^n$ is a
topological algebra w.r.t. the locally convex direct sum topology.
If, in addition, each $A^n$ is finite-dimensional, then the
topology on $A$ is the finest locally convex topology, so that
$A$ is a $\Ptens$-algebra (see Subsection \ref{subsect:topalg}).
In this case we will often identify $\A$ and $A$ and say
that {\em $A=\bigoplus_n A^n$ is a graded algebra}.

Let $A$ be a commutative (ungraded) $\Ptens$-algebra.
Given an $A$-$\Ptens$-module $M$ and $n\in\N$, we can define
the $n$th exterior power of $M$ as in the purely algebraic case.
Namely, consider the antisymmetrization map
$\mathbf a_M\colon\bptens{A}^n M\to\bptens{A}^n M$ defined by
\begin{equation}
\label{asymm}
\mathbf a_M(x_1\otimes\cdots\otimes x_n)=
\frac{1}{n!}\sum_{\sigma\in S_n}
\varepsilon(\sigma)\cdot x_{\sigma^{-1}(1)}\otimes\cdots\otimes x_{\sigma^{-1}(n)}\, .
\end{equation}
It is easy to see that $\mathbf a_M$ is an $A$-$\Ptens$-module morphism
and that $\mathbf a_M^2=\mathbf a_M$.
Hence $\Im\mathbf a_M$ is a direct
$A$-$\Ptens$-module summand of $\bptens{A}^n M$.
We set $\pbigwedge{A}^n M=\Im\mathbf a_M$ (or, equivalently,
$\pbigwedge{A}^n M=\bptens{A}^n M/\Ker\mathbf a_M$), and call the resulting
$A$-$\Ptens$-module the {\em $n$th exterior power} of $M$.
As usual, for each $x_1,\ldots ,x_n\in M$ we denote the element
$\mathbf a_M(x_1\otimes\cdots\otimes x_n)$ of $\pbigwedge{A}^n M$ by
$x_1\wedge\cdots\wedge x_n$.

For each $p,q\in\Z_+$ we have a bilinear continuous map
\begin{equation*}
\textstyle\pbigwedge{A}\nolimits^p M \times\textstyle\pbigwedge{A}\nolimits^q M \to
\textstyle\pbigwedge{A}\nolimits^{p+q} M,\quad
(x,y)\mapsto x\wedge y=\mathbf a(x\otimes y)
\end{equation*}
(here we set $\pbigwedge{A}^0 M=A$ and
use the canonical identifications $A\ptens{A}X=X\ptens{A} A=X$).
As in the purely algebraic case, the above maps make
$\pbigwedge{A} M=\{ \pbigwedge{A}^p M : p\in\Z_+\}$ into a graded
commutative $\Ptens$-algebra called {\em the exterior algebra} of $M$.

Now let $\A$ be a graded commutative $\Ptens$-algebra.
For each $n\in\N$ we have an $A^0$-$\Ptens$-module morphism
\begin{equation}
\label{exterior}
\textstyle\pbigwedge{A^0}\nolimits^n A^1\to A^n,\quad
a_1\wedge\cdots\wedge a_n\mapsto a_1\cdots a_n.
\end{equation}
$\A$ is called {\em exterior} if \eqref{exterior} is an isomorphism
for each $n\in\N$ (cf. \cite{Kunz_new}).
In other words, $\A$ is exterior if the canonical
morphism $\pbigwedge{A^0} A^1\to\A$ is a graded $\Ptens$-algebra isomorphism.
It is easy to see that a morphism $\varphi\colon\A\to\B$ of graded exterior
$\Ptens$-algebras is an isomorphism if and only if it is
an isomorphism in degrees $0$ and $1$.

\begin{example}
\label{example:Stein_exterior}
Let $(V,\O_V)$ be a Stein space, and let $\F$ be a locally free sheaf
of $\O_V$-modules. Set $A=\O(V)$ and $F=\F(V)$.
We have an obvious sheaf-theoretic version of the
antisymmetrization map \eqref{asymm}
$$
\tens{\O_V}\nolimits^n \F \xra{\mathrm a_\F} \tens{\O_V}\nolimits^n \F.
$$
By definition, $\shIm\mathrm a_\F=\bigwedge^n_{\O_V} \F$.
Since the functor $\Gamma$ of global sections takes tensor products over
$\O_V$ to projective tensor products over $A$ (see above),
we see that the morphism $\Gamma(V,\mathrm a_\F)$ of Stein $A$-modules
coincides with $\mathrm a_F\colon \ptens{A}^n F\to\ptens{A}^n F$.
Since $\Gamma$ is exact, we conclude that
$\pbigwedge{A}^n F=\Im\mathrm a_F=\Gamma(V,\shIm\mathrm a_\F)=
\Gamma(V,\bigwedge_{\O_V}^n \F)$. Thus we have an isomorphism
of graded $\Ptens$-algebras
$\pbigwedge{A} F \cong \Gamma(V,\bigwedge_{\O_V} \F)$.
\end{example}

By a {\em differential graded $\Ptens$-algebra} (a {\em DG $\Ptens$-algebra}
for short) we mean a graded $\Ptens$-algebra $\A$ together with a sequence
$\{ d^p\colon A^p\to A^{p+1} : p\in \Z_+\}$ of linear continuous maps such that
$d^{p+1} d^p=0$ for all $p$ (so that $\A$ becomes a cochain complex),
and $d^{p+q}(ab)=d^p(a)b+(-1)^p ad^q(b)$
for each $a\in A^p,\; b\in A^q$.
In particular, $d^0$ is a derivation of $A^0$ with values in $A^1$.
A DG $\Ptens$-algebra is said to be graded commutative
(exterior, etc.) if it has this property when considered as a graded
$\Ptens$-algebra. Morphisms of DG $\Ptens$-algebras
are morphisms of graded $\Ptens$-algebras commuting with differentials.

Let $\Omega^1 A$ be the module of K\"ahler
differentials of a commutative $\Ptens$-algebra $A$.
Then the exterior algebra $\pbigwedge{A}(\Omega^1 A)$ has a unique
structure of DG $\Ptens$-algebra such that the mapping
$d^0\colon A\to\Omega^1 A$ coincides with the universal
derivation $d_A$
(cf. \cite{Bour_HA,Kunz_new}).
The resulting DG $\Ptens$-algebra is denoted
by $\Omega(A)$ and is called the {\em algebra of differential forms}
of $A$. The cohomology groups of
$\Omega(A)$ are called the {\em de Rham cohomology groups} of $A$
and are denoted by $H^n_{DR}(A)$.

The algebra of differential forms has the following universal
property (cf. \cite{Kunz_new}):
{\em For each graded commutative DG $\Ptens$-algebra $\B$
and each $\Ptens$-algebra morphism $\psi\colon A\to B^0$ there exists a unique
DG $\Ptens$-algebra morphism $\varphi\colon\Omega(A)\to\B$
such that $\varphi^0=\psi$.}
In particular, each morphism $\psi\colon A\to B$ of $\Ptens$-algebras
uniquely extends to a morphism $\psi_*\colon\Omega(A)\to\Omega(B)$
of DG $\Ptens$-algebras. Thus the assignment $A\mapsto\Omega(A)$
is a functor from the category of commutative $\Ptens$-algebras
to the category of graded commutative DG $\Ptens$-algebras.

\begin{prop}
\label{prop:DR_Stein}
Let $V$ be a Stein manifold. Then the topological cohomology groups
$H^n_{\tp}(V,\CC)$ coincide with the de Rham cohomology groups
$H^n_{DR}(\O(V))$ of the Fr\'echet algebra $\O(V)$.
\end{prop}
\begin{proof}
For each $n$ denote by $\Omega^n_V$ the sheaf of holomorphic differential
$n$-forms on $V$. By the Poincar\'e lemma and Cartan's theorem B,
the de Rham complex
$$
0 \to \CC \to \O_V \to \Omega^1_V \to \Omega^2_V \to \cdots
$$
is an acyclic resolution of the constant sheaf $\CC$.
Therefore the cohomology groups of
$\Omega(V)=\Gamma(V,\Omega^\bullet_V)$
coincide with the topological cohomology groups $H^n_{\tp}(V,\CC)$.
On the other hand, the embedding $\O(V)\to\Omega(V)$ uniquely
extends to a DG $\Ptens$-algebra morphism $\Omega(\O(V))\to\Omega(V)$
that is an isomorphism in degrees $0$ and $1$
(see Lemma~\ref{lemma:Stein_diff}). Since both the algebras are exterior
(see Example~\ref{example:Stein_exterior}),
we conclude that $\Omega(\O(V))\to\Omega(V)$ is a DG $\Ptens$-algebra
isomorphism. The rest is clear.
\end{proof}

\subsection{Lie algebra actions}
\label{subsect:Lie_act}
Throughout the paper, by a Lie algebra we always mean a finite-dimensional
complex Lie algebra.

Let $\fg$ be a Lie algebra and $M$ a right $\fg$-module. For each $n\in\Z_+$
set $C_n(\fg,M)=M\Tens \bigwedge^n\fg$. The boundary mappings
$d_n\colon C_{n+1}(\fg,M)\to C_n(\fg,M)$
are defined by
\begin{multline*}
d_n(m\otimes X_1\wedge\cdots\wedge X_{n+1})=
\sum_{i=1}^{n+1} (-1)^{i-1} m\cdot X_i\otimes X_1\wedge\cdots\wedge\hat X_i\wedge
\cdots\wedge X_{n+1}\\
+\sum_{1\le i<j\le {n+1}} (-1)^{i+j} m\otimes [X_i,X_j]\wedge X_1\wedge\cdots\wedge\hat X_i
\wedge\cdots\wedge\hat X_j\wedge\cdots\wedge X_{n+1}.
\end{multline*}
(Here, as usual, the notation $\hat X_i$ indicates that $X_i$ is omitted.)
The spaces $C_n(\fg,M)$ together with the mappings $d_n$ form a chain
complex denoted by $C_\cdot(\fg,M)$. The homology groups of this complex
are called the {\em homology groups of $\fg$ with coefficients in $M$}
and are denoted by $H_n^{\Lie}(\fg,M)$.

Now let $U(\fg)$ be the universal enveloping algebra of $\fg$. We consider
$U(\fg)$ as a right $\fg$-module w.r.t. the right regular representation given
by $(a,X)\mapsto aX$. The corresponding chain complex
$V_\cdot(\fg)=C_\cdot(\fg,U(\fg))$
augmented by the counit map $\eps\colon U(\fg)\to\CC$
is exact (see \cite{CE}, Chap. XIII).
Since all the $d_n$'s are morphisms of left $U(\fg)$-modules in this case,
it follows that $V_\cdot(\fg)$ is a free
resolution of the trivial $\fg$-module $\CC$ in the category
of left $U(\fg)$-modules (the {\em Koszul resolution}). If $M$ is a right
$\fg$-module, then $C_\cdot(\fg,M)$ is isomorphic to the tensor product
$M\tens{U(\fg)} V_\cdot(\fg)$. Therefore
$H_n^{\Lie}(\fg,M)=\Tor_n^{U(\fg)}(M,\CC)$ for each $n\in\Z_+$.

Dually, if $M$ is a left $\fg$-module, then the space $C^n(\fg,M)$
of {\em $n$-cochains} is defined as $\Hom_{\CC}(\bigwedge^n\fg,M)$.
Thus $n$-cochains are just alternating multilinear maps of $n$ variables
from $\fg$ with values in $M$. The coboundary mappings
$d^n\colon C^n(\fg,M)\to C^{n+1}(\fg,M)$ are defined by
\begin{multline*}
d^n f(X_1\wedge\cdots\wedge X_{n+1})=
\sum_{i=1}^{n+1} (-1)^{i-1} X_i \cdot f(X_1\wedge\cdots\wedge\hat X_i\wedge
\cdots\wedge X_{n+1})\\
+\sum_{1\le i<j\le n+1} (-1)^{i+j} f([X_i,X_j]\wedge X_1\wedge\cdots\wedge\hat X_i
\wedge\cdots\wedge\hat X_j\wedge\cdots\wedge X_{n+1}).
\end{multline*}
The spaces $C^n(\fg,M)$ together with the mappings $d^n$ form a cochain
complex denoted by $C^\cdot(\fg,M)$ (the {\em Chevalley-Eilenberg complex}).
The cohomology groups of this complex
are called the {\em cohomology groups of $\fg$ with coefficients in $M$}
and are denoted by $H^n_{\Lie}(\fg,M)$.
As in the case of homology groups, we have
$H^n_{\Lie}(\fg,M)=\Ext^n_{U(\fg)}(\CC,M)$ for each $n\in\Z_+$.

\begin{remark}
Recall that each right $\fg$-module $M$ can also be viewed as a left $\fg$-module
w.r.t. the action $X\cdot m=-m\cdot X\; (m\in M,\; X\in\fg)$, and vice versa,
so that the categories of left $\fg$-modules and of right $\fg$-modules are
isomorphic. Thus one can speak about the complex $C_\cdot(\fg,M)$
(resp. $C^\cdot(\fg,M)$) in the case where $M$ is a left (resp. right) $\fg$-module.
\end{remark}

By a {\em right $\fg$-$\Ptens$-module}  we mean a complete Hausdorff
l.c.s. $M$ together with the structure of right $\fg$-module such that
the map $M\to M,\; m\mapsto m\cdot X$ is continuous for each $X\in\fg$.
If we endow $\fg$ with the usual topology of a finite-dimensional
vector space, then the above condition means precisely
that the map $M\Ptens\fg\to M,\; m\otimes X\mapsto m\cdot X$ is continuous.
Similarly, one can speak about left $\fg$-$\Ptens$-modules.
If $M$ is a right (resp. left) $\fg$-$\Ptens$-module, then the strong dual,
$M'$, becomes a left (resp. right) $\fg$-$\Ptens$-module via the action
$\la m,X\cdot m'\ra=\la m\cdot X,m'\ra$
(resp. $\la m'\cdot X,m\ra=\la m',X\cdot m\ra$) for
$m\in M,\; m'\in M',\; X\in\fg$.

\begin{remark}
If we endow $U(\fg)$ with the finest locally convex topology, then
each $\fg$-$\Ptens$-module $M$ becomes a topological $U(\fg)$-module.
Note, however, that $M$ need not be an $U(\fg)$-$\Ptens$-module,
i.e., the action $U(\fg)\times M\to M$ need not be jointly continuous.
\end{remark}

If $M$ is a right (resp. left) $\fg$-$\Ptens$-module, then
the obvious identifications $M\Tens\bigwedge^n\fg=M\Ptens\bigwedge^n\fg$
and $\Hom_{\CC}(\bigwedge^n\fg,M)=\L(\bigwedge^n\fg,M)$ enable us to consider
$C_\cdot(\fg,M)$
(resp. $C^\cdot(\fg,M)$) as a complex in $\mathbf{LCS}$.
If $M$ is a right (resp. left) $\fg$-$\Ptens$-module,
then the complex
$C^\cdot(\fg,M')$ (resp. $C_\cdot(\fg,M')$) is isomorphic to the
strong dual of $C_\cdot(\fg,M)$ (resp. $C^\cdot(\fg,M)$).
This readily follows from the canonical isomorphisms
$\L(\bigwedge^n\fg,M)\cong (\bigwedge^n\fg)'\Ptens M$.

Let $A$ be a $\Ptens$-algebra together with the structure of left
$\fg$-$\Ptens$-module. Suppose that
for each $X\in\fg$ the map $A\to A,\; a\mapsto X\cdot a$
is a derivation. In this case we say that {\em $\fg$ acts on $A$ by derivations}.
The complex $C^\cdot(\fg,A)$ has then a structure of DG $\Ptens$-algebra
(cf. \cite{Dubois}). The multiplication on $C^\cdot(\fg,A)$ comes from
the identification of $C^\cdot(\fg,A)$ with the tensor
product of algebras $\bigwedge\fg^*\Ptens A$. In particular,
if $A$ is commutative, then $C^\cdot(\fg,A)$ is isomorphic (as a graded
$\Ptens$-algebra) to the exterior algebra
$\pbigwedge{A} C^1(\fg,A)$.

\section{The inverse process for Hopf $\Ptens$-algebras}
\label{sect:inverse}
In this section we describe a version of the Cartan-Eilenberg
``inverse process'' (\cite{CE}, Chap.~X) adapted to the Hopf $\Ptens$-algebra
case. Originally, Cartan and Eilenberg applied the inverse process to the
study of homological dimensions of group algebras and universal enveloping
algebras. Subsequently some generalizations were obtained for
cocommutative \cite{Kielp} and commutative \cite{Linckel} Hopf algebras.
Though we believe that the algebraic versions of the results below
are known, we could not find them in the literature in a form suitable for
our purposes. That is why we give complete proofs.

For convenience of the reader, we recall some algebraic definitions
(see \cite{Majid} for details).
Let $\cC$ be a monoidal category,
i.e., a category equipped with a bifunctor $\Tens\colon\cC\times\cC\to\cC$,
a neutral object $I$, and natural isomorphisms
$$
a_{X,Y,Z}\colon (X\Tens Y)\Tens Z\to X\Tens (Y\Tens Z),\quad
l_X\colon I\Tens X\to X,
\quad r_X\colon X\Tens I\to X
$$
satisfying natural coherence (constraint)
conditions (see, e.g., \cite{ML}).
Without loss of generality (by MacLane's coherence theorem \cite[Theorem 15.1]{ML}),
we may assume that $\cC$ is {\em strict}, so that all associativity and unit
isomorphisms are identities.
An {\em algebra} in $\cC$ is an object $A$ together with morphisms
${\mu\colon A\Tens A\to A}$ (multiplication) and $\eta\colon I\to A$ (unit)
such that the diagrams
$$
\xymatrix{
A\Tens A \Tens A \ar[r]^(.56){\mu\otimes\id_A} \ar[d]_{\id_A\otimes\mu}
& A\Tens A \ar[d]^\mu \\
A\Tens A \ar[r]^\mu & A
}
\qquad
\xymatrix{
I\Tens A \ar[r]^{\eta\otimes\id_A} \ar[dr]_{l_A} & A\Tens A \ar[d]^\mu
& A\Tens I \ar[l]_{\id_A\otimes\eta} \ar[dl]^{r_A} \\
& A
}
$$
are commutative.
For example, if $\cC=\LCS$ is the category of complete
locally convex spaces over $\CC$, and $\Tens=\Ptens$
is the bifunctor of the completed
projective tensor product,
then we obtain the definition of $\Ptens$-algebra given in Subsection~\ref{subsect:topalg}.
If $(A,\mu_A,\eta_A)$ and $(B,\mu_B,\eta_B)$ are algebras in $\cC$, then a morphism
$\varphi\colon A\to B$ is an {\em algebra homomorphism} if
$\varphi\mu_A=\mu_B(\varphi\otimes\varphi)$ and $\varphi\eta_A=\eta_B$.

Dually, a {\em coalgebra} in
$\cC$ is an object $C$ together with morphisms
$\Delta\colon C\to C\Tens C$ (comultiplication) and $\eps\colon C\to I$ (counit)
such that the diagrams
$$
\xymatrix{
C\Tens C \Tens C & C\Tens C \ar[l]_(.42){\Delta\otimes\id_C} \\
C\Tens C \ar[u]^{\id_C\otimes\Delta}
& C \ar[l]_\Delta \ar[u]_\Delta
}
\qquad
\xymatrix{
I\Tens C & C\Tens C \ar[l]_{\eps\otimes\id_C} \ar[r]^{\id_C\otimes\eps} & C\Tens I\\
& C \ar[ul]^{l_C^{-1}} \ar[ur]_{r_C^{-1}} \ar[u]_\Delta
}
$$
are commutative.

The monoidal category $\cC$ is {\em braided} if it is
equipped with a natural isomorphism $c_{X,Y}\colon X\Tens Y\to Y\Tens X$
satisfying the relations
$$
(c_{X,Z}\otimes\id_Y)(\id_X\otimes c_{Y,Z})=c_{X\Tens Y,Z}
\quad\text{and}\quad
(\id_Y\otimes c_{X,Z})(c_{X,Y}\otimes\id_Z)=c_{X,Y\Tens Z}
$$
In this case, the tensor product of any two algebras $A, B$ in $\cC$ is an algebra
with multiplication and unit defined as the compositions
\begin{gather*}
A\Tens B\Tens A\Tens B \xra{\id_A\otimes c_{B,A}\otimes\id_B}
A\Tens A\Tens B\Tens B \xra{\mu_A\otimes\mu_B} A\Tens B,\\
I \xra{r_I=l_I} I\Tens I \xra{\eta_A\Tens\eta_B} A\Tens B.
\end{gather*}
A {\em bialgebra} in $\cC$ is an object $H$ equipped with the algebra and the coalgebra
structures such that $\Delta\colon H\to H\Tens H$ and $\eps\colon H\to I$
are algebra homomorphisms. Finally, a {\em Hopf algebra} in $\cC$ is a bialgebra $H$
together with a morphism $S\colon H\to H$ (antipode) such that
the diagram
$$
\xymatrix{
H\Tens H \ar[d]_{S\otimes\id_H} & H \ar[l]_(.4)\Delta \ar[r]^(.4)\Delta
\ar[d]_{\eta\eps} & H\Tens H \ar[d]^{\id_H\otimes S} \\
H\Tens H \ar[r]^(.6)\mu & H & H\Tens H \ar[l]_(.6)\mu
}
$$
is commutative.

\begin{lemma}
\label{lemma:PhiPsi}
Let $H$ be a Hopf algebra in a braided monoidal category $\cC$, and let
$\Phi,\Psi\colon H\Tens H\to H\Tens H$ be given by
$$
\Phi=(\mu\otimes\id_H)(\id_H\otimes\Delta)\quad\text{and}\quad
\Psi=(\mu\otimes\id_H)(\id_H\otimes S\otimes\id_H)(\id_H\otimes\Delta).
$$
Then $\Phi=\Psi^{-1}$.
\end{lemma}
\begin{proof}
The relation $\Phi\Psi=\id_{H\Tens H}$
follows from the commutative diagram
$$
\xymatrix@C+5mm{
{H\Tens H
\ar `u[rr] `[rr]^\Phi [rr]
\ar[r]^{\id\otimes\Delta} }
& {H\Tens H\Tens H \ar[r]^{\mu\otimes\id}}
& {H\Tens H} \\
{H\Tens H\Tens H \ar[r]^(.45){\id\otimes\id\otimes\Delta}
\ar[u]_{\mu\otimes\id}}
& {H\Tens H\Tens H\Tens H \ar[r]^(.55){\id\otimes\mu\otimes\id}
\ar[u]_{\mu\otimes\id\otimes\id}}
& {H\Tens H\Tens H \ar[u]^{\mu\otimes\id}}\\
{H\Tens H\Tens H \ar[r]^(.45){\id\otimes\id\otimes\Delta}
\ar[u]_{\id\otimes S\otimes\id}}
& {H\Tens H\Tens H\Tens H \ar[u]_{\id\otimes S\otimes\id\otimes\id}}\\
{H\Tens H
\ar `d[r] `[rr]_{\id\otimes l_H^{-1}} [rr]
\ar `l[u] `[uuu]^\Psi [uuu]
\ar[u]_{\id\otimes\Delta}
\ar[r]_{\id\otimes\Delta}
}
& {H\Tens H\Tens H \ar[u]_{\id\otimes\Delta\otimes\id}
\ar[r]_{\id\otimes\varepsilon\otimes\id}}
& {H\Tens I\otimes H \ar[uu]^{\id\otimes\eta\otimes\id}
\ar `r[uuu] `[uuu]_{r_H\otimes\id} [uuu]
}
}
$$
Similarly, the commutative diagram
$$
\xymatrix@C+5mm{
**{+<20pt,0pt>} {H\Tens H
\ar `u[r] `[rr]^\Phi [rr]
\ar[r]^{\id\otimes\Delta}
\ar `l[d] `[ddd]_{\id\otimes l_H^{-1}} [ddd]
\ar[d]^{\id\otimes\Delta} }
& H\Tens H\Tens H \ar[r]^{\mu\otimes\id}
\ar[d]_{\id\otimes\id\otimes\Delta}
& H\Tens H \ar[d]_{\id\otimes\Delta}
\ar `r/10pt[d] `[ddd]^\Psi [ddd]
\\
H\Tens H\Tens H \ar[r]^(.45){\id\otimes\Delta\otimes\id}
\ar[dd]^{\id\otimes\varepsilon\otimes\id}
& H\Tens H\Tens H\Tens H \ar[r]^(.55){\mu\otimes\id\otimes\id}
\ar[d]_{\id\otimes\id\otimes S\otimes\id}
& H\Tens H\Tens H \ar[d]_{\id\otimes S\otimes\id}\\
& H\Tens H\Tens H\Tens H \ar[r]^(.55){\mu\otimes\id\otimes\id}
\ar[d]_{\id\otimes\mu\otimes\id}
& H\Tens H\Tens H \ar[d]_{\mu\otimes\id}\\
H\Tens I\Tens H
\ar[r]^{\id\otimes\eta\otimes\id}
\ar `d[r] `[rr]_{r_H\otimes\id} [rr]
& H\Tens H\Tens H \ar[r]^{\mu\otimes\id}
& H\Tens H
}
$$
shows that $\Psi\Phi=\id_{H\Tens H}$.
\end{proof}

Let $A$ be an algebra in $\cC$. Recall that a {\em left $A$-module}
is an object $M$ together with a morphism $\mu_M\colon A\Tens M\to M$
such that the diagrams
$$
\xymatrix{
A\Tens A \Tens M \ar[r]^(.58){\mu_A\otimes\id_M} \ar[d]_{\id_A\otimes\mu_M}
& A\Tens M \ar[d]^{\mu_M} \\
A\Tens M \ar[r]^{\mu_M} & M
}
\qquad
\xymatrix{
I\Tens M \ar[r]^{\eta\otimes\id_M} \ar[dr]_{l_M} & A\Tens M \ar[d]^{\mu_M}\\
& M
}
$$
are commutative. Again, in the case $(\cC,\Tens)=(\LCS,\Ptens)$ we
obtain the definition of $\Ptens$-module given in Subsection~\ref{subsect:topalg}.
If $(M,\mu_M)$ and $(N,\mu_N)$ are left $A$-modules, then a morphism
$\varphi\colon M\to N$ in $\cC$ is an {\em $A$-module morphism} if
$\varphi\mu_M=\mu_N(\id_A\otimes\varphi)$. Right $A$-modules and their
morphisms are defined similarly.

Now let $\cC$ be a braided monoidal category, and let $H$ be a Hopf algebra in $\cC$.
Then $H\Tens H$ has two natural structures
of right $H$-module. The first one is given by the action of $H$ on the right
factor, and the second one arises from the algebra homomorphism
$\Delta\colon H\to H\Tens H$. Thus we obtain the right $H$-modules
$(H\Tens H,\mu^r)$ and $(H\Tens H,\mu^\Delta)$ with the actions
$\mu^r,\mu^\Delta\colon (H\Tens H)\Tens H\to (H\Tens H)$ given, respectively, by
\begin{gather*}
\mu^r\colon (H\Tens H)\Tens H \lriso H\Tens (H\Tens H) \xra{\id_H\otimes\mu}
H\Tens H\\
\mu^\Delta\colon (H\Tens H)\Tens H \xra{\id_{H\otimes H}\otimes\Delta}
(H\Tens H)\Tens (H\Tens H) \xra{\mu_{H\otimes H}} H\Tens H.
\end{gather*}

To simplify notation, we shall often write $H^{(n)}$ to denote the
$n$-fold tensor power $H\Tens\cdots\Tens H$.

\begin{lemma}
\label{lemma:free1}
The morphism $\Phi\colon (H\Tens H,\mu^r)\to (H\Tens H,\mu^\Delta)$
defined in Lemma \ref{lemma:PhiPsi} is an isomorphism of right $H$-modules.
\end{lemma}
\begin{proof}
By Lemma \ref{lemma:PhiPsi}, $\Phi$ is an isomorphism in $\cC$.
To prove that $\Phi$ is a right $H$-module morphism, it is
enough to consider the following commutative diagram:
$$
\xymatrix@C+10pt{
H\Tens H\Tens H \ar[r]^(.45){\id_H\otimes\Delta\otimes\id_H}
\ar `u[r] `[rr]^{\Phi\otimes\id_H} [rr]
\ar[dr]_{\id_H\otimes\Delta\otimes\Delta}
\ar[dd]_{\mu^r=\id_H\otimes\mu_H}
& H\Tens H\Tens H\Tens H \ar[r]^(.55){\mu_H\otimes\id_H\otimes\id_H}
\ar[d]^{\id_{H^{(3)}}\otimes\Delta}
& H\Tens H\Tens H \ar[d]_{\id_{H^{(2)}}\otimes\Delta}
\ar `r[d] `[dd]^{\mu^\Delta} [dd]
\\
& H\Tens H\Tens H\Tens H\Tens H
\ar[r]^(.55){\mu_H\otimes\id_{H^{(3)}}}
\ar[d]_{\id_H\otimes\mu_{H^{(2)}}}
& H\Tens H\Tens H\Tens H
\ar[d]_{\mu_{H^{(2)}}}
\\
H\Tens H \ar[r]^(.45){\id_H\otimes\Delta}
\ar `d[r] `[rr]_{\Phi} [rr]
& H\Tens H\Tens H \ar[r]^{\mu_H\otimes\id_H}
& H\Tens H
}
$$
\end{proof}

Let $H^{\op}$ denote the algebra opposite to $H$, i.e.,
$\mu_{H^{\op}}=\mu_H c_{H,H}$, and let $H^e=H\Tens H^{\op}$.
Then $S\colon H\to H^{\op}$ and $E=(\id_H\Tens S)\Delta\colon H\to H^e$
are algebra homomorphisms \cite{Majid}.
Hence $H^e$ becomes a right $H$-module via $E$.
We denote this module by $H^e_E$.

\begin{lemma}
\label{lemma:free2}
If $H$ has invertible antipode, then
the right $H$-modules $(H\Tens H,\mu^r)$ and $H^e_E$
are isomorphic.
\end{lemma}
\begin{proof}
Since $S$ is an isomorphism in $\cC$, it follows that
$\id_H\Tens S\colon H\Tens H\to H\Tens H^{\op}=H^e$ is an algebra isomorphism.
Hence $H^e_E\cong (H\Tens H,\mu^\Delta)$ as right $H$-modules.
Now it remains to apply Lemma \ref{lemma:free1}.
\end{proof}

From now on, let $(\cC,\Tens)=(\LCS,\Ptens)$ be the category
of complete locally convex spaces over $\CC$.
By a {\em Hopf $\Ptens$-algebra} we mean a Hopf
algebra in $\LCS$ (cf. also \cite{Litv_Hopf,BFGP,Pfl-Schot}).
Given a Hopf $\Ptens$-algebra $H$, we consider $\CC$ as a left $H$-module
via the counit map $\varepsilon\colon H\to\CC$.

\begin{lemma}
\label{lemma:tens_iso}
Let $H$ be a Hopf $\Ptens$-algebra with invertible antipode.
There exists an isomorphism
of left $H^e$-$\Ptens$-modules
\begin{equation}
\label{tens_iso}
\varphi\colon H^e_E\ptens{H}\CC\to H,\qquad
u\otimes 1\mapsto \mu(u).
\end{equation}
\end{lemma}
\begin{proof}
Consider the bilinear map
$$
R\colon H^e\times\CC\to H,\qquad (u,\lambda)\mapsto \lambda\mu(u).
$$
To prove that \eqref{tens_iso} is a well-defined linear map, we have
to show that
\begin{equation}
\label{R_balanced}
R\bigl((a\otimes b)\cdot c,\lambda\bigr)=R(a\otimes b,c\cdot\lambda)
\end{equation}
for each $a,b,c\in H$ and each $\lambda\in\CC$.
To this end, note first that
\begin{equation}
\label{muE}
\mu E=\mu (\id_H\otimes S)\Delta=\eta\varepsilon.
\end{equation}
Since $\mu\colon H^e\to H$ is a left $H^e$-module morphism,
we see that
\begin{equation}
\label{Rleft}
R\bigl((a\otimes b)\cdot c,\lambda\bigr)=\lambda\mu\bigl((a\otimes b)E(c)\bigr)=
\lambda a\mu\bigl(E(c)\bigr) b=\varepsilon(c)\lambda ab.
\end{equation}
On the other hand,
\begin{equation}
\label{Rright}
R(a\otimes b,c\cdot\lambda)=R\bigl(a\otimes b,\varepsilon(c)\lambda\bigr)=
\varepsilon(c)\lambda ab.
\end{equation}
Comparing \eqref{Rleft} and \eqref{Rright}, we obtain \eqref{R_balanced},
as required. Hence $\varphi$ is a well-defined, linear continuous map.
Evidently, $\varphi$ is also a left $H^e$-module morphism.

To construct the inverse of $\varphi$, consider the map
$$
\psi\colon H\to H^e_E\ptens{H}\CC,\qquad a\mapsto (a\otimes 1)\otimes 1.
$$
Clearly, $\varphi\psi=\id_H$. Thus it remains to prove that
$\psi\varphi=\id_{H^e_E\ptens{H}\CC}$, which is equivalent to
\begin{equation}
\label{psi_phi}
u\otimes 1=\bigl(\mu(u)\otimes 1\bigr)\otimes 1
\end{equation}
for each $u\in H^e$.

Take the map $\Phi\colon H\Ptens H\to H\Ptens H$ defined in
Lemma \ref{lemma:PhiPsi}, and set
$$
\Phi'=(\id_H\otimes S)\Phi\colon (H\Ptens H,\mu^r)\to H^e_E.
$$
By Lemmas \ref{lemma:free1} and \ref{lemma:free2}, $\Phi'$ is an
isomorphism of right $H$-$\Ptens$-modules. We have
\begin{multline*}
\Phi'(a\otimes 1)=(\id_H\otimes S)(\mu\otimes\id_H)(\id_H\otimes\Delta)(a\otimes 1)\\
=(\id_H\otimes S)(\mu\otimes\id_H)(a\otimes 1\otimes 1)
=(\id_H\otimes S)(a\otimes 1)=a\otimes 1
\end{multline*}
and hence
$$
\Phi'(a\otimes b)=\Phi'(a\otimes 1\cdot b)=\Phi'(a\otimes 1)E(b)=
(a\otimes 1)E(b).
$$
Since $\Phi'$ is bijective, it is enough to check \eqref{psi_phi}
with $u=(a\otimes 1)E(b)$. Using \eqref{muE} and the fact that
$\mu$ is a left $H^e$-module morphism, we see that
$$
\mu(u)=\mu\bigl((a\otimes 1)E(b)\bigr)=
a\mu\bigl(E(b)\bigr)=
\varepsilon(b)a.
$$
Hence the right-hand side of \eqref{psi_phi} is
$$
\bigl(\mu(u)\otimes 1\bigr)\otimes 1=
\varepsilon(b)(a\otimes 1)\otimes 1,
$$
while the left-hand side of \eqref{psi_phi} is
$$
u\otimes 1=(a\otimes 1)E(b)\otimes 1=
(a\otimes 1)\otimes b\cdot 1=
\varepsilon(b)(a\otimes 1)\otimes 1.
$$
Therefore \eqref{psi_phi} is satisfied, and so $\psi=\varphi^{-1}$.
Hence $\varphi$ is an isomorphism, as required.
\end{proof}

\begin{theorem}
\label{thm:inverse}
Let $H$ be a Hopf $\Ptens$-algebra with invertible antipode, and let
$$
0 \lla \CC \lla P_\bullet
$$
be a projective resolution of $\CC$ in $H\lmod$.
Then the tensor product complex
\begin{equation}
\label{bimodres}
0 \lla H \lliso H^e_E\ptens{H} \CC \lla H^e_E\ptens{H} P_\bullet
\end{equation}
is a projective bimodule resolution of $H$.
\end{theorem}
\begin{proof}
By Lemma \ref{lemma:free2}, $H^e_E$ is a free right $H$-$\Ptens$-module.
Hence the augmented complex \eqref{bimodres} is admissible.
To complete the proof, it remains to apply Lemma \ref{lemma:tens_iso}.
\end{proof}

Let $H$ be a Hopf $\Ptens$-algebra and $M$ an $H$-$\Ptens$-bimodule
(i.e., a left $H^e$-module). We may consider $M$ as a left $H$-$\Ptens$-module via
$E\colon H\to H^e$. Similarly, by considering $M$ as a right $H^e$-$\Ptens$-module,
we obtain a right $H$-$\Ptens$-module structure on $M$.
The resulting left (resp. right) $H$-$\Ptens$-module will be denoted
by ${_E}M$ (resp. $M_E$).

\begin{corollary}
\label{cor:H-Ext}
Let $H$ be a Hopf $\Ptens$-algebra with invertible antipode.
Then for each $M\in H\bimod H$ there exist natural isomorphisms
$$
\H^n(H,M)\cong\Ext^n_H(\CC,{_E}M)\quad\text{and}\quad
\H_n(H,M)\cong\Tor_n^H(M_E,\CC).
$$
\end{corollary}
\begin{proof}
Let $P_\bullet$ be a projective resolution of $\CC$ in $H\lmod$.
In view of Theorem \ref{thm:inverse}, we have
$$
\H^n(H,M)=\HH^n\bigl({_{H^e}}\h(H^e_E\ptens{H} P_\bullet,M)\bigr)\cong
\HH^n\bigl({_H}\h(P_\bullet,{_E}M)\bigr)=\Ext^n_H(\CC,{_E}M).
$$
Similarly,
$$
\H_n(H,M)=\HH_n(M\ptens{H^e} H^e_E\ptens{H}P_\bullet)\cong
\HH_n(M_E\ptens{H} P_\bullet)=\Tor_n^H(M_E,\CC).
$$
\end{proof}

\begin{corollary}
Let $H$ be a Hopf $\Ptens$-algebra with invertible antipode. Then
$\dh_H\CC=\dg H=\db H$.
\end{corollary}

The following two examples are ``continuous versions'' of Cartan-Eilenberg's
result on the Hochschild cohomology of group algebras (\cite{CE}, Chap.~X, \S6).

\begin{example}
Let $G$ be a discrete group. The Banach algebra $\ell^1(G)$ has
a canonical Hopf $\Ptens$-algebra structure uniquely determined by
$$
\Delta(\delta_g)=\delta_g\otimes\delta_g,\quad
\eps(f)=\sum_{g\in G}f(g),\quad
Sf(g)=f(g^{-1}).
$$
(Here $\delta_g$ denotes the function which is $1$ at $g\in G$, $0$ elsewhere.)
Using the bar resolution of $\CC$ in $\ell^1(G)\lmod$ (see \cite{X1}),
it is easy to check that $\Ext^n_{\ell^1(G)}(\CC,X)$ is isomorphic to $H^n_b(G,X)$,
the $n$th bounded cohomology group of $G$ with coefficients in $X$
(\cite{Noskov}; cf. also \cite{Johnson}).
Thus we obtain the following
\end{example}

\begin{corollary}
Let $G$ be a discrete group and $M$ a Banach $\ell^1(G)$-bimodule.
Denote by ${_E}M$ the left $G$-module obtained from $M$ by setting
$gm=\delta_g\cdot m\cdot\delta_{g^{-1}}\; (g\in G,\; m\in M)$.
Then there exist canonical isomorphisms
$$
\H^n(\ell^1(G),M)\cong H^n_b(G,{_E}M).
$$
\end{corollary}

\begin{example}
Let $G$ be a real Lie group. The convolution algebra $\E'(G)$
of compactly supported distributions on $G$ is a Hopf $\Ptens$-algebra
in a natural way (see, e.g., \cite{Litv,BFGP}; cf. also Section \ref{sect:an_func}
below). Let $X$ be a left $\E'(G)$-$\Ptens$-module. As in the previous example,
it can easily be checked
that $\Ext^n_{\E'(G)}(\CC,X)$ is isomorphic to $H^n_c(G,X)$, the $n$th
continuous (or, equivalently, differentiable) cohomology group
of $G$ with coefficients in $X$ (cf. \cite{Guichardet}, Chap.~III, Prop.~1.5).
Thus we obtain the following
\end{example}

\begin{corollary}
Let $G$ be a real Lie group and $M$ an $\E'(G)$-$\Ptens$-bimodule.
Denote by ${_E}M$ the left $G$-module obtained from $M$ by setting
$gm=\delta_g\cdot m\cdot\delta_{g^{-1}}\; (g\in G,\; m\in M)$.
Then there exist canonical isomorphisms
$$
\H^n(\E'(G),M)\cong H^n_c(G,{_E}M).
$$
\end{corollary}

We end this section with an application of the above results to
left amenability in the sense of Lau \cite{Lau_F-alg}.
Recall that a Banach algebra $A$ is said to be {\em amenable} \cite{Johnson} if
$\H^1(A,X^*)=0$ for each Banach $A$-bimodule $X$, i.e., if every derivation from
$A$ to $X^*$ is inner. Suppose $A$ is endowed with an augmentation $\eps_A$
(i.e., a continuous homomorphism $A\to\CC$). Then $A$ is said to be
{\em left amenable} \cite{Lau_F-alg} if for each Banach $A$-bimodule $X$ such that
$a\cdot x=\eps_A(a)x$ for all $a\in A,\; x\in X$, every derivation from
$A$ to $X^*$ is inner.

In the next lemma, we consider $\CC$ as a left Banach $A$-module via
$\eps_A\colon A\to\CC$.

\begin{lemma}
Let $A$ be an augmented Banach algebra. Then $A$ is left amenable
if and only if $\Ext^1_A(\CC,Y^*)=0$ for each right Banach $A$-module $Y$.
\end{lemma}
\begin{proof}
Obviously, the $A$-bimodules in the definition of left amenability
are precisely those of the form $X=\CC\Ptens Y$ where $Y\in\rmod A$.
Hence $X^*\cong\L(\CC,Y^*)$ (see \cite[II.5.21]{X1}), and so
$\H^1(A,X^*)\cong\Ext^1_A(\CC,Y^*)$ (see \cite[III.4.12]{X1}).
The rest is clear.
\end{proof}

\begin{prop}
Let $H$ be a Banach Hopf algebra (i.e., a Hopf $\Ptens$-algebra whose underlying
locally convex space is a Banach space) with invertible antipode.
Then $H$ is left amenable if and only if $H$ is amenable.
\end{prop}
\begin{proof}
The ``if'' part is clear. Conversely, assume $H$ is left amenable, and let $X$
be a Banach $H$-bimodule. By Corollary~\ref{cor:H-Ext}, we have
$\H^1(H,X^*)\cong\Ext^1_H(\CC,{_E}(X^*))$. On the other hand,
it is immediate that ${_E}(X^*)=(X_E)^*$. Now the result follows from
the previous lemma.
\end{proof}

\section{Localizations and weak localizations}

Let $A$ be a Fr\'echet algebra, $X\in\rmod A$, and $Y\in A\lmod$. Then
$X$ and $Y$ are said to be {\em transversal} over $A$ (notation:
$X\perp_A Y$) if $\Tor_0^A(X,Y)$ is Hausdorff, and
$\Tor_n^A(X,Y)=0$ for all $n>0$. This notion was introduced in \cite{Ast}
and has proved to be extremely useful in complex analytic geometry
and operator theory \cite{KV,Ast,Eschm_Put,Demailly}.
We shall need a somewhat stronger condition of transversality type.

\begin{prop}
\label{prop:str_transv}
Let $A$ be a $\Ptens$-algebra, $X\in \rmod A$, and $Y\in A\lmod$.
The following conditions are equivalent:

\smallskip\noindent
{\upshape (i)}
There exists a projective resolution
\begin{equation}
\label{X-P}
0 \lar X \lar P_\bullet
\end{equation}
of $X$ in $\rmod A$ such that the tensored complex
\begin{equation}
\label{XY-PY}
0 \lar X\ptens{A} Y \lar P_\bullet\ptens{A} Y
\end{equation}
is admissible.

\smallskip\noindent
{\upshape (i)$'$}
For each projective resolution
\eqref{X-P} of $X$ in $\rmod A$ the complex
\eqref{XY-PY} is admissible.

\smallskip\noindent
{\upshape (ii)}
There exists a projective resolution
\begin{equation}
\label{Y-Q}
0 \lar Y \lar Q_\bullet
\end{equation}
of $Y$ in $A\lmod$ such that the tensored complex
\begin{equation}
\label{XY-XQ}
0 \lar X\ptens{A} Y \lar X\ptens{A} Q_\bullet
\end{equation}
is admissible.

\smallskip\noindent
{\upshape (ii)$'$}
For each projective resolution
\eqref{Y-Q} of $Y$ in $A\lmod$ the complex
\eqref{XY-XQ} is admissible.
\end{prop}
\begin{proof}
The equivalences (i)$\iff$(i)$'$ and (ii)$\iff$(ii)$'$
readily follow from the fact that every two projective resolutions
of a $\Ptens$-module are homotopy equivalent
(see \cite{X1}).

Let us prove that (i)$\iff$(ii).
Choose a projective resolution
\begin{equation}
\label{A-L}
0 \lar A \lar L_\bullet
\end{equation}
of $A$ in $A\bimod A$. Then the complexes
\begin{gather*}
0 \lar X \lar X\ptens{A} L_\bullet\\
0 \lar Y \lar L_\bullet\ptens{A} Y
\end{gather*}
are projective resolutions of $X\in\rmod A$ and $Y\in A\lmod$,
respectively. Since (i)$\iff$(i)$'$ and (ii)$\iff$(ii)$'$,
we see that both (i) and (ii) are equivalent
to the admissibility of the complex
$$
0 \lar X\ptens{A} Y \lar X\ptens{A} L_\bullet \ptens{A} Y.
$$
Therefore (i)$\iff$(ii).
\end{proof}

\begin{definition}
We say that $X\in\rmod A$ and $Y\in A\lmod$ are {\em strongly transversal} over $A$
if they satisfy the conditions of Proposition \ref{prop:str_transv}.
In this case, we write $X\pperp_A Y$.
\end{definition}

\begin{remark}
Suppose that $A$ is a Fr\'echet algebra.
If we require that \eqref{XY-PY} or \eqref{XY-XQ} be only exact
(but not necessarily admissible), then we come to the usual definition
of transversality (see the beginning of this section).
\end{remark}

\begin{prop}
\label{prop:loc}
Let $\theta\colon A\to B$ be a homomorphism of $\Ptens$-algebras.
Suppose that the map
\begin{equation}
\label{local2}
B\ptens{A} B\to B,\qquad b_1\otimes b_2\mapsto b_1 b_2
\end{equation}
is a topological isomorphism.
Then the following conditions are equivalent:

\smallskip\noindent
\upshape{(i)} $B\pperp_A B$;

\smallskip\noindent
\upshape{(ii)} $B\pperp_A M$ for each $M\in B\lmod$;

\smallskip\noindent
\upshape{(iii)} $M\pperp_A B$ for each $M\in \rmod B$;

\smallskip\noindent
\upshape{(iv)} $B^e\pperp_{A^e} A$.
\end{prop}
\begin{proof}
(ii)$\Longrightarrow$(i), (iii)$\Longrightarrow$(i): this is clear.

To prove the remaining implications,
take a projective bimodule resolution
\eqref{A-L} of $A$ in $A\bimod A$, and note that
\begin{equation}
\label{B-BL}
0 \lar B \lar L_\bullet\ptens{A} B
\end{equation}
is a projective resolution of $B\in A\lmod$.

(i)$\Longrightarrow$(iv).
If (i) holds, then the complex
\begin{equation}
\label{BB-BLB}
0 \lar B\ptens{A} B \lar B\ptens{A} L_\bullet \ptens{A} B
\end{equation}
is admissible. On the other hand, the latter complex is isomorphic to
$$
0\lar B^e\ptens{A^e} A \lar B^e\ptens{A^e} L_\bullet,
$$
and we obtain (iv).

(iv)$\Longrightarrow$(iii). If (iv) holds, then the complex
\eqref{BB-BLB} is admissible. Since $B\ptens{A} B\cong B$ is projective
in $B\lmod$, we see that \eqref{BB-BLB} splits in $B\lmod$.
Hence $M\ptens{B}$\eqref{BB-BLB} is admissible. On the other
hand, $M\ptens{B}$\eqref{BB-BLB} is isomorphic to $M\ptens{A}$\eqref{B-BL},
and we obtain (iii).

The implication (iv)$\Longrightarrow$(ii) is proved similarly.
\end{proof}

\begin{remark}
It is easy to see that (i)$\iff$(iv) without the additional
assumption that \eqref{local2} is an isomorphism.
\end{remark}

The following basic notion was introduced by Taylor \cite{T2};
cf. also \cite{GL} for a purely algebraic version.

\begin{definition}
\label{def:loc}
A homomorphism $\theta\colon A\to B$ of $\Ptens$-algebras is
a {\em localization}\footnote[1]{In Taylor's paper \cite{T2},
such homomorphisms are called {\em absolute localizations},
whereas the term ``localization'' is used for a somewhat wider
class of homomorphisms.} if it satisfies the conditions of
Proposition~\ref{prop:loc}.
In this case, we say (following \cite{NR}) that $B$ is
{\em stably flat} over $A$.
\end{definition}

\begin{remark}
Using condition (iv) of Proposition~\ref{prop:loc}, we see that
{\em $\theta\colon A\to B$ is a localization if and only if the functor
$B\ptens{A}(\,\cdot\,)\ptens{A}B\colon A\bimod A\to B\bimod B$ takes some
(=every) projective
bimodule resolution of $A$ to a projective bimodule resolution of $B$}.
This is exactly the definition given by Taylor \cite{T2}.
\end{remark}

\begin{prop}
\label{prop:loc_tens}
Suppose that $\theta\colon A\to B$ is a localization. Then for each
$M\in B\lmod$ the canonical map $B\ptens{A} M\to M,\; b\otimes x\mapsto b\cdot x$,
is an isomorphism.
\end{prop}
\begin{proof}
Apply the functor $(\,\cdot\,)\ptens{B} M$ to \eqref{local2}.
\end{proof}

A useful property of localizations is that they
``do not change homological relations between modules''.
In particular, if $A\to B$ is a localization,
then $\H^p(B,M)=\H^p(A,M)$ and $\H_p(B,M)=\H_p(A,M)$ for each
$B$-$\Ptens$-bimodule $M$ (see \cite{T2}, Prop.~1.4 and 1.7).
The next proposition is a combination of this fact with
the Cartan-Eilenberg inverse process \cite[XIII.5.1]{CE}.

\begin{prop}
\label{prop:Hoch_Lie}
Let $\fg$ be a finite-dimensional Lie algebra,
and let $U(\fg)$ be its universal
enveloping algebra endowed with the finest locally convex topology.
Suppose that $\theta\colon U(\fg)\to B$ is a localization.
For each $M\in B\bimod B$ denote by ${_E}M$ (resp. $M_E$) the left
(resp., right) $\fg$-module obtained from $M$ by setting
$X\cdot m=\theta(X)\cdot m-m\cdot\theta(X)$
(resp., $m\cdot X=m\cdot\theta(X)-\theta(X)\cdot m$); $X\in\fg,\; m\in M$.
Then there exist vector space isomorphisms
\begin{equation*}
\H^p(B,M)\cong H^p_{\Lie}(\fg,{_E}M),\quad
\H_p(B,M)\cong H_p^{\Lie}(\fg,M_E)\qquad
(p\in\Z).
\end{equation*}
\end{prop}

For later reference, we note the following
\begin{prop}[\cite{T2}]
\label{prop:loc_compos}
Let $A\xra{\theta} B \xra{\lambda} C$ be $\Ptens$-algebra homomorphisms.
Suppose $\theta$ is a localization. Then $\lambda$ is a localization if and only if
$\lambda\theta$ is a localization.
\end{prop}

By an {\em augmented $\Ptens$-algebra} we mean a $\Ptens$-algebra
$A$ together with a homomorphism $\varepsilon_A\colon A\to\CC$.
Homomorphisms of augmented $\Ptens$-algebras are defined in an obvious way.
Given an augmented $\Ptens$-algebra $A$, we consider $\CC$ as
an $A$-module via $\varepsilon_A$.

\begin{definition}
\label{def:wloc}
A homomorphism $\theta\colon A\to B$ of augmented $\Ptens$-algebras is
a {\em weak localization} if $B\pperp_A\CC$, and
the map
\begin{equation}
\label{w_local2}
B\ptens{A} \CC\to\CC,\qquad b\otimes\lambda\mapsto
\varepsilon(b)\lambda
\end{equation}
is a topological isomorphism.
\end{definition}

Setting $M=\CC$ in Proposition \ref{prop:loc_tens}, we get the following.

\begin{prop}
\label{prop:loc_wloc}
Each localization of augmented $\Ptens$-algebras is a weak localization.
\end{prop}

In the case of Hopf $\Ptens$-algebras with invertible antipodes,
the converse is also true.
To see this, let us first observe that if $\theta\colon U\to H$ is a homomorphism
of Hopf $\Ptens$-algebras, then the homomorphisms
$E_H\theta \colon U\to H^e$ and $(\theta\otimes\theta)E_U\colon U\to H^e$
coincide. Indeed,
\begin{multline*}
(\theta\otimes\theta)E_U=(\theta\otimes\theta)(\id_U\otimes S_U)\Delta_U
=(\theta\otimes\theta S_U)\Delta_U
=(\theta\otimes S_H\theta)\Delta_U\\
=(\id_H\otimes S_H)(\theta\otimes\theta)\Delta_U
=(\id_H\otimes S_H)\Delta_H\theta
=E_H\theta.
\end{multline*}
Hence any of the above homomorphisms can be used to make
$H^e$ into a right $U$-$\Ptens$-module. It also follows from the above that
the canonical isomorphisms
\begin{gather}
H^e_{E_H}\ptens{H} H_\theta\to H^e, \quad
x\otimes h\mapsto xE_H(h);\notag\\
\label{HeU}
H^e_{\theta\otimes\theta}\ptens{U^e} U^e_{E_U}\to H^e,\quad
x\otimes w\mapsto
x(\theta\otimes\theta)(w)
\end{gather}
are isomorphisms in $\rmod U$.

\begin{prop}
\label{prop:wloc_loc}
Let $\theta\colon U\to H$ be a homomorphism of Hopf $\Ptens$-algebras
with invertible antipodes.
Then $\theta$ is a localization if and only if it is a weak localization.
\end{prop}
\begin{proof}
The ``only if'' part readily follows from Proposition \ref{prop:loc_wloc}.
If $\theta$ is a weak localization, then the map $H\ptens{U}\CC\to\CC,\;
h\otimes\lambda\mapsto\varepsilon(h)\lambda$ is an isomorphism.
Combining this fact with Lemma \ref{lemma:tens_iso} and
\eqref{HeU}, we obtain a chain of isomorphisms
\begin{multline*}
H\ptens{U} H \lriso H\ptens{U} U \ptens{U} H
\lriso H^e_{\theta\otimes\theta}\ptens{U^e} U
\lriso H^e_{\theta\otimes\theta}\ptens{U^e} U^e_{E_U}\ptens{U}\CC\\
\lriso H^e\ptens{U}\CC
\lriso H^e_{E_H}\ptens{H} H \ptens{U}\CC
\lriso H^e_{E_H}\ptens{H}\CC
\lriso H.
\end{multline*}
It is easy to check that the composition of the above isomorphisms
takes each $h_1\otimes h_2\in H\ptens{U} H$
to $h_1 h_2\in H$. Thus we have shown that the canonical map
$H\ptens{U} H\to H$ is an isomorphism.

Now let $P_\bullet$ be a projective resolution of $\CC$ in $U\lmod$,
and let $\overline{P}_\bullet$ denote the augmented complex
$P_\bullet\to\CC\to 0$.
By Theorem \ref{thm:inverse}, the complex
$Q_\bullet=U^e_{E_U}\ptens{U} P_\bullet$
is a projective bimodule resolution of $U$.
In order to prove that $\theta$ is a localization, it remains to show that
the augmented tensor product complex
$H\ptens{U}\overline{Q}_\bullet\ptens{U}H
=H^e_{\theta\otimes\theta}\ptens{U^e}\overline{Q}_\bullet$
is admissible.

Since $\theta$ is a weak localization, we see that
$L_\bullet=H\ptens{U} P_\bullet$ is a projective resolution of $H\ptens{U}\CC\cong\CC$
in $H\lmod$. Using again Theorem \ref{thm:inverse}, we conclude that
$H^e_{E_H}\ptens{H} L_\bullet$ is a projective bimodule resolution of $H$.
In particular, the augmented complex $H^e_{E_H}\ptens{H}\overline{L}_\bullet$
is admissible. Now it follows from \eqref{HeU} that
$$
H^e_{E_H}\ptens{H}\overline{L}_\bullet
\cong H^e_{E_H}\ptens{H} H\ptens{U}\overline{P}_\bullet
\cong H^e\ptens{U}\overline{P}_\bullet
\cong H^e_{\theta\otimes\theta}\ptens{U^e} U^e_{E_U}\ptens{U}\overline{P}_\bullet
\cong H^e_{\theta\otimes\theta}\ptens{U^e}\overline{Q}_\bullet.
$$
Therefore $H^e_{\theta\otimes\theta}\ptens{U^e}\overline{Q}_\bullet$
is admissible, as required.
\end{proof}

We end this section with the following simple observation.
\begin{lemma}
\label{lemma:dense_loc}
Let $\theta\colon A\to B$ be a homomorphism of $\Ptens$-algebras
(resp. of augmented $\Ptens$-algebras) with dense range.
Then $\theta$ is a localization (resp. weak localization) if and only if
$B^e\pperp_{A^e} B$ (resp. $B\pperp_A\CC$).
\end{lemma}
\begin{proof}
Since $\Im\theta$ is dense in $B$, the map $X\ptens{A} Y
\to X\ptens{B} Y,\; x\tens{A}y\mapsto x\tens{B}y$ is a topological
isomorphism for each $X\in \rmod B$ and each $Y\in B\lmod$.
In particular, \eqref{local2} (resp. \eqref{w_local2}) is a topological
isomorphism. The rest is clear.
\end{proof}

\section{Localizations of $U(\fg)$ and duality}
\label{sect:loc_dual}

Following \cite{BFGP} (cf. also \cite{Litv_Hopf}),
we say that a Hopf $\Ptens$-algebra is
{\em well-behaved} if its underlying locally convex space is
either a nuclear Fr\'echet space or a nuclear (DF)-space.
Recall (see, e.g., \cite{Groth}) that the strong dual of a nuclear
Fr\'echet space is a complete nuclear (DF)-space, and vice versa.
Moreover, if $E$ is either a nuclear Fr\'echet space or a
complete nuclear (DF)-space, then there is a canonical topological
isomorphism $E'\Ptens E'\cong (E\Ptens E)'$. Therefore for each
well-behaved Hopf $\Ptens$-algebra $H$ the strong dual, $H'$, is also
a well-behaved Hopf $\Ptens$-algebra in a natural way.
More precisely, the multiplication (resp. comultiplication) on $H'$
is the dual of the comultiplication (resp. multiplication) on $H$,
the antipode of $H'$ is the dual of that of $H$, etc.
Note that $H$ is commutative (resp. cocommutative) if and only if
$H'$ is cocommutative (resp. commutative).
For example, if $G$ is a real Lie group, then the algebra
$C^\infty(G)$ of smooth functions is a nuclear commutative
Fr\'echet Hopf algebra,
and its dual is the Hopf algebra $\E'(G)$ of compactly supported distributions.
For later reference, recall that the comultiplication, the counit, and the
antipode of $C^\infty(G)$ are given, respectively, by
\begin{equation}
\label{grp_Hopf}
(\Delta f)(x,y)=f(xy),\quad \eps(f)=f(e), \quad (Sf)(x)=f(x^{-1}).
\end{equation}
Here we identify $C^\infty(G)\Ptens C^\infty(G)$ with $C^\infty(G\times G)$
(see, e.g., \cite{Groth}, Chap.~II, \S 3, no.~3).

Another important example is $U(\fg)$, the universal enveloping algebra
of a finite-dimensional Lie algebra $\fg$. If we endow $U(\fg)$ with the
finest locally convex topology, then it becomes a cocommutative nuclear
(DF) Hopf $\Ptens$-algebra. Recall that the comultiplication,
the counit, and the antipode of $U(\fg)$ are uniquely determined by
$$
\Delta(X)=X\otimes 1+1\otimes X,\quad
\eps(X)=0,\quad
S(X)=-X\quad (X\in\fg).
$$
The strong dual of $U(\fg)$ is topologically isomorphic to the Fr\'echet algebra
of formal power series $\CC[[z_1,\ldots ,z_n]]$ with the topology of convergence
of each coefficient (cf. \cite[Prop.~2.7.5]{Dixmier} and
\cite[Theorem 22.1]{Treves}; cf. also Lemma~\ref{lemma:dual_[U]} below).

Many other examples of well-behaved Hopf $\Ptens$-algebras
can be found in \cite{Litv}, \cite{Litv_Hopf},
\cite{Litv_hypergr}, \cite{BFGP}, and \cite{Pfl-Schot}.

Let $\fg$ be a Lie algebra. Suppose we are given a
Hopf $\Ptens$-algebra homomorphism
of $U(\fg)$ to a well-behaved Hopf $\Ptens$-algebra $H$.
In this section we formulate some conditions on the dual algebra,
$H'$, that are sufficient for $H$ to be stably flat over $U(\fg)$.

\subsection{Homotopy of commutative $\Ptens$-algebras}
\label{subsect:hmt}
In this subsection we briefly discuss a continuous version of
the notion of homotopy between morphisms of commutative
algebras. This notion was introduced by Chen \cite{Chen_hmt} in the purely
algebraic case. All the definitions and the results in this subsection are
straightforward adaptations of \cite{Chen_hmt} to the $\Ptens$-case.

Throughout this subsection, all $\Ptens$-algebras are assumed
to be commutative.

\begin{definition}[cf. \cite{Chen_hmt}, \cite{Chen_cov}]
A $\Ptens$-algebra $C$ is called {\em exact} if it possesses at least
one nonzero augmentation $C\to\CC$, and if there exists a derivation
$\partial\colon C\to N$ to some $C$-$\Ptens$-module $N$ such that the
sequence
$$
0 \rar \CC \xra{\eta_C} C \xra{\partial} N \rar 0
$$
splits in $\mathbf{LCS}$. A derivation $\partial$ with the above property
is called {\em split exact}.
\end{definition}

Basic examples of exact algebras are the algebra of smooth functions
$C^\infty(I)$ on an interval $I\subset\R$,
the algebra $\O(U)$ of holomorphic functions
on a simply connected domain $U\subset\CC$, the algebra $\CC[[z]]$ of formal power series,
the polynomial algebra $\CC[z]$ (with the finest locally convex topology), etc.
In each of the above examples, the usual derivation
$\frac{d}{dz}\colon A\to A$ is split exact.

\begin{definition}[cf. \cite{Chen_hmt}]
\label{def:hmt}
Two morphisms $\varphi_0,\varphi_1\colon A\to B$ of commutative
$\Ptens$-algebras are said to be {\em homotopic} if there exists an exact
algebra $C$, a morphism $\varPhi\colon A\to C\Ptens B$ and two
augmentations $\eps_0,\eps_1\colon C\to\CC$ such that
$$
\varphi_i=(\eps_i\otimes\id_B)\,\varPhi\quad (i=0,1).
$$
\end{definition}

For instance, if $X$ and $Y$ are smooth manifolds and $f_0, f_1\colon X\to Y$
are smooth homotopic mappings, then the induced morphisms
$f_0^*, f_1^*\colon C^\infty(Y)\to C^\infty(X)$ are homotopic in the above
sense. To see this, it suffices to set $C=C^\infty[0,1]$ and to reverse the
arrows in the usual definition of a (smooth) homotopy between $f_0$ and $f_1$.

Two $\Ptens$-algebras $A,B$ are called {\em homotopy equivalent} if there
exist morphisms $\varphi\colon A\to B$ and $\psi\colon B\to A$ such that
$\psi\varphi$ is homotopic to $\id_A$ and $\varphi\psi$ is homotopic to $\id_B$.
A $\Ptens$-algebra is said to be {\em contractible}\footnote[1]{We use the
word ``contractible'' following Chen \cite{Chen_hmt}; it should be noted, however,
that the notion of ``contractible algebra'' has an absolutely different
meaning in the cohomology theory of locally convex algebras (see, e.g., \cite{X1}).}
if it is homotopy equivalent to $\CC$. Equivalently, $A$ is contractible iff
there exist an exact algebra $C$, a morphism $\varPhi\colon A\to C\Ptens A$,
and augmentations $\eps_0,\eps_1\colon C\to\CC$ and $\eps_A\colon A\to\CC$
such that $(\eps_1\otimes\id_A)\,\varPhi=\id_A$
and $(\eps_0\otimes\id_A)\,\varPhi=\eta_A\eps_A$.
For example, the algebra of smooth functions on a contractible smooth
manifold is contractible. It is also easy to prove that the polynomial
algebra $\CC[z_1,\ldots ,z_n]$, the algebra of formal power series
$\CC[[z_1,\ldots ,z_n]]$, the algebra of entire functions $\O(\CC^n)$ etc.
are contractible.

\begin{theorem}[\cite{Chen_hmt}]
If two morphisms $\varphi_0,\varphi_1\colon A\to B$ of commutative
$\Ptens$-algebras are homotopic,
then the induced morphisms
$\varphi_{0,*}, \varphi_{1,*}\colon\Omega(A)\to\Omega(B)$
are chain homotopic (as morphisms of complexes in $\mathbf{LCS}$).
\end{theorem}

We omit the proof, because it is an obvious modification of the proof
from \cite{Chen_hmt} to the $\Ptens$-case.

\begin{corollary}
\label{cor:split_Omega}
If $A$ is a contractible $\Ptens$-algebra, then the augmented de Rham complex
$0 \to \CC \xra{\eta_A} \Omega(A)$ splits in $\mathbf{LCS}$.
\end{corollary}

\subsection{Lie algebra actions and parallelizability}
\label{subsect:parall}
\begin{definition}
Let $A$ be a commutative $\Ptens$-algebra, and let $\fg$ be a Lie algebra
acting on $A$ by derivations. We say that $A$ is {\em $\fg$-parallelizable} if
the derivation
$$
d^0\colon A\to C^1(\fg,A),\quad a\mapsto (X\mapsto Xa)
$$
is universal, i.e., if $(C^1(\fg,A),d^0)$ is the module of K\"ahler differentials for $A$.
\end{definition}

\begin{prop}
\label{prop:par_crit}
$A$ is $\fg$-parallelizable if and only if the identity map of $A$
extends to a DG $\Ptens$-algebra isomorphism between
$C^\cdot(\fg,A)$ and $\Omega(A)$.
\end{prop}
\begin{proof}
The ``if'' part is clear. To prove the converse, recall that the universal
property of $\Omega(A)$ yields a unique DG $\Ptens$-algebra morphism
$\varphi\colon\Omega(A)\to C^\cdot(\fg,A)$ such that $\varphi^0=\id_A$.
If $A$ is $\fg$-parallelizable, then
$\varphi^1\colon\Omega^1 A\to C^1(\fg,A)$ is an isomorphism.
Since both $\Omega(A)$ and $C^\cdot(\fg,A)$ are exterior,
we conclude that $\varphi$ is an isomorphism (see Subsection \ref{subsect:DGA}).
\end{proof}

Now suppose that $H$ is a well-behaved cocommutative Hopf $\Ptens$-algebra,
$\fg$ is a Lie algebra,
and $\theta\colon U(\fg)\to H$ is a Hopf $\Ptens$-algebra homomorphism.
We consider $H$ as a right $\fg$-$\Ptens$-module via $\theta$ by setting
$x\cdot X=x\theta(X)$ for each $x\in H,\; X\in\fg$.
The strong dual space, $H'$, is then a left $\fg$-$\Ptens$-module
in a natural way (see Subsection~\ref{subsect:Lie_act}).
Namely, the action of $\fg$ on $H'$ is given by
\begin{equation}
\label{act_theta}
\langle X\cdot a,x\rangle = \langle a,x\theta(X)\rangle\quad
(a\in H',\; X\in\fg,\; x\in H).
\end{equation}
It is easy to check that $\fg$ acts
on $H'$ by derivations. Indeed, for each $a,b\in H',\; X\in\fg,\; x\in H$ we obtain
\begin{multline}
\label{der_act}
\langle X\cdot ab,x\rangle
=\langle ab,x\theta(X)\rangle
=\langle a\otimes b,\Delta(x\theta(X))\rangle
=\langle a\otimes b,\Delta(x)\Delta(\theta(X))\rangle\\
=\langle a\otimes b,\Delta(x)\cdot\theta\otimes\theta(\Delta(X))\rangle
=\langle a\otimes b,\Delta(x)(\theta(X)\otimes 1+1\otimes\theta(X))\rangle.
\end{multline}
For each $x_1,x_2\in H$ we have
\begin{multline*}
\langle a\otimes b,(x_1\otimes x_2)(\theta(X)\otimes 1)\rangle
=\langle a\otimes b,x_1\theta(X)\otimes x_2\rangle\\
=\langle Xa,x_1\rangle\langle b,x_2\rangle
=\langle Xa\otimes b,x_1\otimes x_2\rangle.
\end{multline*}
Therefore $\langle a\otimes b,u(\theta(X)\otimes 1)\rangle
=\langle Xa\otimes b,u\rangle$ for each $u\in H\Ptens H$.
Similarly, $\langle a\otimes b,u(1\otimes\theta(X))\rangle
=\langle a\otimes Xb,u\rangle$ for each $u\in H\Ptens H$.
Setting $u=\Delta(x)$ and substituting in \eqref{der_act},
we see that
$$
\langle X\cdot ab,x\rangle
=\langle Xa\otimes b+a\otimes Xb,\Delta(x)\rangle
=\langle Xa\cdot b+a\cdot Xb,x\rangle.
$$
Hence $\fg$ acts on $H'$ by derivations.

In what follows, we say that the action defined by \eqref{act_theta}
is {\em determined
by $\theta$}. We shall sometimes refer to $\theta$ explicitly
by writing $X\cdot_\theta a$ instead of $X\cdot a$ or $Xa$.

\begin{theorem}
\label{thm:par_contr_loc}
Let $\fg$ be a Lie algebra, and let
$H$ be a well-behaved cocommutative Hopf $\Ptens$-algebra.
Suppose $\theta\colon U(\fg)\to H$ is a Hopf $\Ptens$-algebra homomorphism
with dense range.
Assume that $H'$ is $\fg$-parallelizable (w.r.t. the action determined by $\theta$)
and contractible. Then $\theta$ is a localization.
\end{theorem}
\begin{proof}
Since $H$ is cocommutative, we have $S^2=\id_H$ (for a categorical
proof of this classical fact, see \cite{Street}, Chap.~9).
In particular, $S$ is invertible.
In view of Proposition \ref{prop:wloc_loc}, it suffices to show that
$\theta$ is a weak localization. Set $U=U(\fg)$, and consider the
Koszul resolution
$$
0 \lar \CC \xla{\eps_U} V_\cdot(\fg)
$$
of the trivial $\fg$-module $\CC$
(see Subsection \ref{subsect:Lie_act}).
Clearly, the chain complexes
$H\ptens{U} V_\cdot(\fg)$ and $C_\cdot(\fg,H)$ are isomorphic.
Due to Lemma~\ref{lemma:dense_loc}, we need only
check that the augmented complex
\begin{equation*}
0 \lar \CC \xla{\eps_H} C_\cdot(\fg,H)
\end{equation*}
splits in $\mathbf{LCS}$.
Since the above complex consists of reflexive spaces,
it splits if and only if the dual complex
$$
0 \rar \CC \xra{\eta_{H'}} C^\cdot(\fg,H')
$$
splits. Now it remains to apply Corollary \ref{cor:split_Omega}
and Proposition \ref{prop:par_crit}.
\end{proof}

\section{Power series envelopes of $U(\fg)$}
\label{sect:power_env}
Our next task is to show that the strong dual algebras of some locally convex
completions of $U(\fg)$ (for $\fg$ nilpotent) are indeed $\fg$-parallelizable.
To this end, recall some facts on the ``formal
power series completion'' $[U(\fg)]$ of $U(\fg)$ (see \cite{Good_filt}).

Let $\fg$ be a nilpotent Lie algebra and let $I\subset U(\fg)$ be the ideal
generated by $\fg$. Recall that the quotient algebra $U(\fg)/I^n$ is
finite-dimensional for each $n$ (see, e.g., \cite[2.5.1]{Dixmier}).
Endow each $U(\fg)/I^n$ with the usual locally convex topology of a finite-dimensional
vector space, and set $[U(\fg)]=\varprojlim U(\fg)/I^n$. Clearly, $[U(\fg)]$
is a nuclear Fr\'echet-Arens-Michael algebra. We have a canonical homomorphism
\begin{equation}
\label{U->[U]}
\theta\colon U(\fg)\to [U(\fg)],\quad x\mapsto (x+I^n).
\end{equation}
Since $\fg$ is nilpotent, it follows that $\bigcap_n I^n=\{ 0\}$
(see, e.g., \cite[XIV.4.1]{Hochbook}),
so that \eqref{U->[U]} is injective. For notational convenience, we shall
often write $U$ instead of $U(\fg)$ and $[U]$ instead of $[U(\fg)]$, and we shall
identify $U$ with its canonical image in $[U]$.

It is easy to show that $[U(\fg)]$ has a natural structure of Hopf $\Ptens$-algebra
such that \eqref{U->[U]} is a Hopf algebra homomorphism.
Indeed, let $K=I\Tens U+U\Tens I\subset U\Tens U$
be the augmentation ideal of $U\Tens U$. Evidently, we have $\Delta(I)\subset K$,
and so $\Delta(I^n)\subset K^n$ for each $n$.
Therefore we obtain an algebra homomorphism
\begin{equation}
\label{U->lim}
[U]=\varprojlim U/I^n \to \varprojlim (U\Tens U)/K^n,\quad
(x+I^n)\mapsto (\Delta(x)+K^n).
\end{equation}
Since $U\Tens U$ is isomorphic to $U(\fg\times\fg)$, and since $\fg\times\fg$
is nilpotent together with $\fg$, it follows that
$\dim (U\Tens U)/K^n <\infty$ for each $n$. Hence we can
endow $\varprojlim (U\Tens U)/K^n$ with a locally convex topology in
the same way as we did for $[U]$. Thus \eqref{U->lim} becomes a $\Ptens$-algebra
homomorphism.

For each $n$ denote by $\tau_n\colon U\to U/I^n$ the quotient map, and
set
$$
\pi_n={\tau_n\otimes\tau_n}\colon U\Tens U\to (U/I^n)\Tens (U/I^n).
$$
We clearly have
$$
K^n=\sum_{i+j=n} I^i\Tens I^j,
$$
and so $K^{2n}\subset\Ker\pi_n$. Hence there exists a homomorphism
$$
(U\Tens U)/K^{2n} \to (U/I^n)\Ptens (U/I^n),\quad
y+K^{2n}\mapsto \pi_n(y).
$$
Taking the inverse limit and using the fact that the projective tensor
product commutes with reduced inverse limits \cite[41.6]{Kothe},
we obtain a homomorphism
\begin{equation*}
\varprojlim (U\Tens U)/K^n \to [U]\Ptens [U].
\end{equation*}
(This is even an isomorphism, since $\Ker\pi_n\subset K^n$ for each $n$.)
Composing with \eqref{U->lim}, we get a $\Ptens$-algebra
homomorphism
$$
[\Delta] \colon [U] \to [U]\Ptens [U].
$$
It is easy to check that $[\Delta]$ extends $\Delta$ in the sense
that the diagram
$$
\xymatrix{
[U] \ar[r]^(.4){[\Delta]} & [U]\Ptens [U]\\
U \ar[r]^(.4){\Delta} \ar[u]^\theta & U\Ptens U \ar[u]_{\theta\otimes\theta}
}
$$
is commutative. Since $\theta$ has dense range, the coassociativity
of $[\Delta]$ readily follows from that of $\Delta$.

Arguing as above, it is easy to construct an antipode $[S]\colon [U]\to [U]$
and a counit $[\eps]\colon [U]\to\CC$ in such a way that $[U]$ becomes a Hopf $\Ptens$-algebra
and $\theta\colon U\to [U]$ becomes a Hopf $\Ptens$-algebra homomorphism.

A somewhat more explicit construction of $[U]$ was suggested by Goodman \cite{Good_filt}.
Fix a positive filtration $\mathscr F$ on $\fg$, i.e., a decreasing chain of subspaces
$$
\fg=\fg_1\supset\fg_2\supset\cdots\supset\fg_\ell\supset\fg_{\ell+1}=0,\quad
[\fg_i,\fg_j]\subset\fg_{i+j}.
$$
The smallest $\ell$ such that $\fg_{\ell+1}=0$ is called the {\em length}
of the filtration.

An example of a positive filtration is the lower central series of $\fg$ defined
inductively by $\fg_{n+1}=[\fg,\fg_n]$.

Given $X\in\fg$, $X\ne 0$, the
{\em $\mathscr F$-weight} of $X$ is defined by $w(X)=\max\{ n: X\in\fg_n\}$.
A basis $(e_i)$ of $\fg$ is called an {\em $\mathscr F$-basis} if
$w(e_i)\le w(e_{i+1})$ for all $i$, and $\fg_n=\spn\{ e_i: w(e_i)\ge n\}$ for all $n$.
Given an $\mathscr F$-basis $(e_i)$, we set $w_i=w(e_i)$ for each $i$.
For each multi-index $\alpha=(\alpha_1,\ldots ,\alpha_N)\in\Z_+^N$ ($N=\dim\fg$),
set $|\alpha|=\sum_i \alpha_i$ and $w(\alpha)=\sum_i w_i\alpha_i$.
By the Poincar\'e-Birkhoff-Witt theorem, the elements
$e^\alpha=e_1^{\alpha_1}\cdots e_N^{\alpha_N}$ form a basis of $U(\fg)$.
For each $n$, set
\begin{equation}
\label{J_n}
J_n=\spn\{ e^\alpha : w(\alpha)\ge n\}\subset U(\fg).
\end{equation}
Then we have
\begin{equation}
\label{filt_U}
U(\fg)=J_0\supset I=J_1\supset J_2\supset\cdots,\quad
J_i J_j\subset J_{i+j},
\end{equation}
so that $\{ J_n\}$ is a decreasing filtration on $U(\fg)$ satisfying
$\bigcap_n J_n=\{ 0\}$. In particular, each $J_n$ is an ideal of $U(\fg)$.
Goodman \cite{Good_filt} defines $[U(\fg)]_{\mathscr F}$ as the completion
of $U(\fg)$ w.r.t. the topology determined by the filtration $\{ J_n\}$.
Thus we have an algebraic isomorphism
$[U(\fg)]_{\mathscr F}=\varprojlim U(\fg)/J_n$.

If we endow each $U(\fg)/J_n$ with the usual topology of a finite-dimensional
vector space, then it is easily seen that $[U(\fg)]_{\mathscr F}$ is isomorphic
(as a topological algebra) to the algebra $[U(\fg)]$ introduced above.
Indeed,  setting $C=\max w_i$, we see that $w(\alpha)\le C|\alpha|$ for each
$\alpha\in\Z_+^N$, and so $J_{Cn}\subset I^n$ for each $n$.
On the other hand, we have $I^n=J_1^n\subset J_n$ for each $n$.
Therefore the filtrations $\{ J_n\}$ and $\{ I^n\}$ are equivalent, and so the
algebras $[U(\fg)]_{\mathscr F}$ and $[U(\fg)]$ are isomorphic.

As a locally convex space, $[U(\fg)]$ is isomorphic to the space
of all formal power series $x=\sum_\alpha c_\alpha e^\alpha$ endowed with
the topology of convergence of each coefficient (cf. \cite{Good_filt}).
More precisely, the topology on $[U(\fg)]$ can be generated by the sequence
of seminorms $\{\|\cdot\|_n: n\in\Z_+\}$ defined by
\begin{equation*}
\| x\|_n=\sum_{w(\alpha)\le n} |c_\alpha|\quad\text{for each}\quad
x=\sum_\alpha c_\alpha e^\alpha\in [U].
\end{equation*}

For each multi-index $\alpha\in\Z_+^N$ set $e_\alpha=e^\alpha/\alpha!$.
Then $\Delta(e_\gamma)=\sum_{\alpha+\beta=\gamma} e_\alpha\otimes e_\beta$
(see \cite[2.7.2]{Dixmier}), and the same relation clearly holds for $[\Delta]$.

\begin{lemma}
\label{lemma:dual_[U]}
The mapping $\varkappa\colon [U]'\to \CC[z_1,\ldots ,z_N]$ defined by the rule
\begin{equation}
\label{dual_[U]}
f\mapsto \sum_\alpha f(e_\alpha) z^\alpha
\end{equation}
is an algebra isomorphism. Moreover, $\varkappa$ is a topological isomorphism
w.r.t. the strong topology on $[U]'$ and the finest l.c.
topology on $\CC[z_1,\ldots ,z_N]$.
\end{lemma}
\begin{proof}
The continuity of $f$ implies that there exists $n\in\Z_+$ such that
$f(e_\alpha)=0$ whenever $w(\alpha)\ge n$. Hence the sum in the right-hand
side of \eqref{dual_[U]} is finite, and $\varkappa$ is well defined.
Since the $e_\alpha$'s generate a dense subspace of $[U]$, we see that
$\varkappa$ is injective. Conversely, for every polynomial
$p=\sum_\alpha \lambda_\alpha z^\alpha$ the mapping $f\colon [U]\to\CC,\;
f(\sum_\alpha c_\alpha e_\alpha)=\sum_\alpha c_\alpha \lambda_\alpha$ is a
continuous linear functional on $[U]$ satisfying $\varkappa(f)=p$.
Hence $\varkappa$ is bijective. A direct computation (see \cite[2.7.5]{Dixmier})
shows that $\varkappa$ is an algebra homomorphism.
Finally, since the topology on the strong dual of a countable inverse limit
of finite-dimensional spaces is the finest l.c. topology (see, e.g.,
\cite[Theorem 22.1]{Treves}),
we see that $\varkappa$ is a topological isomorphism.
\end{proof}

\begin{definition}
\label{def:power_env}
Let $\fg$ be a nilpotent Lie algebra. By a {\em power series
envelope} of $U(\fg)$ we mean a Hopf $\Ptens$-algebra $H$ together with
Hopf $\Ptens$-algebra homomorphisms $\theta_1\colon U(\fg)\to H$ and
$\theta_2\colon H\to [U(\fg)]$ such that both $\theta_1$ and
$\theta_2$ are injective with dense ranges,
and the composition
$$
U(\fg) \xra{\theta_1} H \xra{\theta_2} [U(\fg)]
$$
coincides with the canonical homomorphism $\theta$ defined by
\eqref{U->[U]}.
\end{definition}

\begin{remark}
Since $\theta$ is injective with dense range, the conditions
``$\theta_1$ is injective'' and ``$\theta_2$ has dense range'' are satisfied
automatically. Note also that, since $U(\fg)$ is cocommutative and $\theta_1$
has dense range, $H$ is also cocommutative. For the same reason, we have
$S^2=\id_H$ in $H$.
\end{remark}

It is immediate from the definition that the ``smallest'' power series
envelope of $U(\fg)$ is $U(\fg)$ itself, and the ``largest'' one is $[U(\fg)]$.

\begin{theorem}
\label{thm:conv_prl}
Let $\fg$ be a nilpotent Lie algebra, and let
$H$ be a well-behaved power series envelope of $U(\fg)$.
Then $H'$ is $\fg$-parallelizable.
\end{theorem}
\begin{proof}
Fix a positive filtration $\mathscr F$
on $\fg$, and choose an $\mathscr F$-basis $(e_i)$ of $\fg$.
Using Lemma \ref{lemma:dual_[U]}, we may identify $[U]'$
and $\CC[z_1,\ldots ,z_N]$. For each $i=1,\ldots ,N$
set $x_i=\theta_2'(z_i)\in H'$. Since $\theta_2$ is injective, it follows that
$\Im\theta_2'$ is dense in $H'$ w.r.t. the weak$^*$ topology $\sigma(H',H)$.
Using the semireflexivity of $H$, we see that $\Im\theta_2'$ is dense in $H'$ w.r.t.
the strong topology as well. Hence $x_1,\ldots ,x_N$ generate
a dense subalgebra of $H'$.

Set $A=H'$, and consider the free $A$-module $A^N$
with the standard $A$-basis $(u_i)$,
i.e., $u_i=(0,\ldots ,1,\ldots ,0)$
with $1$ in the $i$th coordinate, $0$ elsewhere. Denote by $(e^i)\subset\fg^*$ the basis
dual to $(e_i)$ (i.e., $e^i(e_j)=\delta_{ij}$ for all $i,j$).
Identifying the $A$-modules $C^1(\fg,A)$ and $A\Tens\fg^*$, we see
that the elements $v_i=1\otimes e^i$
($i=1,\ldots ,N$) form an $A$-basis of $C^1(\fg,A)$.

Now consider the $A$-module morphism $\varphi\colon A^N\to C^1(\fg,A)$ taking
each $u_i$ to $d^0(x_i)$. Let $(\varphi_{ij})$ be the matrix of $\varphi$
w.r.t. the bases $(u_i)$ and $(v_i)$, respectively. Applying the identity
$\varphi(u_j)=\sum_i\varphi_{ij} v_i$ to $e_i$, we see that
$$
\varphi_{ij}=\varphi(u_j)(e_i)=d^0(x_j)(e_i)=e_i\cdot x_j.
$$
Given $a\in H'$, denote by $\bar a$ the restriction of $a$ to $U(\fg)$
(i.e., $\bar a=\theta_1'(a)$). Then for each $y\in U(\fg)$ we have
\begin{equation}
\label{ij}
\langle \bar\varphi_{ij} , y\rangle
=\langle \varphi_{ij} , \theta_1(y)\rangle
=\langle e_i\cdot x_j,\theta_1(y)\rangle
=\langle  x_j,\theta_1(ye_i)\rangle
=\langle \bar x_j ,ye_i \rangle.
\end{equation}
We claim that the matrix $(\varphi_{ij})$ is upper triangular with $1$'s
on the main diagonal. Indeed, using \eqref{filt_U}, we see that
$ye_i\in J_{w_i}$ for each $y\in U(\fg)$; moreover, $ye_i\in J_{w_i+1}$
for each $y\in I=J_1$. On the other hand,
it is immediate from \eqref{dual_[U]} that
$\bar x_j=z_j|_{U(\fg)}$ vanishes on $J_{w_j+1}$. Hence
$\langle \bar x_j,ye_i\rangle = 0$ for all $i>j$ and all $y\in U(\fg)$,
$\langle \bar x_i,ye_i\rangle = 0$ for all $y\in I$,
and $\langle \bar x_i,e_i\rangle = 1$.
Together with \eqref{ij}, this gives $\bar\varphi_{ij}=0$ for each $i>j$,
and $\bar\varphi_{ii}=1$.
Finally, since $\Im\theta_1$ is dense in $H$, it follows that $\theta_1'$
is injective, and so the latter relations hold with $\bar\varphi_{ij}$
replaced by $\varphi_{ij}$.
Therefore the matrix of $\varphi$ has the required form,
so that $\varphi$ is an isomorphism.

For each $i=1,\ldots ,N$, let $p_i\colon A^N\to A$ be the projection on the
$i$th direct summand. Evidently, $\partial_i=p_i\varphi^{-1} d^0$
is a derivation of $A$. It is immediate from the definition of $\varphi$
that $\partial_i(x_j)=p_i(u_j)=\delta_{ij}$ for each $i,j$. Hence the conditions
of Lemma \ref{lemma:Omega_free} are satisfied, and so
$\partial=(\partial_1,\ldots ,\partial_N)\colon A\to A^N$ is a universal derivation.
Since $\varphi$ is an isomorphism, we conclude that
$d^0=\varphi\partial$ is a universal derivation as well,
i.e., $A$ is $\fg$-parallelizable.
\end{proof}

Combining the above theorem with Theorem \ref{thm:par_contr_loc},
we obtain the following.

\begin{corollary}
\label{cor:conv_loc}
Let $\fg$ be a nilpotent Lie algebra, and let
$H$ be a well-behaved power series envelope of $U(\fg)$
such that $H'$ is contractible.
Then $H$ is stably flat over $U(\fg)$.
\end{corollary}

\begin{corollary}
\label{cor:[U]_loc}
For each nilpotent Lie algebra $\fg$,
$[U(\fg)]$ is stably flat over $U(\fg)$.
\end{corollary}
\begin{proof}
By Lemma \ref{lemma:dual_[U]}, the algebra dual to $[U(\fg)]$ is
isomorphic to $\CC[z_1,\ldots ,z_N]$ and hence is contractible.
Now it remains to apply Corollary~\ref{cor:conv_loc}.
\end{proof}

\section{Arens-Michael envelopes of universal enveloping algebras}

In this section we prove that for each positively graded,
finite-dimensional Lie algebra $\fg$ the Arens-Michael envelope
of $U(\fg)$ is stably flat over $U(\fg)$.
First we recall some facts on Arens-Michael envelopes.

\subsection{Arens-Michael envelopes}
Arens-Michael envelopes of topological algebras (under a different name)
were introduced by Taylor (\cite{T1}, Definition 5.1).
Here we follow the terminology of Helemskii's book \cite{X2}.

\begin{definition}[\cite{X2}, Chap. V]
\label{def:AM}
Let $A$ be a topological algebra. A pair $(\wh{A},\iota_A)$ consisting of
an Arens-Michael algebra $\wh{A}$ and a continuous homomorphism
$\iota_A\colon A\to\wh{A}$ is called {\em the Arens-Michael envelope} of $A$
if for each Arens-Michael algebra $B$ and for each continuous
homomorphism $\varphi\colon A\to B$ there exists a unique
continuous homomorphism $\wh{\varphi}\colon\wh{A}\to B$ making
the following diagram commutative:
\begin{equation*}
\xymatrix{
\wh{A} \ar@{-->}[r]^{\wh{\varphi}} & B \\
A \ar[u]^{\iota_A} \ar[ur]_\varphi
}
\end{equation*}

In the above situation, we say that {\em $\wh{\varphi}$ extends $\varphi$}
(though $\iota_A$ is not injective in general; see Remark \ref{rem:AM=0} below).
\end{definition}

\begin{remark}
In the above definition, it suffices to consider only homomorphisms
with values in Banach algebras. This is immediate from the fact that
each Arens-Michael algebra is an inverse limit of Banach algebras
(see, e.g., \cite{X2}, Chap. V).
\end{remark}

Clearly, the Arens-Michael envelope is unique in the sense that
if $(\wh{A},\iota_A)$ and $(\overline{A},j_A)$ are Arens-Michael envelopes
of $A$, then there exists a unique isomorphism $j\colon\wh{A}\to\overline{A}$
of topological algebras such that the following diagram is commutative:
$$
\xymatrix{
\wh{A} \ar@{-->}[rr]^j && \overline{A} \\
& A \ar[ul]^{\iota_A} \ar[ur]_{j_A}
}
$$

Recall (see \cite{T1} and \cite[Chap. V]{X2}) that the Arens-Michael envelope
of a topological algebra $A$ always exists and
can be obtained as the completion of $A$
w.r.t. the family of all continuous submultiplicative seminorms on $A$.
This implies, in particular, that $\iota_A\colon A\to\wh{A}$ has dense range.
It is easy to see that the correspondence $A\mapsto\wh{A}$ is a functor
from the category $\mathbf{TA}$ of topological algebras to the category
$\mathbf{AM}$ of Arens-Michael algebras. In what follows, we call it
{\em the Arens-Michael functor}. Clearly, the Arens-Michael functor
is the left adjoint to the forgetful
functor from $\mathbf{AM}$ to $\mathbf{TA}$.

If $A$ is not equipped with a topology, then by an Arens-Michael
envelope of $A$ we mean the Arens-Michael envelope of the finest
locally convex algebra $A_s$ (see Section \ref{sect:prelim}).

Here are two basic examples due to Taylor \cite{T2}.
\begin{example}
The Arens-Michael envelope of the polynomial algebra $\CC[z_1,\ldots ,z_n]$
is the algebra of entire functions $\O(\CC^n)$ endowed with the
compact-open topology.
\end{example}

\begin{example}
\label{example:env_free}
Let $F_n$ be the free $\CC$-algebra on $n$ generators
$\zeta_1,\ldots ,\zeta_n$. Given a $k$-tuple
$\alpha=(\alpha_1,\ldots ,\alpha_k)$ of integers from $[1,n]$, we set
$\zeta_\alpha=\zeta_{\alpha_1}\cdots\zeta_{\alpha_k}\in F_n$ and $|\alpha|=k$.
It is convenient to agree that the identity of $F_n$ corresponds
to the tuple of length zero ($k=0$).
The algebra $\mathscr F_n$ of ``free power series''
consists of all formal expressions
$a=\sum_\alpha \lambda_\alpha \zeta_\alpha$ satisfying the condition
$$
\| a\|_\rho=\sum_\alpha |\lambda_\alpha| \rho^{|\alpha|} <\infty
\quad\mbox{for all}\quad 0<\rho <\infty\, .
$$
The system of seminorms
$\{\|\cdot\|_\rho : 0<\rho < \infty\}$ makes $\mathscr F_n$ into a
Fr\'echet-Arens-Michael algebra.
Evidently, $F_n$ is a subalgebra of $\mathscr F_n$.
Taylor \cite{T2} proved that
$\mathscr F_n$ is the Arens-Michael envelope of $F_n$.
Note that in the case $n=1$ we have $F_1=\CC[z]$ and
$\mathscr F_1\cong\O(\CC)$.
\end{example}

\begin{remark}
\label{rem:AM=0}
It should be noted that the Arens-Michael envelope
can be trivial even in very simple cases.
For example, let $A$ be the Weyl algebra, i.e., the algebra with
two generators $p,q$ subject to the relation $[p,q]=1$.
It is a standard exercise from spectral theory (see, e.g., \cite{X2}, Prop.~2.1.21)
to show that $A$ has no nonzero submultiplicative seminorms.
Hence $\wh{A}=0$.

Another example of this kind is given in \cite{X2}, Chap. V.
\end{remark}

\begin{remark}
If $\fg$ is a finite-dimensional Lie algebra, then (in contrast to
the previous example) the homomorphism
$\iota_{U(\fg)}\colon U(\fg)\to\wh{U}(\fg)$ is injective. This readily
follows from the fact that finite-dimensional representations
(and, a fortiori, Banach space representations) of $\fg$ separate the points
of $U(\fg)$ (see, e.g., \cite[2.5.7]{Dixmier}).
\end{remark}

The next proposition shows that the Arens-Michael functor commutes with quotients.

\begin{prop}
\label{prop:env_quot}
Let $A$ be a topological algebra and $I$ a two-sided ideal of $A$.
Denote by $J$ the closure of $\iota_A(I)$ in $\wh{A}$.
Then $J$ is a two-sided ideal of $\wh{A}$, and the
homomorphism $A/I\to \wh{A}/J$ induced by $\iota_A\colon A\to\wh{A}$
extends to a topological algebra isomorphism
$$
\wh{A/I}\cong (\wh{A}/J)\sptilde\, .
$$
\end{prop}
\begin{proof}
Since $\iota_A$ has dense range, we see that $J$ is indeed an ideal of $\wh{A}$,
and so $(\wh{A}/J)^\sim$ is an Arens-Michael algebra.
Consider the homomorphism $\bar\iota\colon A/I\to (\wh{A}/J)^\sim$ taking
each $a+I\in A/I$ to $\iota_A(a)+J$.
We claim that $((\wh{A}/J)^\sim,\bar\iota)$ is the Arens-Michael
envelope of $A/I$. Indeed, each homomorphism $\varphi$ from $A/I$
to an Arens-Michael algebra $C$ determines a homomorphism $\psi\colon A\to C$
vanishing on $I$. Each such homomorphism extends to a homomorphism
$\wh{\psi}\colon \wh{A}\to C$ that vanishes on $J$ and hence gives a homomorphism
$\wh{\varphi}\colon (\wh{A}/J)^\sim\to C$. It is now elementary to check that
$\wh{\varphi} \bar\iota=\varphi$.
The uniqueness of $\wh{\varphi}$ is immediate from the fact that
$\bar\iota$ has dense range.
\end{proof}

Since each separated quotient of a Fr\'echet space
is complete, we obtain the following
\begin{corollary}
\label{cor:env_quot}
Under the conditions of Proposition \ref{prop:env_quot}, assume that
$\wh{A}$ is a Fr\'echet algebra. Then $\wh{A/I}\cong\wh{A}/J$.
\end{corollary}

\begin{corollary}
\label{cor:AM_fingen}
If $A$ is a finitely generated algebra, then $\wh{A_s}$ is a nuclear Fr\'echet
algebra.
\end{corollary}
\begin{proof}
Since $A$ is finitely generated, it is isomorphic to a quotient of the
free algebra $F_n$ for some $n$. By Corollary \ref{cor:env_quot},
$\wh{A}$ is isomorphic to a quotient of $\wh{F}_n=\mathscr F_n$
(see Example \ref{example:env_free}). Since $\mathscr F_n$ is a nuclear
Fr\'echet space \cite{Luminet_PI}, so is $\wh{A}$.
\end{proof}

Another useful property of the Arens-Michael functor is that it commutes
with projective tensor products.

\begin{prop}
\label{prop:AM_Ptens}
Let $A,B$ be $\Ptens$-algebras. Then there exists a topological algebra isomorphism
$$
(A\Ptens B)\sphat\,\cong\wh{A}\Ptens\wh{B}\, .
$$
\end{prop}
\begin{proof}
Set $\iota=\iota_A\otimes\iota_B\colon A\Ptens B\to\wh{A}\Ptens\wh{B}$.
Clearly, $\iota$ is a continuous homomorphism.
Suppose $\varphi\colon A\Ptens B\to C$ is a homomorphism to some
Arens-Michael algebra $C$. Then $\varphi_1\colon A\to C,\;
\varphi_1(a)=\varphi(a\otimes 1)$ and
$\varphi_2\colon B\to C,\;
\varphi_2(b)=\varphi(1\otimes b)$ extend to continuous homomorphisms
$\wh{\varphi}_1\colon\wh{A}\to C$ and
$\wh{\varphi}_2\colon\wh{B}\to C$, i.e., we have $\wh{\varphi}_1\iota_A=\varphi_1$
and $\wh{\varphi}_2\iota_B=\varphi_2$.
Let $\wh{\varphi}\colon\wh{A}\Ptens\wh{B}\to C$ be the linear continuous
map associated to the bilinear map $\wh{A}\times\wh{B}\to C,\;
(a,b)\mapsto\wh{\varphi}_1(a)\wh{\varphi}_2(b)$.
Evidently, we have $\wh{\varphi}\iota=\varphi$.
Since $\iota$ has dense range, we conclude that
$\wh{\varphi}$ is an algebra homomorphism. For the same reason, $\wh{\varphi}$
is a unique homomorphism extending $\varphi$. Hence
$(\wh{A}\Ptens\wh{B},\iota)$ is the Arens-Michael envelope of
$A\Ptens B$.
\end{proof}

\begin{prop}
\label{prop:AM_op}
Let $A$ be a topological algebra. Then $(A^{\op})\sphat\,\cong\wh{A}^{\op}$.
\end{prop}
\begin{proof}
It suffices to use the 1-1 correspondence between
continuous homomorphisms $A^{\op}\to B$
and continuous homomorphisms $A\to B^{\op}$.
\end{proof}

\begin{corollary}
Let $A$ be a $\Ptens$-algebra. Then $(A^e)\sphat\,\cong (\wh{A})^e$.
\end{corollary}

The next proposition shows that the Arens-Michael functor can also
be considered as a functor from the category $\mathbf{HTA_{\Ptens}}$ of
Hopf $\Ptens$-algebras to the category $\mathbf{HAM_{\Ptens}}$ of Hopf $\Ptens$-algebras
that are Arens-Michael algebras.

\begin{prop}
\label{prop:AM_Hopf}
Let $H$ be a Hopf $\Ptens$-algebra. Then
there exists a unique Hopf $\Ptens$-algebra structure on $\wh{H}$
such that $\iota_H\colon H\to\wh{H}$ becomes a Hopf $\Ptens$-algebra homomorphism.
Moreover, if $L$ is both a Hopf $\Ptens$-algebra
and an Arens-Michael algebra, and $\varphi\colon H\to L$ is a Hopf $\Ptens$-algebra
homomorphism, then so is $\wh{\varphi}\colon\wh{H}\to L$.
\end{prop}
\begin{proof}
To obtain $\Delta_{\wh{H}}$, $\varepsilon_{\wh{H}}$, and
$S_{\wh{H}}$, it suffices to apply the Arens-Michael functor
to $\Delta_H$, $\varepsilon_H$, and $S_H$, respectively, and to use
Propositions \ref{prop:AM_Ptens} and \ref{prop:AM_op}.
The Hopf algebra axioms (such as the coassociativity of
$\wh{\Delta}_{\wh{H}}$ etc.) are then readily verified by applying
the Arens-Michael functor to the appropriate
commutative diagrams involving $H$.

To prove that $\wh{\varphi}$ respects comultiplication, it is enough to
show that $(\wh{\varphi}\otimes\wh{\varphi}) \Delta_{\wh{H}}\iota_H=
\Delta_L\wh{\varphi}\iota_H$. We have
\begin{multline*}
(\wh{\varphi}\otimes\wh{\varphi}) \Delta_{\wh{H}}\iota_H=
(\wh{\varphi}\otimes\wh{\varphi})(\iota_H\otimes\iota_H) \Delta_H=
(\varphi\otimes\varphi) \Delta_H=
\Delta_L\varphi=
\Delta_L \wh{\varphi}\iota_H.
\end{multline*}
A similar argument shows that $\wh{\varphi} S_{\wh{H}}=S_L \wh{\varphi}$
and $\eps_L \wh{\varphi} = \eps_{\wh{H}}$. Hence $\wh{\varphi}$ is a Hopf
$\Ptens$-algebra homomorphism.
\end{proof}

\begin{example}
\label{example:AM_U(g)}
Let $\fg$ be a finite-dimensional Lie algebra and $U(\fg)$ the universal
enveloping algebra of $\fg$. Then it follows from Proposition \ref{prop:AM_Hopf}
and Corollary \ref{cor:AM_fingen} that $\wh{U}(\fg)$ is a well-behaved
(see the beginning of Section~\ref{sect:loc_dual}) Hopf $\Ptens$-algebra.
Denote by $\iota_{\fg}\colon\fg\to\wh{U}(\fg)$ the restriction of
$\iota_{U(\fg)}$ to $\fg$. Then it is easy to see that
$\wh{U}(\fg)$ is characterized by the following universal
property: {\em for each Arens-Michael algebra $A$ and each Lie algebra homomorphism
$\varphi\colon\fg\to A$ there exists a unique $\Ptens$-algebra homomorphism
$\psi\colon\wh{U}(\fg)\to A$ such that $\psi\iota_{\fg}=\varphi$}.
In particular, for each Lie algebra homomorphism
$f\colon\fg\to\fh$ there exists a unique $\Ptens$-algebra homomorphism
$\wh{U}(f)\colon\wh{U}(\fg)\to\wh{U}(\fh)$ such that $\wh{U}(f)\iota_{\fg}=\iota_{\fh} f$. Moreover,
Proposition \ref{prop:AM_Hopf} implies that $\wh{U}(f)$ is in fact
a Hopf $\Ptens$-algebra homomorphism
(cf. \cite{Bour_Lie-II}, Chap.~II, \S1, no.~4).
\end{example}

\subsection{Arens-Michael envelopes of filtered and graded algebras}
In this subsection we describe Arens-Michael envelopes of locally finite
graded algebras. As a corollary, we show that the Arens-Michael envelope
of the universal enveloping algebra of a nilpotent Lie algebra $\fg$
is a power series envelope (see Definition~\ref{def:power_env})
provided $\fg$ admits a positive grading.

Recall that a {\em decreasing filtration} on an algebra $A$
is a chain of linear subspaces
$$
A=A_0\supset A_1\supset A_2\supset\ldots\quad\text{satisfying}\quad
A_i A_j\subset A_{i+j}.
$$
The filtration is called {\em separated} if $\bigcap_n A_n=\{ 0\}$
and is said to be {\em of finite type} if $\dim A_n/A_{n+1} <\infty$
for all $n$. In the sequel all filtrations are assumed to have these properties.

As in Section \ref{sect:power_env}, we endow each $A/A_n$ with
the usual locally convex topology of a finite-dimensional vector space, and
set $[A]=\varprojlim A/A_n$.

The following proposition is immediate from the definition of $[A]$.
\begin{prop}
\label{prop:[A]_filt}
For each $n\in\Z_+$ let $V_n$ be a linear complement of $A_{n+1}$
in $A_n$. Fix a norm on each $V_n$. Then, as a locally convex
space,  $[A]$ is isomorphic to the space of all formal series
$\{ a=\sum_i v_i : v_i\in V_i\}$ endowed with the family of seminorms
$\{ \|\cdot\|_n : n\in\Z_+\}$ defined by
$$
\| a\|_n = \sum_{i=0}^n \| v_i\| \quad\text{for each}\quad
a=\sum_i v_i.
$$
\end{prop}

Since each $A/A_n$ is a finite-dimensional
(hence Banach) algebra, we see that $[A]$ is an Arens-Michael algebra.
Therefore the canonical homomorphism
\begin{equation}
\label{theta}
\theta\colon A\to [A],\quad x\mapsto (x+A_n)_{n\in\Z_+}
\end{equation}
uniquely extends to a homomorphism
\begin{equation}
\label{thetahat}
\wh{\theta}\colon\wh{A}\to [A],\quad
\wh{\theta}\iota_A=\theta.
\end{equation}

\begin{prop}
Let $A$ be an algebra. Suppose that $A$ admits a decreasing, separated
filtration of finite type. Then the canonical homomorphism
${\iota_A\colon A\to\wh{A}}$ is injective. In other words,
submultiplicative seminorms separate the points of $A$.
\end{prop}
\begin{proof}
The condition $\bigcap_n A_n=\{ 0\}$ implies that $\theta$ is injective.
Since $\theta=\wh{\theta}\iota_A$, we conclude that
$\iota_A$ is also injective.
\end{proof}

Our next task is to show that $\wh{\theta}\colon\wh{A}\to [A]$
is also injective provided
the filtration on $A$ comes from a grading.

Let $A=\bigoplus_{n\ge 0} A^n$ be a graded algebra (see Subsection~\ref{subsect:DGA}).
We assume that $A$ is {\em locally finite}, i.e., $\dim A^n<\infty$ for
each $n$. Setting $A_n=\bigoplus_{i\ge n} A^i$, we obtain a decreasing,
separated filtration of finite type on $A$.

The following is a direct consequence of Proposition \ref{prop:[A]_filt}.

\begin{prop}
\label{prop:[A]_grad}
Let $A=\bigoplus_{n\ge 0} A^n$ be a locally finite graded algebra.
Then, as a $\Ptens$-algebra, $[A]$ is isomorphic to the direct
product $\prod_n A^n$ endowed with the multiplication
\begin{equation}
\label{[A]_grad}
(a_i)\cdot (b_j) = (c_k),\quad c_k=\sum_{i+j=k} a_i b_j.
\end{equation}
\end{prop}

In order to describe the Arens-Michael envelope of $A$ as
a certain ``power series algebra'', it will be convenient to use
``vector-valued K\"othe spaces'', which are more or less straightforward
generalizations of classical K\"othe spaces (see, e.g., \cite{Pietsch}).

Let $E=\{ E_i : i\in\N\}$ be a countable family of Hausdorff locally
convex spaces. For each $i$ denote by $\mathfrak N(E_i)$ the set of all
continuous seminorms on $E_i$.

\begin{definition}
\label{def:power}
An {\em $E$-power set} is a family $P$ of functions
$p\colon\N\to\bigcup_i\mathfrak N(E_i)$ such that $p_i=p(i)\in\mathfrak N(E_i)$
for each $i$, and the following conditions are satisfied:
\begin{itemize}
\item[1)] for each $i\in\N$ the family of seminorms
$\{ p_i : p\in P\}$ generate the original topology on $E_i$;
\item[2)] for each $p,q\in P$ there exists $r\in P$ such that
$r_i(x)\ge\max\{ p_i(x),q_i(x)\}$ for each $i\in\N$ and each $x\in E_i$.
\end{itemize}
\end{definition}

\begin{definition}
Given a family $E=\{ E_i : i\in\N\}$ of Hausdorff l.c.s.'s and an $E$-power
set $P$, define the {\em vector-valued K\"othe space}
$\lambda(P,E)$ by
$$
\lambda(P,E)=\Bigl\{ x=(x_i)\in\prod_i E_i :
\| x\|_p=\sum_i p_i(x_i)<\infty\quad\forall p\in P\Bigr\}.
$$
\end{definition}

\begin{remark}
If $E_i=\CC$ for each $i$, then we come to the classical notion of
K\"othe sequence space.
\end{remark}

Evidently, $\lambda(P,E)$ is a Hausdorff locally convex space w.r.t. the family
of seminorms $\{\|\cdot\|_p : p\in P\}$.

\begin{prop}
$\lambda(P,E)$ is complete iff all the $E_i$'s are complete.
\end{prop}

We omit the proof, because it is a straightforward modification of the classical fact that
$\ell^1$ is complete.

\medskip
Now let $A=\bigoplus_{n\ge 0} A^n$ be a locally finite graded algebra.
As usual, we endow each $A^n$ with the usual topology of a finite-dimensional
vector space.
\begin{definition}
A graded submultiplicative seminorm on $A$ is
a function $p\colon\N\to\bigcup_n\mathfrak N(A^n)$ such that
$p_n=p(n)\in \mathfrak N(A^n)$ for all $n\in\Z_+$, and
$p_{i+j}(ab)\le p_i(a) p_j(b)$ for all $i,j\in\Z_+$ and all $a\in A^i,\; b\in A^j$.
\end{definition}
If $p$ is a graded submultiplicative seminorm on $A$, then the associated
seminorm $\|\cdot\|_p\colon A\to\R_+$ defined by
$\| a\|_p=\sum_i p_i(a_i)$ for each $a=\sum_i a_i,\; a_i\in A^i$
is submultiplicative in the usual sense.
Therefore graded submultiplicative seminorms on $A$ are in $1$-$1$ correspondence
with submultiplicative seminorms $\|\cdot\|$ on $A$ satisfying the condition
$\| a\|=\sum_i \| a_i\|$ for each $a=\sum_i a_i,\; a_i\in A^i$.

\medskip
Denote by $P$ the collection of all graded submultiplicative seminorms on $A$.

\begin{lemma}
$P$ is an $A$-power set.
\end{lemma}
\begin{proof}
To check condition 1) of Definition \ref{def:power}, it suffices to show that
for each $n$ there exists $p\in P$ such that $p_n$ is a norm on $A^n$.
Fix a submultiplicative norm on the finite-dimensional algebra $A/A_{n+1}$,
denote by $\tau_{n+1}\colon A\to A/A_{n+1}$ the quotient map, and
set $p_i(a)=\|\tau_{n+1}(a)\|$ for each $i\in\Z_+$ and each $a\in A^i$.
Evidently, $p$ is a graded submultiplicative seminorm on $A$.
Since $\|\cdot\|$ is a norm on $A/A_{n+1}$, and since
$A^n\cap\Ker\tau_{n+1}=\{ 0\}$, we conclude that $p_n$ is a norm on $A^n$.

Given $p,q\in P$, the function $r=\max\{ p,q\}$ (i.e., $r_i(a)=\max\{ p_i(a),q_i(a)\}$
for each $a\in A^i$ and each $i\in\Z_+$) clearly belongs to $P$.
Hence condition 2) of Definition
\ref{def:power} is also satisfied, so that $P$ is an $A$-power set.
\end{proof}

\begin{theorem}
\label{thm:AM_grad}
Let $A=\bigoplus_{n\ge 0} A^n$ be a locally finite graded algebra,
and let $P$ be the set of all graded submultiplicative seminorms on $A$.
Denote by $\iota_A$ the canonical embedding of $A$ into $\lambda(P,A)$
that is the identity on each $A^n$. Then $\lambda(P,A)$ is a subalgebra
of $\prod_n A^n$, and $(\lambda(P,A),\iota_A)$ is the Arens-Michael envelope of $A$.
\end{theorem}
\begin{proof}
Given $a=(a_i)$ and $b=(b_j)$ in $\lambda(P,A)$, we must show that
the element $c=ab\in\prod_n A^n$ defined by \eqref{[A]_grad}
belongs to $\lambda(P,A)$. For each $p\in P$ we have
\begin{align*}
\sum_k p_k(c_k)
\le\sum_k \sum_{i+j=k} p_k(a_i b_j)
& \le\sum_k\sum_{i+j=k} p_i(a_i)p_j(b_j)\\
&=\sum_i p_i(a_i) \sum_j p_j(b_j)
=\| a\|_p \| b\|_p.
\end{align*}
Hence $ab\in\lambda(P,A)$, and $\| ab\|_p\le\| a\|_p \| b\|_p$.
This implies, in particular, that $\lambda(P,A)$ is an Arens-Michael
algebra, and $\iota_A$ is an algebra homomorphism.

Now let $\varphi\colon A\to B$ be a homomorphism to some Arens-Michael
algebra $B$. Fix a submultiplicative seminorm $\|\cdot\|$ on $B$,
and define $p\colon\Z_+\to\bigcup_n\mathfrak N(A^n)$ by $p_i(a_i)=\|\varphi(a_i)\|$
for each $i\in\Z_+$ and each $a_i\in A^i$. Evidently, $p$ is a graded submultiplicative
seminorm on $A$, and $\|\varphi(a)\|\le \|\iota_A(a)\|_p$ for each $a\in A$.
This implies that $\varphi$ is continuous w.r.t. the topology induced on $A$ from
$\lambda(P,A)$. Since $A$ is dense in $\lambda(P,A)$, we see that there
exists a unique continuous homomorphism $\wh{\varphi}\colon\lambda(P,A)\to B$
extending $\varphi$. Hence $\lambda(P,A)$ is the Arens-Michael envelope of $A$.
\end{proof}

\begin{corollary}
\label{cor:inj_grad}
Let $A=\bigoplus_{n\ge 0} A^n$ be a locally finite graded algebra,
and let $\theta\colon A\to [A]$ be the canonical
homomorphism \eqref{theta}. Then the induced homomorphism
$\wh{\theta}\colon\wh{A}\to [A]$ (see \eqref{thetahat}) is injective.
\end{corollary}
\begin{proof}
If we identify $[A]$ with $\prod_n A^n$ via Proposition \ref{prop:[A]_grad}
and $\wh{A}$ with $\lambda(P,A)$ via Theorem \ref{thm:AM_grad},
then $\wh{\theta}$ becomes the natural inclusion of $\lambda(P,A)$
into $\prod_n A^n$.
\end{proof}

Now let $\fg=\bigoplus_{n=1}^\ell\fg^n$ be a positively graded,
finite-dimensional Lie algebra. As in the case of associative algebras (see above),
we may define a filtration $\mathscr F=\{\fg_n\}$ on $\fg$ by setting
$\fg_n=\bigoplus_{i\ge n}\fg^i$. It is easy to show that the universal enveloping
algebra $U(\fg)$ has a grading $U(\fg)=\bigoplus_{n\ge 0} U(\fg)^n$
such that the associated filtration on $U(\fg)$ coincides with \eqref{J_n}.
Indeed, let $T\fg$ be the tensor algebra of $\fg$, and let
$L$ be the two-sided ideal of $T\fg$ generated by elements of the
form $x\otimes y-y\otimes x-[x,y];\; x,y\in\fg$.
Then we have $U(\fg)\cong T\fg /L$. If $\fg$ is graded, then we can define
a grading on $T\fg$ by
$$
(T\fg)^n=\bigoplus_{i_1+\cdots +i_k=n} \fg^{i_1}\otimes\cdots\otimes\fg^{i_k}.
$$
Thus $T\fg$ becomes a locally finite graded algebra, and $L$ becomes a graded
ideal of $T\fg$. Therefore $U(\fg)=T\fg /L$ is also a locally finite graded algebra.
We have
$$
U(\fg)^n=\sum_{i_1+\cdots +i_k=n} \fg^{i_1}\ldots\fg^{i_k}.
$$
Choose an $\mathscr F$-basis $(e_i)$ of $\fg$ consisting of homogeneous elements,
and set
$$
V^n=\spn\{ e^\alpha : w(\alpha)=n\}.
$$
By the Poincar\'e-Birkhoff-Witt theorem, we have $U(\fg)=\bigoplus_n V^n$.
On the other hand, it is clear that $V^n\subset U(\fg)^n$.
Since $U(\fg)=\bigoplus_n U(\fg)^n$, we conclude that $V^n=U(\fg)^n$ for all $n$.
Hence the associated filtration of $U(\fg)$ has the form
$$
U(\fg)_n=\bigoplus_{m\ge n} U(\fg)^m=\spn\{ e^\alpha : w(\alpha)\ge n\}=J_n
$$
(see \eqref{J_n}).

Now Proposition \ref{prop:AM_Hopf}, Example \ref{example:AM_U(g)},
and Corollary \ref{cor:inj_grad} imply the following.
\begin{prop}
\label{prop:U_conv}
Let $\fg$ be a positively graded Lie algebra. Then $\wh{U}(\fg)$
together with the homomorphisms $\iota_{U(\fg)}\colon U(\fg) \to \wh{U}(\fg)$
and $\wh{\theta}\colon\wh{U}(\fg)\to [U(\fg)]$ is a power series
envelope of $U(\fg)$.
\end{prop}

\subsection{The contractibility of $\wh{U}'(\fg)$}

Let $\fg$ be a positively graded Lie algebra.
In order to prove that $\wh{U}(\fg)$ is stably flat over $U(\fg)$,
it now remains to show that the strong dual, $\wh{U}'(\fg)$, of
$\wh{U}(\fg)$ is a contractible $\Ptens$-algebra
(see Corollary \ref{cor:conv_loc}).
To this end, it will be convenient to use the following Lie algebra
version of contractibility.

\begin{definition}
\label{def:contr_Lie}
We say that a finite-dimensional Lie algebra $\fg$ is {\em contractible} if there
exists a smooth mapping $h\colon [0,1]\times\fg\to\fg$ such that
\begin{itemize}
\item[(i)] for each $t\in [0,1]$ the map $h_t\colon\fg\to\fg,\;
h_t(X)=h(t,X)$ is a Lie algebra homomorphism;
\item[(ii)] $h_0=0$ and $h_1=\id_{\fg}$.
\end{itemize}
\end{definition}

\begin{example}
Each positively graded Lie algebra $\fg=\fg_1\oplus\cdots\oplus\fg_\ell$
is contractible. To see this, it suffices to set $h_t(X)=t^n X$ for each $X\in\fg_n$
and each $t\in [0,1]$.
\end{example}

\begin{example}
Let $\fg$ be the $2$-dimensional Lie algebra with basis $X,Y$ and commutation
relation $[X,Y]=Y$. Take a function $f\in C^\infty(\R)$ such that $f(t)=0$ for each $t\le 0$
and $f(t)=1$ for each $t\ge 1$, and define $h_t\colon\fg\to\fg$ by
$h_t(X)=f(2t)X$ and $h_t(Y)=f(2t-1)Y$. It is easy to check that $h_t$ satisfies
the conditions of Definition \ref{def:contr_Lie}, and so $\fg$ is contractible.
\end{example}

\begin{remark}
It is easy to prove that each contractible Lie algebra is solvable.
Indeed, suppose that $\fg$ is not solvable, and consider the Levi decomposition
$\fg=\mathfrak r\oplus\mathfrak l$ ($\mathfrak r=\mathrm{rad}\,\fg$,
$\mathfrak l$ is a semisimple subalgebra, $\mathfrak l\ne 0$).
It is clear that a semidirect summand (i.e., a retract in the category
of Lie algebras) of a contractible Lie algebra is contractible.
Thus it suffices to show that $\mathfrak l$ is not contractible.
Since $\mathfrak l$ is a direct sum of simple algebras, we need only
prove that a simple Lie algebra is not contractible. Assume towards a
contradiction that $\fg$ is both simple and contractible, and let
$h_t\colon\fg\to\fg$ be a contracting homotopy from Definition~\ref{def:contr_Lie}.
Since $\fg$ is simple, each $h_t$ is either $0$ or an automorphism.
Replacing, if necessary, the segment $[0,1]$ by $[t_0,1]$ where
$t_0=\max\{ t : h_t=0\}$, we may assume that $h_t$ is an automorphism
for all $t>0$. Let $B(\cdot , \cdot)$ denote the Killing form on $\fg$.
By Cartan's criterion, $B$ is nondegenerate.
Since $h_t$ is an automorphism, we have $B(h_t(X),h_t(Y))=B(X,Y)$ for
all $X,Y\in\fg$ and all $t>0$. Letting $t\to 0$, we obtain $B\equiv 0$,
which is a contradiction.
\end{remark}

\begin{remark}
It should be noted that not every nilpotent Lie algebra is contractible.
For example, let $\fg$ be the $7$-dimensional Lie algebra with basis
$X_1,\ldots ,X_7$ and commutation relations
\begin{gather*}
[X_1,X_i]=X_{i+1}\quad (i=2,\ldots ,6),\\
[X_2,X_3]=-X_6,\; [X_3,X_4]=X_7,\; [X_2,X_4]=[X_2,X_5]=-X_7
\end{gather*}
(see \cite{Favre}). Let $\{ h_t : t\in [0,1]\}$ be a continuous family
of endomorphisms of $\fg$, and let $(h_{ij}(t))$ be the matrix of $h_t$
w.r.t. the basis $X_1,\ldots ,X_7$. A routine calculation shows that
if $h_{11}(t_0)\ne 0$ and $h_{22}(t_0)\ne 0$ at some point $t_0$, then
$h_{ii}(t_0)=1$ for all $i=1,\ldots ,7$.
This clearly implies that $\fg$ is not contractible.
\end{remark}

Our next goal is to prove that the contractibility of $\fg$ implies
that of $\wh{U}'(\fg)$. We need some facts on topological vector
spaces. Most of them are standard and
can be easily deduced from \cite{Sch} and \cite{Groth}.

Let $E,\; F$, and $G$ be locally convex spaces (l.c.s's).
Consider the vector space
$\mathfrak B(E\times F,G)$ of all separately continuous bilinear
mappings from $E\times F$ to $G$. We endow this space with the topology of
\textit{bibounded convergence} (i.e., the topology of uniform convergence
on direct products of bounded sets).
There is a natural mapping
\begin{equation}
\label{freeze}
\L(E,\L(F,G))\to\mathfrak B(E\times F,G)
\end{equation}
defined by the rule $\varphi\mapsto ((x,y)\mapsto \varphi(x)(y))$.
Obviously, this mapping is topologically injective.
A bilinear map $\varPhi\colon E\times F\to G$ belongs to the image of the
mapping \eqref{freeze} iff for each $0$-neighborhood $U\subset G$ and each
bounded set $B\subset F$ there exists a $0$-neighborhood
$V\subset E$ such that $\varPhi(V\times B)\subset U$.
Such bilinear maps are usually called {\itshape $F$-hypocontinuous}.
If $E$ is barreled, then each separately continuous map
of $E\times F$ to $G$ is $F$-hypocontinuous (see \cite[III.5.2]{Sch}),
so the mapping \eqref{freeze} is surjective in this case. Therefore,
{\em for each
barreled l.c.s. $E$ and arbitrary l.c.s.'s $F$ and $G$ we have
a topological isomorphism}
\begin{equation}
\label{freeze2}
\L(E,\L(F,G))\cong\mathfrak B(E\times F,G).
\end{equation}

Recall also (see \cite{Groth}, Chapitre~II, Th\'eor\`eme~6 or \cite[IV.9.4]{Sch})
that {\em for each complete barreled nuclear l.c.s. $E$ and each complete l.c.s.
$F$ there exists a natural topological isomorphism
\begin{equation}
\label{nucl_tens}
E\Ptens F\to \L(E',F)
\end{equation}
defined by $x\otimes y\mapsto (x'\mapsto \langle x,x'\rangle y)$.}

\begin{lemma}
\label{lemma:duality}
Let $E$ be either a nuclear Fr\'echet space or a complete nuclear $(DF)$-space,
and let $F$ be a complete nuclear barreled l.c.s. Then for each complete
l.c.s. $G$ there exists a topological isomorphism
\begin{equation}
\label{tvsiso}
\L(E,F\Ptens G) \lriso \L(F',E'\Ptens G)
\end{equation}
taking each $u\colon E\to F\Ptens G$ to $v\colon F'\to E'\Ptens G$
such that
\begin{equation}
\label{tvsiso2}
\langle v(y'), x\otimes z'\rangle = \langle u(x), y'\otimes z'\rangle\, ,
\end{equation}
for each $x\in E,\; y'\in F',\; z'\in G'$.
\end{lemma}
\begin{proof}
Applying \eqref{nucl_tens} and \eqref{freeze2},
we obtain topological isomorphisms
\begin{equation}
\label{tvs1}
\L(E,F\Ptens G)\cong \L(E,\L(F',G))
\cong \mathfrak B(E\times F', G) \cong
\mathfrak B(F'\times E, G).
\end{equation}
Since $F$ is complete and nuclear, it is semireflexive \cite[IV.5]{Sch},
and hence $F'$ is barreled (see \cite{Sch}, IV.5.5).
Further, the assumptions on $E$ imply that $E$ is reflexive,
and $E'$ is barreled and nuclear \cite{Groth}.
Using again \eqref{freeze2} and \eqref{nucl_tens},
we see that
\begin{align}
\mathfrak B(F'\times E, G) &\cong
\L(F',\L(E,G))\notag\\
\label{tvs2}
&\cong\L(F',\L(E'',G))\cong
\L(F',E'\Ptens G).
\end{align}
Combining \eqref{tvs1} and \eqref{tvs2}, we obtain the required
isomorphism \eqref{tvsiso}. Relation \eqref{tvsiso2} is then readily
verified.
\end{proof}

Recall that for each smooth manifold $M$ and each
complete l.c.s. $X$ there exists a topological
isomorphism $C^\infty(M)\Ptens X\cong C^\infty(M,X)$
taking an elementary tensor $f\otimes x$ to the function $t\mapsto f(t)x$
(see \cite{Groth}, Chap.~II, \S 3, no.~3). Applying the previous lemma
to $G=C^\infty(M)$, we obtain the following.

\begin{corollary}
\label{cor:duality}
Let $E$ and $F$ be locally convex spaces satisfying the conditions
of Lemma \ref{lemma:duality}, and let $M$ be a smooth manifold.
Then there exists a topological isomorphism
$$
\L(E,C^\infty(M,F))\cong\L(F',C^\infty(M,E'))
$$
taking each $u\colon E\to C^\infty(M,F)$ to
$v\colon F'\to C^\infty(M,E')$ such that
$$
\langle v(y')(t), x\rangle = \langle y', u(x)(t)\rangle\, ,
$$
for each $x\in E,\; y'\in F',\; t\in M$.
\end{corollary}

\begin{theorem}
\label{thm:U'_contr}
Let $\fg$ be a contractible, finite-dimensional Lie algebra.
Then $\wh{U}'(\fg)$ is contractible as a commutative $\Ptens$-algebra.
\end{theorem}
\begin{proof}
Set $I=[0,1]$ and suppose that
$h\colon I\times\fg\to\fg$ is a smooth map satisfying the conditions
of Definition \ref{def:contr_Lie}. Note that the space $C^\infty(I,\fg)$
is a Lie algebra w.r.t. the pointwise multiplication. It is readily seen that
the map
$$
F\colon\fg\to C^\infty(I,\fg),\quad F(X)(t)=h(t,X)
$$
is a Lie algebra homomorphism. Using the universal property
of $\wh{U}=\wh{U}(\fg)$ (see Example \ref{example:AM_U(g)}) and the obvious fact that
$C^\infty(I,\wh{U})$ is an Arens-Michael algebra, we obtain a unique
continuous homomorphism $\varPsi\colon\wh{U}\to C^\infty(I,\wh{U})$
that fits into the commutative diagram
$$
\xymatrix{
\wh{U} \ar[r]^(.4)\varPsi & C^\infty(I,\wh{U})\\
\fg \ar[u]^{\iota_{\fg}} \ar[r]_(.4)F
& C^\infty(I,\fg) \ar[u]_{C^\infty(I,\iota_{\fg})}
}
$$
For each $t\in I$ define $\psi_t\colon\wh{U}\to\wh{U}$ by
$\psi_t(x)=\varPsi(x)(t)$. Evidently, $\psi_t$ is an algebra homomorphism.
Using the above diagram, it is readily seen that $\psi_t$
extends $h_t$ in the sense that $\psi_t\iota_{\fg}=\iota_{\fg} h_t$.
Hence $\psi_t=\wh{U}(h_t)$ (see Example \ref{example:AM_U(g)}), and
so $\psi_t$ is a Hopf $\Ptens$-algebra homomorphism.
Since $h_1=\id_\fg$ and $h_0=0$, we see that $\psi_1=\id_{\wh{U}}$
and $\psi_0=\eta_{\wh{U}}\eps_{\wh{U}}$.

Now set $A=\wh{U}'$ and let $\varPhi\colon A\to C^\infty(I,A)$ be the map
corresponding to $\varPsi$ under the isomorphism
$$
\L(\wh{U},C^\infty(I,\wh{U}))\cong \L(A,C^\infty(I,A))
$$
(see Corollary \ref{cor:duality}). For each $a\in A$, $t\in I$, and $x\in\wh{U}$
we have
$$
\langle \varPhi(a)(t), x\rangle = \langle a,\varPsi(x)(t)\rangle=
\langle a,\psi_t(x)\rangle = \langle a\circ\psi_t,x\rangle.
$$
Hence $\varPhi(a)(t)=a\circ\psi_t$. In other words, for each $t\in I$ the map
$\varphi_t\colon A\to A$ defined by $\varphi_t(a)=\varPhi(a)(t)$
is the dual of $\psi_t$.
Since $\psi_t$ is a $\Ptens$-coalgebra homomorphism, we conclude that
$\varphi_t$ is a $\Ptens$-algebra homomorphism. Hence so is $\varPhi$.
Note also that $\varphi_1=\id_A$ and $\varphi_0=(\eta_{\wh{U}}\eps_{\wh{U}})'=
\eta_A\eps_A$.

Now it is easy to check that $\varPhi\colon A\to C^\infty(I,A)\cong C^\infty(I)\Ptens A$
yields a homotopy between $\id_A$ and $\eta_A\eps_A$
(see Definition~\ref{def:hmt}). Indeed, consider the
augmentations $\eps_i\colon C^\infty(I)\to\CC,\; \eps_i(f)=f(i)$ ($i=0,1$).
Then for each $a\in A$ we have
\begin{gather*}
\bigl((\eps_0\otimes\id_A)\,\varPhi\bigr)(a)=\varPhi(a)(0)=\varphi_0(a)=(\eta_A\eps_A)(a),\\
\bigl((\eps_1\otimes\id_A)\,\varPhi\bigr)(a)=\varPhi(a)(1)=\varphi_1(a)=a.
\end{gather*}
Hence $(\eps_0\otimes\id_A)\,\varPhi=\eta_A\eps_A$, and
$(\eps_1\otimes\id_A)\,\varPhi=\id_A$.
Therefore $\id_A$ is homotopic to $\eta_A\eps_A$, i.e., $A$ is contractible.
\end{proof}

Now, applying Corollary \ref{cor:conv_loc}, Proposition \ref{prop:U_conv},
and Theorem \ref{thm:U'_contr}, we obtain the following.

\begin{theorem}
\label{thm:AM_loc}
Let $\fg$ be a finite-dimensional, positively graded Lie algebra.
Then $\wh{U}(\fg)$ is stably flat over $U(\fg)$.
\end{theorem}

We end this section with an application of the above theorem
to computing injective homological dimensions of $\wh{U}(\fg)$-modules.
To this end, we need a formula of ``Poincar\'e duality''
type. Let $\fg$ be a Lie algebra of dimension $n$.
Recall (see, e.g., \cite[6.10]{Knapp}) that for each left $\fg$-module $V$
there exist vector space isomorphisms
\begin{equation*}
H^p_{\Lie}(\fg,V)\cong H_{n-p}^{\Lie}(\fg,V\Tens(\textstyle\bigwedge^n\fg)^*)\qquad
(p\in\Z).
\end{equation*}
If $\fg$ is nilpotent, then it is easily seen that the action of $\fg$
on $\bigwedge^n\fg$ is trivial. (To see this, it suffices to take a basis
$(e_i)$ of $\fg$ with the property that
$[e_i,e_j]\in\spn\{ e_k : k\ge\max\{i,j\}\}$
and to observe that each $e_i$ acts on $e_1\wedge\ldots\wedge e_n$ trivially.)
Therefore the above formula takes the form
\begin{equation*}
H^p_{\Lie}(\fg,V)\cong H_{n-p}^{\Lie}(\fg,V)\qquad
(p\in\Z).
\end{equation*}
Combining this with Proposition~\ref{prop:Hoch_Lie}, we obtain the following.

\begin{corollary}
\label{cor:poinc_U(g)}
Let $\fg$ be a finite-dimensional, positively graded Lie algebra,
and let $n=\dim\fg$. Then for each $M\in\wh{U}(\fg)\bimod\wh{U}(\fg)$
there exist vector space isomorphisms
\begin{equation*}
\H^p(\wh{U}(\fg),M)\cong\H_{n-p}(\wh{U}(\fg),M)\qquad (p\in\Z).
\end{equation*}
\end{corollary}

\begin{corollary}
\label{cor:injdim}
Let $\fg$ be a finite-dimensional, positively graded Lie algebra.
Then
\begin{itemize}
\item[{\upshape (i)}]
$\injdh_{\wh{U}(\fg)} M=\dim\fg$ for each $M\in\wh{U}(\fg)\lmod$,
$M\ne 0$;
\item[{\upshape (ii)}]
$\dh_{\wh{U}(\fg)} M=\dim\fg$ for each Banach $M\in\wh{U}(\fg)\lmod$,
$M\ne 0$.
\end{itemize}
In particular, there are no nonzero injective
$\wh{U}(\fg)$-$\Ptens$-modules.
\end{corollary}
\begin{proof}
This is an immediate consequence of \cite[Theorem 2.1]{Pir_inj2},
\cite[Corollary 4.1.3]{Pir_Nova},
and Corollary~\ref{cor:poinc_U(g)}.
\end{proof}

\section{Weighted completions of universal enveloping algebras}
\label{sect:weight}
In this section we describe one more class of Fr\'echet Hopf algebras
that are stably flat completions of universal enveloping algebras. These algebras
were introduced by Goodman in \cite{Good_filt} and \cite{Good_hreps}.
Each of them is a power series envelope of $U(\fg)$
(see Definition \ref{def:power_env}) and consists of power series
$x\in [U]$ subject to certain growth conditions.

Recall some definitions and notation from \cite{Good_filt} and \cite{Good_hreps}.
Let $\fg$ be a nilpotent Lie algebra, and let $N=\dim\fg$.
Choose a positive filtration $\F$ on $\fg$,
and fix an $\F$-basis $(e_i)$ for $\fg$ (see Section \ref{sect:power_env}).
A sequence $\M=\{ M_\alpha:
\alpha\in\Z_+^N\}$ of positive numbers
is an {\em $\F$-weight sequence} if $M_0=1$ and $M_\gamma\le M_\alpha M_\beta$
whenever $w(\gamma)\ge w(\alpha)+w(\beta)$. Given an $\F$-weight sequence $\M$,
consider the space
$$
U(\fg)_\M=\Bigl\{ x=\sum_\alpha c_\alpha e^\alpha\in [U] :
\| x\|_r=\sum_\alpha |c_\alpha| \alpha! M_\alpha r^{w(\alpha)}<\infty
\;\forall r>0\Bigr\}.
$$
Clearly, $U(\fg)_\M$ is a Fr\'echet space w.r.t. the topology defined by the
family of seminorms $\{\|\cdot\|_r : r>0\}$. Using the Grothendieck-Pietsch criterion
(see, e.g., \cite{Pietsch}), it is easy to see that $U(\fg)_\M$ is nuclear.
Goodman \cite{Good_filt} proved that
$U(\fg)_\M$ is a subalgebra of $[U]$, and the multiplication in $U(\fg)_\M$
is (jointly) continuous w.r.t. the above topology. Note, however, that $U(\fg)_\M$
need not be an Arens-Michael algebra.

\begin{example}
\label{example:M_fact}
Let $\fg$ be an abelian Lie algebra endowed with the trivial filtration $\F$
(i.e., $\F_1=\fg$ and $\F_2=0$), and let $M_\alpha=|\alpha|^{-|\alpha|}$.
Then it is easy to see that $U(\fg)_\M$ is isomorphic to the algebra
$\O(\CC^N)$ of entire functions on $\CC^N$. Indeed, $\O(\CC^N)$ is topologized
by the family of seminorms $\|\cdot\|'_r\; (r>0)$ defined by
$\| f\|'_r=\sum_\alpha |c_\alpha| r^{|\alpha|}$ for each
$f(z)=\sum_\alpha c_\alpha z^\alpha\in\O(\CC^N)$. We clearly have
$\alpha!\le |\alpha|^{|\alpha|}$ for each $\alpha\in\Z_+^N$. On the other hand,
the Cauchy estimates applied to the entire function $z\mapsto\exp(\sum_i z_i)$
imply that $|\alpha|^{|\alpha|}\le C^{|\alpha|} \alpha!$ for some constant $C>0$.
Since $w(\alpha)=|\alpha|$ in this case, we obtain
$\| f\|'_{r/C} \le \| f\|_r \le \| f\|'_r$ for all polynomials $f$
and all $r>0$. Hence the families of seminorms $\|\cdot\|_r$ and $\|\cdot\|'_r$
are equivalent, and so $U(\fg)_\M$ and $\O(\CC^N)$ are isomorphic.
\end{example}

\begin{example}
\label{example:M_1}
Let $\fg$ be an abelian Lie algebra endowed with the trivial filtration,
and let $M_\alpha=1$ for all $\alpha$. Then $U(\fg)_\M$ is isomorphic to the
algebra of entire functions on $\CC^N$ of exponential order $\le 1$
and minimal type (cf. \cite{Rash}).
\end{example}

An $\F$-weight sequence $\M$ is {\em entire} \cite{Good_hreps} if it satisfies
the following two conditions:
\begin{gather}
\label{entire}
\sum_\alpha M_\alpha r^{w(\alpha)} <\infty\quad\text{for all } r>0;\\
\sup_{\alpha,\beta\ne 0} \{ A^{w(\alpha)/w(\beta)} M_\beta^{1/w(\beta)}
M_\alpha^{-1/w(\alpha)}\} <\infty
\quad\text{for some } A>0.\notag
\end{gather}

For instance, the weight sequence of Example \ref{example:M_fact}
is entire \cite{Good_hreps}, while that of Example \ref{example:M_1} is not entire.

If $\M$ is an entire $\F$-weight sequence,
then the dual of $U(\fg)_\M$ admits an explicit description
as a certain function algebra \cite{Good_hreps}.
Namely, let $G$ be the connected, simply connected complex Lie group
corresponding to $\fg$. Since $\fg$ is nilpotent,
the exponential map $\exp\colon\fg\to G$
is biholomorphic. The {\em homogeneous norm} on $G$ is defined
by
$$
|g|=\max_i |t_i|^{1/w_i}\quad\text{for each } g=\exp(\sum_i t_i e_i)\in G.
$$
Given $z\in\CC$, define a linear map $\delta_z\colon\fg\to\fg$ by
$\delta_z(e_i)=z^{w_i}e_i$. We use the same symbol $\delta_z$ to denote
the corresponding holomorphic self-map of $G$ satisfying $\delta_z\circ\exp=
\exp\circ\delta_z$. It is immediate that $\delta_1=\id_G$,
$\delta_0(g)=e$ for all $g\in G$ (here $e$ is the identity of $G$),
$\delta_z\delta_{z'}=\delta_{zz'}$,
$\delta_z^{-1}=\delta_{z^{-1}}$ for each $z\ne 0$,
and $|\delta_z g|=|z||g|$ for each $z\in\CC,\; g\in G$.

Given an entire $\F$-weight sequence $\M$, define the weight function
$W_\M$ on $G$ by
$$
W_\M(g)=\sum_\alpha M_\alpha |g|^{w(\alpha)}.
$$
Condition \eqref{entire} implies that $W_\M$ is finite on $G$.
For example, if $\fg$ is abelian and $M_\alpha=|\alpha|^{-|\alpha|}$
(see Example \ref{example:M_fact}), then $W_\M$ satisfies the estimate
\begin{equation}
\label{exp_grth}
\exp(N|g|/C)\le W_\M(g) \le \exp(N|g|).
\end{equation}

Given $r>0$, consider the space
$$
A_{\M,r}(G)=\Bigl\{ f\in\O(G) :
N_r(f)=\sup_{g\in G}\frac{|f(g)|}{W_\M(\delta_r g)}<\infty\Bigr\}.
$$
Evidently, $A_{\M,r}(G)$ is a Banach space w.r.t. the norm $N_r$.
Note that $W_\M(\delta_s g)\le W_\M(\delta_r g)$ whenever $0\le s\le r$.
This implies that $A_{\M,s}(G)\subset A_{\M,r}(G)$ for each $s\le r$, and
$N_r(f)\le N_s(f)$ for each $f\in A_{\M,s}(G)$. Therefore we may consider
the locally convex space
$$
A_\M(G)=\varinjlim A_{\M,r}(G).
$$
Goodman \cite{Good_hreps} proved that $A_\M(G)$ is a subalgebra
of $\O(G)$ (under pointwise multiplication), and the multiplication is
jointly continuous w.r.t. the inductive limit topology on $A_\M(G)$.

For example, if $\fg$ is abelian and $M_\alpha=|\alpha|^{-|\alpha|}$
(Example \ref{example:M_fact}), then it follows from \eqref{exp_grth} that
$A_\M(G)$ is the algebra of entire functions on $G=\fg$ of exponential
order $\le 1$.

Denote by $\P(G)$ the algebra of polynomial functions on $G$
(i.e., functions $f$ such that $f\circ\exp$ is a polynomial on $\fg$).
This is a dense subalgebra of $A_\M(G)$ (see \cite{Good_hreps}).
Using the identification $\P(G)\Tens\P(G)\cong\P(G\times G)$, one can show
that $\P(G)$ has a Hopf algebra structure given by \eqref{grp_Hopf}
(cf. \cite{Good_filt}, Prop.~2.1).
The algebra $U(\fg)$ acts on $\P(G)$ via left-invariant differential
operators, and this leads to
a canonical Hopf algebra pairing $U(\fg)\times\P(G)\to\CC$
defined by $\langle a,f\rangle = (af)(e)$ for $a\in U(\fg),\; f\in\P(G)$
(cf. \cite[XVI.3]{Hochbook}).
Goodman \cite{Good_hreps} proved that this pairing extends
to a pairing $U(\fg)_\M\times A_\M(G)\to\CC$ and defines a topological
isomorphism between $U(\fg)_\M$ and the strong dual space of $A_\M(G)$.
Since $U(\fg)_\M$ is a nuclear Fr\'echet space, it follows that
the multiplication on $A_\M(G)$ yields (by duality) a comultiplication
$U(\fg)_\M\to U(\fg)_\M\Ptens U(\fg)_\M$ that extends the comultiplication
of $U(\fg)$. Similarly, the multiplication on $U(\fg)_\M$ yields a comultiplication
$A_\M(G)\to A_\M(G)\Ptens A_\M(G)$ that extends the comultiplication
of $\P(G)$. It is also easy to see that the antipode and the counit of $U(\fg)$
(resp. $\P(G)$) extend by continuity to $U(\fg)_\M$
(resp. $A_\M(G)$), so that $U(\fg)_\M$ (resp. $A_\M(G)$) becomes a
Hopf $\Ptens$-algebra containing $U(\fg)$ (resp. $\P(G)$)
as a dense Hopf subalgebra. Thus $U(\fg)_\M$ and $A_\M(G)$ are
well-behaved Hopf $\Ptens$-algebras dual to each other.

The above properties of $U(\fg)_\M$ imply the following.

\begin{prop}
\label{prop:wght_cnv}
Let $\fg$ be a nilpotent Lie algebra with a positive filtration $\F$,
and let $\M$ be an entire $\F$-weight sequence. Then $U(\fg)_\M$
is a well-behaved power series envelope of $U(\fg)$.
\end{prop}

\begin{prop}
\label{prop:A(G)_contr}
$A_\M(G)$ is contractible.
\end{prop}
\begin{proof}
Given a function $f\colon G\to\CC$ and $z\in\CC$, define
$f_z\colon G\to\CC$ by $f_z(g)=f(\delta_z g)$. Using the obvious identity
$W_\M(\delta_z g)=W_\M(\delta_{|z|} g)$, we obtain
\begin{align}
N_r(f_z)=\sup_g\frac{|f(\delta_z g)|}{W_\M(\delta_r g)}
&=\sup_h\frac{|f(h)|}{W_\M(\delta_r\delta_z^{-1}h)}\notag\\
\label{N-transl}
&=\sup_h\frac{|f(h)|}{W_\M(\delta_{r|z|^{-1}}h)}
=N_{r|z|^{-1}}(f)
\end{align}
for each $r>0$ and each $z\ne 0$.
Therefore for each $f\in A_\M(G)$ we have $f_z\in A_\M(G)$, and the
mapping $A_\M(G)\to A_\M(G),\; f\mapsto f_z$
is continuous. Note also that $f_1=f$ and $f_0=f(e)1$ for each
$f\in A_\M(G)$.

For each $f\in\P(G)$, the function $(z,g)\mapsto f_z(g)$ is clearly a polynomial
on $\CC\times G$. Therefore we have an algebra homomorphism
$$
\varPhi_0\colon\P(G)\to\P(\CC,\P(G))\cong\P(\CC\times G),\quad
\varPhi_0(f)(z)=f_z.
$$
We use the same symbol $\varPhi_0$ to denote the composition
of the above homomorphism with the canonical embedding
$\P(\CC,\P(G))\hookrightarrow\O(\CC,A_\M(G))$.

We claim that $\varPhi_0$ is continuous w.r.t. the topology on $\P(G)$
inherited from $A_\M(G)$ and the compact-open topology on $\O(\CC,A_\M(G))$.
Indeed, let $\|\cdot\|$ be a continuous seminorm on $A_\M(G)$ and let $R>0$.
Then the rule
$$
\| u\|_R=\sup\{ \| u(z)\| : |z|\le R\}
$$
defines a continuous seminorm on $\O(\CC,A_\M(G))$. Furthermore,
the compact-open topology on $\O(\CC,A_\M(G))$ is generated by all seminorms
of this form. Therefore to prove the continuity of $\varPhi_0$ we have to show that
for each continuous seminorm $\|\cdot\|$ on $A_\M(G)$ and each $R>0$
the seminorm $f\mapsto\|\varPhi_0(f)\|_R$ is continuous on $\P(G)$.
Since $\|\cdot\|$ is continuous on $A_\M(G)$, we see that for each $r>0$ there exists
$C>0$ such that $\| f\|\le CN_{rR}(f)$ for all $f\in A_{\M,rR}(G)$.
Now let $f$ be in $\P(G)$. Using \eqref{N-transl}
and the fact that $N_r\le N_s$ whenever $s\le r$, we obtain
$$
\| \varPhi_0(f)\|_R=\sup_{|z|\le R} \| f_z\|
\le C\sup_{|z|\le R} N_{rR}(f_z)
=C\sup_{0<|z|\le R} N_{rR|z|^{-1}}(f)
=CN_r(f).
$$
This means that the seminorm $f\mapsto\|\varPhi_0(f)\|_R$ is continuous on $\P(G)$
w.r.t. the topology inherited from $A_\M(G)$. Therefore $\varPhi_0$ is continuous.
Since $\P(G)$ is dense in $A_\M(G)$ (see \cite{Good_hreps}), we see that $\varPhi_0$
extends to a continuous homomorphism
$$
\varPhi\colon A_\M(G)\to \O(\CC,A_\M(G))\cong \O(\CC)\Ptens A_\M(G).
$$
Let $\eps\colon A_\M(G)\to\CC,\; f\mapsto f(e)$ denote the counit of $A_\M(G)$.
We claim that $\varPhi$ is a homotopy between $\id_{A_\M(G)}$ and $\eta\eps$
(see Definition~\ref{def:hmt}). Indeed,
for each $z\in\CC$ the mappings $f\mapsto f_z$ and $f\mapsto\varPhi(f)(z)$
from $A_\M(G)$ to itself are continuous, and they coincide on $\P(G)$.
Hence $\varPhi(f)(z)=f_z$ for each $f\in A_\M(G)$ and each $z\in\CC$.
In particular, $\varPhi(f)(1)=f$ and $\varPhi(f)(0)=f(e)1=(\eta\eps)(f)$.
In other words,
$$
(\eps_1\otimes\id_{A_\M(G)})\,\varPhi=\id_{A_\M(G)}\quad\text{and}\quad
(\eps_0\otimes\id_{A_\M(G)})\,\varPhi=\eta\eps,
$$
where the augmentations $\eps_k\colon\O(\CC)\to\CC$ ($k=0,1$)
are defined by $\eps_k(f)=f(k)$.
Since $\O(\CC)$ is an exact algebra, we conclude that
$\id_{A_\M(G)}$ is homotopic to $\eta\eps$, and so
$A_\M(G)$ is contractible.
\end{proof}

Now Proposition \ref{prop:wght_cnv}, Proposition \ref{prop:A(G)_contr},
and Corollary \ref{cor:conv_loc} imply the following.

\begin{theorem}
\label{thm:weight_loc}
Let $\fg$ be a nilpotent Lie algebra with a positive filtration $\F$,
and let $\M$ be an entire $\F$-weight sequence. Then
$U(\fg)_\M$ is stably flat over $U(\fg)$.
\end{theorem}

\section{Algebras of analytic functionals\\
and hyperenveloping algebras}
\label{sect:an_func}

Let $\fg$ be a Lie algebra, and let $G$ denote the corresponding connected,
simply connected complex Lie group.
In this section we prove that the {\em hyperenveloping algebra}
$\fF(\fg)$ (see \cite{Rash})
is always stably flat over $U(\fg)$. We also show that
the {\em algebra of analytic functionals}
$\A(G)$ (see \cite{Litv}) is stably flat over $U(\fg)$ if and only if $\fg$
is solvable.

First recall some definitions. Let $G$ be a complex Lie group.
The Fr\'echet algebra $\O(G)$ of holomorphic functions on $G$ has a canonical
structure of Hopf $\Ptens$-algebra given by \eqref{grp_Hopf}.
Since $\O(G)$ is nuclear,
the strong dual space, $\O(G)'$, is a Hopf $\Ptens$-algebra and, in addition,
a nuclear (DF)-space. It is denoted by $\A(G)$ and is called {\em the algebra
of analytic functionals} on $G$ (see \cite{Litv}).
The product of $\alpha,\beta\in\A(G)$ is called the
{\em convolution} and is denoted by $\alpha*\beta$.
By definition, we have
$\langle \alpha * \beta,f\rangle=\langle \alpha\otimes\beta,\Delta f\rangle$
for each $\alpha,\;\beta\in\A(G)$ and each $f\in\O(G)$.

Consider the algebra $\O_e$ of germs
of holomorphic functions at the identity $e\in G$. We endow $\O_e$ with
its usual inductive limit topology, i.e., $\O_e=\varinjlim\O(U)$, where
$U$ runs through the collection of all neighborhoods of $e$.
Relative to this topology, $\O_e$ becomes a nuclear, complete (DF)-space
(see \cite{Groth}, Chap.~II, \S2, no.~3).
Moreover, the multiplication in $\O_e$ is jointly continuous, so that
$\O_e$ is a $\Ptens$-algebra.

By localizing \eqref{grp_Hopf} at the identity, we obtain a Hopf
$\Ptens$-algebra structure on $\O_e$ (cf. \cite[4.2]{Litv_hypergr}
and \cite[3.2.3]{VO}).
More exactly, take a neighborhood $U$ of $e$,
choose a neighborhood $V\ni e$ such that $V^2\subset U$, and consider the
map
$$
\Delta_{UV}\colon\O(U)\to\O(V\times V)\cong\O(V)\Ptens\O(V),\quad
(\Delta_{UV}f)(x,y)=f(xy).
$$
Composing with the restriction map $\O(V)\Ptens\O(V)\to\O_e\Ptens\O_e$
and taking the direct limit over $U\ni e$, we obtain a comultiplication
$\Delta\colon\O_e\to\O_e\Ptens\O_e$. The counit and the antipode
are defined similarly using \eqref{grp_Hopf}.
Since all Lie groups with the same Lie algebra are locally isomorphic,
the Hopf algebra structure on $\O_e$ depends only on $\fg$.

By definition, the {\em hyperenveloping algebra} $\fF(\fg)$ is the
strong dual algebra of $\O_e$. (Note that the original definition of
$\fF(\fg)$ given by Rashevskii in \cite{Rash} was different; we follow
the approach suggested by Litvinov \cite{Litv,Litv_hypergr}.)
Since $\O_e$ is a nuclear (DF)-space, $\fF(\fg)$ is a nuclear Fr\'echet space.

Let $\m_e$ be the ideal of $\O_e$ consisting of all germs vanishing at $e$.
Consider the formal completion $\wh{\O}_e=\varprojlim\O_e/\m_e^n$.
We endow each quotient $\O_e/\m_e^n$ with the
standard topology of a finite-dimensional vector space, so that
$\wh{\O}_e$ becomes a nuclear Fr\'echet algebra. Moreover,
the comultiplication and the antipode of $\O_e$ extend to
$\wh{\O}_e$ (cf. \cite[3.2.3]{VO}), so that $\wh{\O}_e$
has a canonical structure of Hopf $\Ptens$-algebra.

There is a natural Hopf algebra pairing
between $U(\fg)$ and $\wh{\O}_e$ defined as follows (for details,
see \cite[3.2]{VO}).
For each $X\in\fg$, let $\wt{X}$ denote the corresponding left-invariant
vector field on $G$. For each open set $U\subset G$ we use the same symbol
$\wt{X}$ to denote the corresponding derivation of $\O(U)$.
Taking the direct limit over $U\ni e$, we see that $\wt{X}$ determines
a derivation of $\O_e$ which we also denote by $\wt{X}$.
It is easy to see that $\wt{X}(\m_e^n)\subset \m_e^{n-1}$ for each $n$,
so that $\wt{X}$ extends to a derivation of $\wh{\O}_e$ (again denoted by $\wt{X}$).
The resulting map $\fg\to\Der\wh{\O}_e,\; X\mapsto\wt{X}$, yields
an algebra homomorphism $\rho\colon U(\fg)\to\End_{\CC}\wh{\O}_e$. Thus
$U(\fg)$ acts on $\wh{\O}_e$ via ``formal left-invariant differential
operators'' (cf. Section~\ref{sect:weight}).
The canonical pairing between $U(\fg)$ and $\wh{\O}_e$ defined by
$\la a,\f\ra = [\rho(a)\f](e)$ for each $a\in U(\fg),\; \f\in\wh{\O}_e$,
gives an algebraic isomorphism between $\wh{\O}_e$ and the algebraic
dual of $U(\fg)$ \cite[3.2.3]{VO}. If we endow $U(\fg)$ with the finest
locally convex topology, then $\wh{\O}_e$ becomes the topological dual of $U(\fg)$,
and the strong dual topology on $\wh{\O}_e$ coincides with the inverse limit
topology introduced above (cf. the beginning of Section~\ref{sect:loc_dual}).

The restriction maps
$$
\O(G) \to \O_e \to \wh{\O}_e
$$
are obviously Hopf $\Ptens$-algebra homomorphisms.
Taking the dual maps, we obtain Hopf $\Ptens$-algebra homomorphisms
\begin{equation}
\label{chain_emb_1}
U(\fg) \xra{\lambda} \fF(\fg) \to \A(G).
\end{equation}
Note that $\O_e\to\wh{\O}_e$ is always injective with dense range,
so that ${U(\fg)\to\fF(\fg)}$ has the same property.
The restriction map $\O(G)\to\O_e$ is injective provided $G$ is connected,
and has dense range provided $G$ is a Stein group.
Therefore for each connected Stein group (in particular, for each connected,
simply connected complex Lie group) both the maps in \eqref{chain_emb_1}
are injective with dense ranges (cf. \cite{Litv}).

Let $\tau\colon U(\fg)\to\A(G)$ denote the composition
of the above maps. It follows from the definition of the duality between $U(\fg)$
and $\wh{\O}_e$ that $\la \tau(X),f\ra=(\wt{X}f)(e)$
for all $X\in\fg,\; f\in\O(G)$.
It is also easy to see that for each $X\in\fg$ the action of $X$
on $\A(G)'=\O(G)$ determined by $\tau$
(see Subsection \ref{subsect:parall}) coincides with the derivation $\wt{X}$. Indeed,
given $x\in G$, denote by $\delta_x\in\A(G)$ the functional which is evaluation
at $x$. Then for each $X\in\fg$, $f\in\O(G)$, and $x\in G$
we have
\begin{multline*}
(X\cdot_\tau f)(x)=\langle X\cdot_\tau f,\delta_x\rangle
=\langle f,\delta_x \tau(X)\rangle\\
=\langle \Delta f,\delta_x\otimes\tau(X)\rangle
=\frac{d}{dt}\Bigl |_{t=0} f(x\exp tX) = (\tilde Xf)(x),
\end{multline*}
i.e., $X\cdot_\tau f=\wt{X} f$, as required.

Similarly, for each $X\in\fg$ the action of $X$
on $\fF(\fg)'=\O_e$ determined by the canonical homomorphism
$\lambda\colon U(\fg)\to\fF(\fg)$
coincides with $\wt{X}$. Indeed, given $f\in\O_e$,
denote by $\hat f$ the canonical image of $f$ in $\wh{\O}_e$;
then for each $X\in\fg$ and each $a\in U(\fg)$ we have
\begin{multline*}
\langle X\cdot_\lambda f,\lambda(a)\rangle = \langle f,\lambda(aX)\rangle =
\la \hat f, aX\ra = [\rho(aX)\hat f](e)\\
= [\rho(a)\wt{X}\hat f](e) = \langle (\wt{X}f)\sphat\, ,a\rangle
=\la \wt{X}f,\lambda(a)\ra.
\end{multline*}
Since $\Im\lambda$ is dense in $\fF(\fg)$, this implies
$X\cdot_\lambda f=\wt{X} f$, as required.

\begin{prop}
\label{prop:O(G)_par}
Let $G$ be a Stein group with Lie algebra $\fg$.
Then $\O(G)$ is $\fg$-parallelizable.
\end{prop}
\begin{proof}
By Lemma \ref{lemma:Stein_diff}, we may identify the $\O(G)$-module
$\Omega^1(\O(G))$ of K\"ahler differentials with the module $\Omega^1(G)$
of holomorphic $1$-forms on $G$ in such a way that the exterior (de Rham)
derivative $d\colon\O(G)\to \Omega^1(G)$ becomes a universal derivation.
Denote by $\Vect(G)$ the Lie algebra of
holomorphic vector fields on $G$. In what follows, we identify $\fg$ with the
Lie subalgebra of $\Vect(G)$ consisting of left-invariant vector fields.
Each $\omega\in\Omega^1(G)$
can be viewed as an $\O(G)$-module morphism $\Vect(G)\to\O(G)$.
Hence the rule $\omega\mapsto\omega|_{\fg}$ determines a linear map
$\varphi\colon\Omega^1(G)\to C^1(\fg,\O(G))$ which is easily seen to
be an $\O(G)$-module morphism.
Evidently, $\varphi(df)(\wt{X})=\wt{X}f$ for each $X\in\fg$,
$f\in\O(G)$, i.e., $\varphi d=d^0$. It remains to show that $\varphi$
is an isomorphism.

Given $\omega\in\fg^*$, denote by $\tilde\omega\in\Omega^1(G)$ the corresponding
left-invariant $1$-form on $G$.
Let $\psi\colon C^1(\fg,\O(G))\to\Omega^1(G)$ be the unique $\O(G)$-module
morphism taking each $\omega\otimes 1\in C^1(\fg,\O(G))$ to $\tilde\omega$.
Recall that for each $\omega\in\fg^*$ the value of the
left-invariant form $\tilde\omega$ at a left-invariant vector field
$\wt{X}$ is the constant function equal to $\langle \omega, X\rangle$
(see, e.g., \cite[3.12]{Warner}). This means precisely that
$\tilde\omega|_{\fg}=\omega\otimes 1$, i.e., $\varphi\psi=\id_{C^1(\fg,\O(G))}$.

Let $\omega_1,\ldots ,\omega_n$ be a basis of $\fg^*$. Then
$\tilde\omega_1(x),\ldots ,\tilde\omega_n(x)$ is clearly a basis of
the cotangent space, $T_x^* G$, for each $x\in G$. Hence the forms
$\tilde\omega_1,\ldots ,\tilde\omega_n$ generate $\Omega^1(G)$ as an
$\O(G)$-module. This implies, in particular, that $\psi$ is surjective. Since
$\varphi\psi=\id$, we see that $\varphi$ and
$\psi$ are inverse to one another. This completes the proof.
\end{proof}

Combining this with Propositions~\ref{prop:DR_Stein} and \ref{prop:par_crit},
we obtain the following well-known fact.

\begin{corollary}
\label{cor:top-Lie-Stein}
Let $G$ be a Stein group with Lie algebra $\fg$. Then
$H^p(\fg,\O(G))\cong H^p_{\tp}(G,\CC)$ for each $p$.
\end{corollary}

The following result is an analytic version of \cite[Prop. 7.2]{T2}.

\begin{theorem}
\label{thm:A(G)_loc}
Let $\fg$ be a Lie algebra, and let
$G$ be the corresponding connected, simply connected complex Lie group.
Then $\A(G)$ is stably flat over $U(\fg)$
if and only if $\fg$ is solvable.
\end{theorem}
\begin{proof}
Suppose that $\fg$ is solvable. Then
$G$ is also solvable and so is biholomorphic
with $\CC^n$ (see, e.g., \cite{Bour_Lie-II}, Chap.~III, \S9, no.~6).
Hence $\O(G)\cong\O(\CC^n)$ is
a contractible $\Ptens$-algebra (cf. Subsection~\ref{subsect:hmt}).
On the other hand, $\O(G)$ is $\fg$-parallelizable by
Proposition~\ref{prop:O(G)_par}. Now it remains to apply
Theorem~\ref{thm:par_contr_loc}.

Now suppose that $\fg$ is not solvable. Consider the Levi decomposition
$\fg=\mathfrak r\oplus\mathfrak l$ ($\mathfrak r=\mathrm{rad}\,\fg$,
$\mathfrak l$ is a semisimple
subalgebra, $\mathfrak l\ne 0$). Let $R$ and $L$ be the corresponding analytic
subgroups of $G$. Since $G$ is simply connected, it is isomorphic to the
semidirect product $R\rtimes L$ (see, e.g., \cite{VO-sem}, Chap.~6).
Since $L$ is semisimple, we have $H^m_{\tp}(L,\CC)\ne 0$ for
$m=\dim L$ (see, e.g., \cite{Matsushima}, Lemme~5).
Therefore $H^m_{\tp}(G,\CC)\ne 0$. On the other hand,
$H^m_{\tp}(G,\CC)\cong H^m(\fg,\O(G))$ by
Corollary~\ref{cor:top-Lie-Stein}, and so
the augmented standard complex
$$
0\to\CC\to C^\cdot(\fg,\O(G))
$$
is not exact. Using the reflexivity argument
(cf. the proof of Theorem~\ref{thm:par_contr_loc}), we see that the dual complex
$$
0 \lar \CC \lar C_\cdot(\fg,\A(G))
$$
does not split in $\mathbf{LCS}$.
Therefore $U(\fg)\to \A(G)$ is not a localization, i.e., $\A(G)$ is not stably
flat over $U(\fg)$.
\end{proof}

We now turn to the hyperenveloping algebra $\fF(\fg)$.
Recall that a commutative algebra is called {\em local} if it has
a unique maximal ideal. By a {\em local $\Ptens$-algebra} we mean
a commutative $\Ptens$-algebra $A$ which is local in the above sense
and such that the maximal ideal of $A$ is closed and has codimension $1$.
For example, $\O_e$ is a local $\Ptens$-algebra with maximal ideal
$\m_e=\{ f\in\O_e : f(e)=0\}$ (see above), and the same is true for $\wh{\O}_e$.

We need the following simple lemma.
\begin{lemma}
\label{lemma:simple}
Let $A$ be an algebra, $I\subset A$ a left ideal,
and $E$ a finite-dimensional vector space. Then for each
$T\in\Hom_{\CC}(E,A)$ the following conditions are equivalent:
\begin{itemize}
\item[(i)] $T\in I\cdot\Hom_{\CC}(E,A)$;
\item[(ii)] $\Im T\subset I$.
\end{itemize}
\end{lemma}
\begin{proof}
The implication (i)$\Longrightarrow$(ii) is clear.
To prove the converse, take a basis $(e_i)$ of $E$, and let
$(e^i)$ be the dual basis of $E^*$. Identifying $\Hom_{\CC}(E,A)$ and
$E^*\Tens A$, we see that $T=\sum_i e^i\otimes a_i$, where $a_i=T(e_i)\in I$.
Setting $T_i=e^i\otimes 1$, we obtain $T=\sum_i a_i T_i\in I\cdot\Hom_{\CC}(E,A)$,
as required.
\end{proof}

\begin{lemma}
\label{lemma:loc_act}
Let $A$ be a local $\Ptens$-algebra with maximal ideal $\m$,
and let $\fg$ be a Lie algebra
acting on $A$ by derivations. Suppose there exists a linear map
$\chi\colon\fg^*\to A$ such that
\begin{itemize}
\item[(i)] $\Im\chi$ generates a dense subalgebra of $A$;
\item[(ii)] $X\cdot \chi(\omega)=\langle \omega,X\rangle 1\mod\m$
for each $X\in\fg$ and each $\omega\in\fg^*$.
\end{itemize}
Then $A$ is $\fg$-parallelizable.
\end{lemma}
\begin{proof}
We proceed in much the same way as in the proof of Theorem \ref{thm:conv_prl}.
Consider the $A$-module morphism $\varphi\colon A\Tens\fg^*\to C^1(\fg,A)$ uniquely
determined by $1\otimes\omega\mapsto d^0(\chi(\omega))$. Our objective
is to prove that $\varphi$ is an isomorphism. To this end, note that,
since $A$ is local and both $A\Tens\fg^*$ and
$C^1(\fg,A)$ are free and finitely generated, we need only prove that
the induced map
\begin{equation}
\label{barphi}
\bar\varphi\colon
A\Tens\fg^*/\m\cdot A\Tens\fg^*\to C^1(\fg,A)/\m\cdot C^1(\fg,A)
\end{equation}
is a vector space isomorphism (see \cite{Bour_AC}, Chap.~II, \S 3, no.~2).

Since $\m$ is closed and has codimension $1$, there exists a continuous
homomorphism $\eps\colon A\to\CC$ such that $\m=\Ker\eps$.
Hence we can identify $A\Tens\fg^*/\m\cdot A\Tens\fg^*$ and $\fg^*$
via the map
\begin{equation}
\label{g^*-1}
\alpha\colon\fg^*\to A\Tens\fg^*/\m\cdot A\Tens\fg^*,\quad
\omega\mapsto 1\otimes\omega+\m\cdot A\Tens\fg^*.
\end{equation}
The inverse map is given by $a\otimes\omega+\m\cdot A\Tens\fg^*\mapsto
\eps(a)\omega$.

Next consider the map
\begin{equation}
\label{g^*-2}
\beta\colon C^1(\fg,A)/\m\cdot C^1(\fg,A)\to\fg^*,\quad
T+\m\cdot C^1(\fg,A)\mapsto \eps T.
\end{equation}
Lemma \ref{lemma:simple} implies that $\beta$ is well defined and bijective.
Indeed, the map
taking each $\omega\in\fg^*$ to
$\omega\otimes 1+\m\cdot C^1(\fg,A)$ is readily seen to be an inverse
of $\beta$.

Now it is easy to see
that the map $\bar\varphi$
defined by \eqref{barphi}
corresponds to the identity mapping of $\fg^*$
under the identifications \eqref{g^*-1} and \eqref{g^*-2}.
Indeed, for each $\omega\in\fg^*$ we have $(\bar\varphi\alpha)(\omega)=
d^0(\chi(\omega))+\m\cdot C^1(\fg,A)$, and hence
$$
\langle (\beta\bar\varphi\alpha)(\omega),X \rangle =
\eps(d^0(\chi(\omega))X)=\eps(X\cdot \chi(\omega))=
\eps(\langle\omega,X\rangle 1)=\langle\omega,X\rangle
$$
for every $X\in\fg$. Therefore $\beta\bar\varphi\alpha=\id_{\fg^*}$, and so
$\bar\varphi$ is an isomorphism.
By the above remarks, so is $\varphi$.

Now define a derivation $d\colon A\to A\Tens\fg^*$ by $d=\varphi^{-1}d^0$.
Note that $d(\chi(\omega))=1\otimes\omega$ for each $\omega\in\fg^*$.
Choose a basis $(e_i)$ of $\fg$, and let $(e^i)$ be the dual basis
of $\fg^*$. We may identify $A\Tens\fg^*$ and $A^n$ ($n=\dim\fg$)
via the map
$$
\psi\colon (a_1,\ldots ,a_n)\in A^n \mapsto \sum_i a_i\otimes e^i\in A\Tens\fg^*.
$$
Let $\partial=\psi^{-1} d = (\partial_1,\ldots ,\partial_n)\colon A \to A^n$
be the derivation corresponding to $d$ under the above identification.
Since $d(\chi(e^j))=1\otimes e^j$ for each $j$, it follows that
$\partial_i (\chi(e^j))=\delta_{ij}$ for each $i,j$.
Thus the elements $x_i=\chi(e^i)\in A$ and the derivations $\partial_i\in\Der A$
satisfy the conditions of Lemma \ref{lemma:Omega_free}.
Therefore $\partial$ is a universal derivation.
Since both $\varphi$ and $\psi$ are isomorphisms, we conclude that
$d^0=\varphi d=\varphi\psi\partial$ is a universal derivation as well,
i.e., $A$ is $\fg$-parallelizable.
\end{proof}

\begin{theorem}
\label{thm:hyper_loc}
Let $\fg$ be a Lie algebra. Then
$\fF(\fg)$ is stably flat over $U(\fg)$.
\end{theorem}
\begin{proof}
In view of Theorem \ref{thm:par_contr_loc}, it suffices to check that
$\O_e$ is contractible and $\fg$-parallelizable.

To prove the contractibility of $\O_e$, it suffices to do this for the algebra
$\O_0$ of holomorphic germs at the origin $0\in\CC^n$.
For each $r>0$ denote by $D^n_r$ the polydisc in $\CC^n$ of radius $r$,
i.e.,
$$
D^n_r=\{ z=(z_1,\ldots ,z_n)\in\CC^n : |z_i|<r\;\forall i=1,\ldots ,n\}.
$$
Consider the homomorphism
\begin{gather*}
\varPhi_r\colon\O(D^n_r)\to\O(D^1_2)\Ptens\O(D^n_{r/2})\cong
\O(D^1_2\times D^n_{r/2}),\\
(\varPhi_r f)(z,w)=f(zw).
\end{gather*}
Let $\eps_k\colon\O(D^1_2)\to\CC\; (k=0,1)$ be given by $\eps_k(f)=f(k)$.
We clearly have
\begin{gather}
\bigl[(\eps_0\otimes\id)(\varPhi_r f)\bigr](w)=f(0),\notag\\
\label{contr_O_0}
\bigl[(\eps_1\otimes\id)(\varPhi_r f)\bigr](w)=f(w).
\end{gather}
Composing $\varPhi_r$ with the restriction map
$\O(D^1_2)\Ptens\O(D^n_{r/2})\to \O(D^1_2)\Ptens\O_0$ and taking next
the inductive limit over $D_r\ni 0$, we obtain a homomorphism
$\varPhi\colon\O_0\to\O(D^1_2)\Ptens\O_0$. Relations
\eqref{contr_O_0} imply that $(\eps_0\otimes\id_{\O_0})\,\varPhi=\eta\eps_{\O_0}$
and $(\eps_1\otimes\id_{\O_0})\,\varPhi=\id_{\O_0}$, where
$\eps_{\O_0}\colon\O_0\to\CC$ is the evaluation at $0$.
Since $\O(D^1_2)$ is an exact algebra, we see that
$\varPhi$ is a homotopy between $\id_{\O_0}$ and
$\eta\eps_{\O_0}$, and so $\O_0$ is contractible.

We now turn to the $\fg$-parallelizability of $\O_e$.
Since the exponential
map is biholomorphic in a neighborhood of $0\in\fg$,
it follows that for each $\omega\in\fg^*$
there exists a unique function $f_\omega$ holomorphic
in a neighborhood of $e\in G$ such that $f_\omega(\exp\xi)=\omega(\xi)$
for all sufficiently small $\xi\in\fg$. Consider the map $\chi\colon\fg^*\to\O_e$
taking each $\omega\in\fg^*$ to the germ of $f_\omega$ at $e$.

We claim that $\chi$ satisfies the conditions of Lemma \ref{lemma:loc_act}.
To prove this,
fix a basis $(e_i)$ of $\fg$, and
let $x^1,\ldots ,x^n$ be the corresponding canonical coordinates of the first kind
on a suitable neighborhood of $e\in G$.
Recall that they are defined by the rule
$x^j(\exp\sum_i t^i e_i)=t^j$ for all $j=1,\ldots ,n$.
For each $\omega\in\fg^*$ we have
$$
f_\omega\bigl(\exp\sum_i t^i e_i\bigr)=
\sum_i \omega(e_i) t^i,\quad\mbox{i.e.,}\quad
f_\omega=\sum_i \omega(e_i) x^i.
$$
Therefore $\Im\chi$ consists of all germs of linear
functions in $x^1,\ldots ,x^n$. Since polynomials in $x^1,\ldots ,x^n$ form
a dense subalgebra of $\O_e$, we conclude that condition~(i)
of Lemma~\ref{lemma:loc_act} is satisfied.

Finally, for each $X\in\fg$ and each $\omega\in\fg^*$ we have
$$
(X\cdot\chi(\omega))(e)=(\wt{X} f_\omega)(e)=
\frac{d}{dt}\Bigl |_{t=0} f_\omega(\exp tX) = \omega(X).
$$
Thus we see that condition (ii) of Lemma \ref{lemma:loc_act} is also satisfied.
Hence $\O_e$ is $\fg$-parallelizable.
Now the result follows from Theorem \ref{thm:par_contr_loc}.
\end{proof}

\begin{remark}
A similar argument applied to the local algebra $A=\wh{\O}_e$ shows
that $A$ is contractible and $\fg$-parallelizable. Hence the standard cochain
complex $0\to \CC \to C^\cdot (\fg,A)$ splits in $\mathbf{LCS}$.
Taking the topological dual, we recover the classical fact that the Koszul
complex $0 \lar \CC \lar V_\cdot(\fg)$ is exact.
\end{remark}

\section{Relations between various completions of $U(\fg)$}

In this final section, we explain how the completions of $U(\fg)$
considered above are related to one another, and formulate some
open problems.

Let $\fg$ be a nilpotent Lie algebra, and let $G$ be the
corresponding connected, simply connected complex Lie group.
Choose a positive filtration $\F$ on $\fg$,
and let $\M$ be an entire $\F$-weight sequence (see Section~\ref{sect:weight}).
The algebra $A_\M(G)$ is a subalgebra
of $\O(G)$, and it is easy to see that the inclusion map $A_\M(G)\to\O(G)$
is continuous. Indeed, given $r>0$ and a compact set $K\subset G$, let
$C=\sup_{g\in K} W_\M(\delta_r g)$.
Then for each $f\in A_{\M,r}(G)$ we have
$$
\| f\|_K=\sup_{g\in K} |f(g)| \le C \sup_{g\in K}\frac{|f(g)|}{W_\M(\delta_r g)}
\le CN_r(f).
$$
Hence the inclusion of $A_\M(G)$ into $\O(G)$ is continuous.

Consider the Hopf algebra $\P(G)$ of polynomial functions on $G$. This is
a Hopf $\Ptens$-algebra w.r.t. the finest locally convex topology.
We have a chain of canonical inclusion/restriction maps
\begin{equation}
\label{chain_emb_func}
\P(G) \to A_\M(G) \to \O(G) \to \O_e \to \wh{\O}_e
\end{equation}
which are obviously Hopf $\Ptens$-algebra homomorphisms.
To examine the dual of this chain,
choose an $\F$-basis $e_1,\ldots ,e_N$ of $\fg$, and
recall that there is a duality $\la\cdot ,\cdot\ra_\varkappa$
between the formal completion
$[U(\fg)]$ of $U(\fg)$ and the polynomial algebra $\CC[z_1,\ldots ,z_N]$
defined by $\la a,\varphi\ra_\varkappa=\varkappa^{-1}(\varphi)(a)$
for all $a\in [U(\fg)],\; \varphi\in\CC[z_1,\ldots ,z_N]$
(see Lemma~\ref{lemma:dual_[U]}). Thus the dual of the inclusion map
$i\colon U(\fg)_\M\to [U(\fg)]$ can be viewed as a $\Ptens$-algebra homomorphism
from $\CC[z_1,\ldots ,z_N]$ to $A_\M(G)$.
For each $u\in U(\fg)_\M$ and each $\varphi\in\CC[z_1,\ldots ,z_N]$
we have
$$
\la u,i'(\varphi)\ra=\la i(u),\varphi\ra_\varkappa,
$$
where the brackets $\la\cdot ,\cdot\ra$ on the left-hand side
denote the duality between $U(\fg)_\M$ and $A_\M(G)$.

We claim that $i'\colon\CC[z_1,\ldots ,z_N]\to A_\M(G)$
becomes the canonical inclusion of $\P(G)$ into $A_\M(G)$ if we identify
$\CC[z_1,\ldots ,z_N]$ with $\P(G)$ using the canonical coordinates of
the second kind. Indeed, for each $\alpha\in\Z_+^N$ and each $\psi\in\O(G)$
we have (in the standard multi-index notation)
\begin{equation}
\label{e^alpha}
\la e^\alpha,\psi\ra=
[e^\alpha\psi](e)=D^\alpha_z\psi(0)
\end{equation}
w.r.t. the canonical coordinates of the second kind (see \cite{Litv}, Lemma~7.2).
On the other hand,
it follows from \eqref{dual_[U]} that
$\la e^\alpha,\varphi\ra_\varkappa=
D^\alpha_z \varphi(0)$
for each $\varphi\in\CC[z_1,\ldots ,z_N]=\P(G)$.
Together with \eqref{e^alpha}, this gives
$\la e^\alpha,\varphi\ra_\varkappa=\la e^\alpha, \varphi\ra$, and so
$
\la e^\alpha, i'(\varphi)\ra = \la i(e^\alpha),\varphi\ra_\varkappa =
\la e^\alpha, \varphi\ra_\varkappa
=\la e^\alpha, \varphi\ra
$
for each $\alpha\in\Z_+^N$. This implies that
$i'(\varphi)=\varphi$ for each $\varphi\in\P(G)$, which proves the claim.

Thus the sequence dual to \eqref{chain_emb_func} has the form
\begin{equation}
\label{chain_emb_2}
U(\fg) \to \fF(\fg) \to \A(G) \to U(\fg)_\M \xra{i} [U(\fg)].
\end{equation}
All the maps here are injective Hopf $\Ptens$-algebra
homomorphisms with dense ranges. Combining
Theorems~\ref{thm:hyper_loc}, \ref{thm:A(G)_loc}, \ref{thm:weight_loc},
Corollary~\ref{cor:[U]_loc}, and Proposition~\ref{prop:loc_compos}, we see that
all the morphisms in \eqref{chain_emb_2} are localizations.

\begin{prop}
Let $G$ be a connected, simply connected complex Lie group with Lie algebra
$\fg$, and let $\tau\colon U(\fg)\to\A(G)$ be the canonical homomorphism
(see \eqref{chain_emb_1}). Then there exists a unique Hopf
$\Ptens$-algebra homomorphism $j\colon \A(G)\to\wh{U}(\fg)$ such that
$\iota_{U(\fg)}=j\circ\tau$. (In other words, the canonical morphism
$\iota_{U(\fg)}\colon U(\fg)\to\wh{U}(\fg)$ factors through $\A(G)$).
\end{prop}
\begin{proof}
Given a Banach algebra $A$ and a homomorphism $\varphi\colon U(\fg)\to A$,
we shall construct a continuous homomorphism $\bar\varphi\colon\A(G)\to A$
such that the diagram
\begin{equation}
\label{A-factor}
\xymatrix{
\A(G) \ar[r]^{\bar\varphi} & A\\
U(\fg) \ar[u]^\tau \ar[ur]_{\varphi}
}
\end{equation}
is commutative.

For each $a\in A$ denote by $L_a\colon A\to A$ the left multiplication by $a$.
Consider the representation
$$
\pi_0\colon \fg\to\L(A),\quad \pi_0(X)=L_{\varphi(X)}.
$$
Since $G$ is simply connected, $\pi_0$ determines a holomorphic representation
$\pi\colon G\to GL(A)$ such that $\exp\circ\pi_0=\pi\circ\exp$
(see, e.g., \cite{Bour_Lie-II}, Chap.~III, \S6, no.~1).
For each $x\in A$ and each $y\in A'$,
let $\pi_{x,y}\in\O(G)$ denote the corresponding
matrix element of $\pi$ defined by $\pi_{x,y}(g)=\langle y,\pi(g)x\rangle$.
By \cite[Prop.~3.5]{Litv}, $\pi$ uniquely extends to a continuous representation
$\bar\pi\colon\A(G)\to\L(A)$ such that
$\langle y,\bar\pi(a')x\rangle = \langle a',\pi_{x,y}\rangle$ for all
$a'\in\A(G),\; x\in A,\; y\in A'$.

Consider the map $\eps_1\colon\L(A)\to A,\;\eps_1(T)=T(1)$,
and define $\bar\varphi\colon\A(G)\to A$ by $\bar\varphi=\eps_1\bar\pi$.
We claim that $\bar\varphi$ makes diagram \eqref{A-factor} commutative.
Indeed, for each $X\in\fg$ and each $y\in A'$ we have
\begin{multline*}
\langle y,\bar\varphi\tau(X)\rangle
=\langle y,\bar\pi(\tau(X))1\rangle
=\langle \tau(X),\pi_{1,y}\rangle
=(\wt{X}\pi_{1,y})(e)\\
=\frac{d}{dt}\Bigl |_{t=0} \pi_{1,y}(\exp tX)
=\frac{d}{dt}\Bigl |_{t=0} \langle y,\pi(\exp tX)1\rangle
=\frac{d}{dt}\Bigl |_{t=0} \langle y,\exp\pi_0(tX)1\rangle\\
=\frac{d}{dt}\Bigl |_{t=0} \langle y,\exp t\varphi(X)\rangle
=\langle y,\varphi(X)\rangle,
\end{multline*}
i.e., $\bar\varphi\tau=\varphi$. Hence diagram \eqref{A-factor} is commutative.
Since $\Im\tau$ is dense in $\A(G)$, we conclude that $\bar\varphi$
is an algebra homomorphism. For the same reason, $\bar\varphi$ is a unique
linear continuous map making \eqref{A-factor} commutative.

The above construction can easily be extended to the case where $A$ is
an Arens-Michael algebra. Indeed, we have $A=\varprojlim\{ A_\nu,\sigma^\mu_\nu\}$
for some inverse system $\{ A_\nu,\sigma^\mu_\nu\}$ of Banach algebras.
For each $\nu$, let $\sigma_\nu\colon A\to A_\nu$ denote the canonical map.
Given a homomorphism $\varphi\colon U(\fg)\to A$, we can extend the
homomorphism $\sigma_\nu\varphi\colon U(\fg)\to A_\nu$ to a
homomorphism $\bar\varphi_\nu\colon\A(G)\to A_\nu$
satisfying $\bar\varphi_\nu\tau=\sigma_\nu\varphi$.
Since such an extension is unique, we clearly have
$\sigma^\mu_\nu\bar\varphi_\mu=\bar\varphi_\nu$ whenever $\mu\succ\nu$.
Setting $\bar\varphi=\varprojlim\bar\varphi_\nu$, we obtain a $\Ptens$-algebra
homomorphism making \eqref{A-factor} commutative.

Now set $A=\wh{U}(\fg)$ and $\varphi=\iota_{U(\fg)}\colon U(\fg)\to\wh{U}(\fg)$.
Then the above construction yields a unique $\Ptens$-algebra homomorphism
$j=\bar\varphi\colon\A(G)\to\wh{U}(\fg)$ satisfying $j\tau=\iota_{U(\fg)}$.
Since $\Im\tau$ is dense in $\A(G)$, and since $\iota_{U(\fg)}$ is a Hopf
$\Ptens$-algebra homomorphism, we conclude that so is $j$. This completes the proof.
\end{proof}

The above theorem implies that for each nilpotent Lie algebra $\fg$
the chain of inclusions \eqref{chain_emb_2} can be completed as follows:
\begin{equation}
\label{chain_emb_3}
\xymatrix@R-20pt{
& & & U(\fg)_\M \ar[dr] \\
U(\fg) \ar[r] & \fF(\fg) \ar[r] & \A(G) \ar[ur] \ar[dr]_j && [U(\fg)]\\
& & & \wh{U}(\fg) \ar[ur]_{\wh{\theta}}
}
\end{equation}

The following summarizes the main results of the previous sections.

\begin{theorem}
Suppose $\fg$ is a positively graded Lie algebra,
$G$ is the corresponding connected, simply connected complex Lie group,
and $\M$ is an entire weight sequence on $\fg$.
Then all the arrows in \eqref{chain_emb_3}
are Hopf $\Ptens$-algebra localizations.
\end{theorem}

We end this section with some open problems.

\medskip\noindent
{\bfseries Problem 1.}
Is the canonical map $U(\fg)\to\wh{U}(\fg)$ a localization
for every nilpotent Lie algebra $\fg$?

\medskip
By Proposition \ref{prop:loc_compos}, we can replace $U(\fg)$
in the above problem by either $\fF(\fg)$ or $\A(G)$ (assuming that
$G$ is connected and simply connected).

\medskip\noindent
{\bfseries Problem 2.}
Let $\fg$ be a nilpotent Lie algebra.
\begin{itemize}
\item[1)] Is the canonical homomorphism $\wh{\theta}\colon\wh{U}(\fg)\to [U(\fg)]$
injective?
\item[2)] Is the algebra $\wh{U}'(\fg)$ contractible?
\end{itemize}

A positive answer to Problem 2 would imply a positive solution of Problem~1
(see Corollary~\ref{cor:conv_loc}).

\begin{remark}
The diagram dual to \eqref{chain_emb_3} has the form
$$
\xymatrix@R-20pt{
& & & A_\M(G) \ar[dl] \\
\wh{\O}_e & \O_e \ar[l] & \O(G) \ar[l] && \P(G) \ar[ul] \ar[dl]^{\wh{\theta}'} \\
& & & \wh{U}'(\fg) \ar[ul]^{j'}
}
$$
Recall that all the maps here (except for $\wh{\theta}'$ and $j'$) are the usual
set-theoretic inclusions/restrictions of function algebras.
Since $j$ and $\wh{\theta}$
have dense ranges, it follows that $\wh{\theta}'$ and $j'$ are injective.
Hence the algebra $\wh{U}'(\fg)$ can be viewed as a certain algebra
of holomorphic functions on $G$ containing the polynomials.
Thus Question~1) of Problem~2 has a positive solution if and only if
the polynomials are dense in $\wh{U}'(\fg)$.
\end{remark}

\vspace*{20mm}
\begin{flushleft}
\scshape\small
Department of Differential Equations and Functional Analysis\\
Faculty of Science\\
Peoples' Friendship University of Russia\\
Mikluho-Maklaya 6\\
117198 Moscow\\
RUSSIA

\medskip
{\itshape Address for correspondence:}\\

\medskip\upshape
Krupskoi 8--3--89\\
Moscow 119311\\
Russia

\medskip
{\itshape E-mail:} {\ttfamily pirkosha@sci.pfu.edu.ru, pirkosha@online.ru}
\end{flushleft}
\end{document}